\documentclass{article}

\usepackage{pinlabel}
\usepackage[UKenglish]{babel}

\usepackage{amssymb}
\usepackage{amsfonts}
\usepackage{amsmath}
\usepackage{amsthm}
\usepackage{enumerate}

\usepackage{hyperref}

\usepackage[all]{xy}
\usepackage{graphicx}
\usepackage{geometry}

\usepackage{tikz}
\usetikzlibrary{shapes,shadows,arrows}

\makeatletter
\renewcommand{\@biblabel}[1]{\hspace{-0.8em}
        \begin{tabular*}{0.1\textwidth}{l}
        [#1] \\
\end{tabular*}}
\makeatother

\title{$G_2-$instantons over asymptotically cylindrical manifolds}

\author{Henrique N. S\'a Earp}
\date{\today}

\newcommand{\tsum}{\tilde\#_S}

\newcommand{\ii}{^{(i)}}
\newcommand{\func}[1]{\textnormal {#1\,}}
\newcommand{\tublim}{\overset{C^\infty}{\underset{S\rightarrow\infty}{\longrightarrow}}}

\newcommand{\R}{\mathbb{R}}
\newcommand{\wleq}{\overset{w}{\leq}}
\newcommand{\bl}{\bar\lambda}
\newcommand{\vol}{\func{Vol}}
\newcommand{\dvol}{\:d\func{Vol}_\omega}
\newcommand{\dvtub}{\:d\text{V}_{tub}}

\newcommand{\oo}{\mathcal{O}_X}
\newcommand{\pic}{\func{Pic}}

\newcommand{\pp}{\mathbb{P}}
\newcommand{\rank}{\func{rank}}
\newcommand{\img}{\func{img}}
\newcommand{\mc}[1]{\mathcal{#1}}
\newcommand{\e}[1]{\mathbf{e}^{#1}}

\newtheorem{theorem}{Theorem}

\newtheorem{claim}[theorem]{Claim}

\newtheorem{corollary}[theorem]{Corollary}

\newtheorem{definition}[theorem]{Definition}

\newtheorem{lemma}[theorem]{Lemma}
\newtheorem{notation}[theorem]{Notation}
\newtheorem{proposition}[theorem]{Proposition}
\newtheorem{remark}[theorem]{Remark}

\begin{document}


\maketitle

\begin{abstract}
    A concrete model for a 7-dimensional gauge theory under special holonomy  is proposed, within
the paradigm of Donaldson and Thomas \cite{Donaldson-Thomas}, over the asymptotically cylindrical $G_2-$manifolds provided by Kovalev's solution to a noncompact version of the Calabi conjecture \cite{kovalevzinho,kovalevzao}. 

One obtains a solution to the $G_2-$instanton equation  from the associated Hermitian Yang-Mills problem, to which the methods of Simpson et al. are applied, subject to a crucial \emph{asymptotic stability} assumption over the `boundary at infinity'. 
\end{abstract}

\tableofcontents

\setcounter{section}{0}
\section*{Introduction}
\addcontentsline{toc}{section}{Introduction}

This article inaugurates a concrete path for the study of 7-dimensional gauge theory in the context of  $G_2-$manifolds. Its core is devoted to the solution of the Hermitian Yang-Mills (HYM) problem for holomorphic bundles over a noncompact Calabi-Yau 3-fold, under certain stability assumptions.

This initiative  fits in the wider context of gauge theory in higher
dimensions, following the seminal works of S. Donaldson, R. Thomas, G. Tian
and others. The common thread to such generalisations is the presence of a closed differential form on the base manifold $M$, inducing an
analogous notion of anti-self-dual connections, or \emph{instantons}, on
bundles over $M$. In the case at hand, $G_{2}-$manifolds are $7-$dimensional
Riemannian manifolds with holonomy in the exceptional Lie group $G_{2}$,
which translates exactly into the presence of such a closed structure. This
allows one to make sense of $G_{2}-$\emph{instantons} as the
energy-minimising gauge classes of connections, solutions to the
corresponding Yang-Mills equation.

While similar theories in dimensions four \cite{4-manifolds} and six \cite{Thomas} have led to the remarkable invariants associated to instanton moduli spaces, little is currently known about the $7-$dimensional case. Indeed, no $G_{2}-$instanton has yet been constructed\footnote{since this article's first version, T. Walpuski \cite{Walpuski} obtained instantons over D. Joyce's $G_2-$manifolds  \cite{Joyce}. }, much less a Casson-type invariant rigorously defined. This is due not least to the success and
attractiveness of the previous theories themselves, but partly also to the
relative scarcity of working examples of $G_{2}-$manifolds \cite{Bryant}\cite{Bryant-Salamon}\cite{Joyce}.
\\

In 2003, A. Kovalev provided an original construction of compact manifolds $M$
with holonomy $G_{2}$ \cite{kovalevzinho,kovalevzao}. These are obtained by gluing two smooth
asymptotically cylindrical Calabi-Yau 3-folds $W'$ and $W''$, truncated
sufficiently far along the noncompact end, via an additional
`twisted'\ circle component $S^{1}$, to obtain a compact Riemannian 7-manifold $M= \left(
W'\times S^{1}\right) \widetilde{\#}\left( W''\times S^{1}\right)$, with possibly `long neck' and $\func{Hol}(M)=G_2$. This
opened a clear, three-step path in the theory of $G_{2}-$instantons: 
\begin{description}
        \item[(1)]

to obtain HYM connections over the $W\ii$, which pull
back to instantons $A\ii$ over  $W\ii\times S^{1}$,

\item[(2)] 
to extend the twisted sum to bundles $\mc{E}\ii\to W\ii$, in order to glue instantons  $A'$ and $A''$ and obtain a $G_{2}-$instanton as $A=A'\widetilde{\#}A''$, say, over the compact base $M  $ \cite{floer}\cite{Taubes}, and 

\item[(3)] 
to study the moduli space of such instantons and compute
invariants in some particular cases.
\end{description}

The present work  is  devoted to the proof of  \emph{Theorem \ref{Thm MAIN THEOREM}}, which guarantees the existence of HYM metrics on suitable holomorphic bundles over such a noncompact 3-fold $W$, and an explicit construction satisfying the relevant hypotheses is also provided, as an example, in \emph{Subsection \ref{Subsec Examples of ASB}}. Thus it completes step (1) of the
above strategy, while postponing the gluing theory for a sequel, provisionally cited as \cite{G2II} and briefly outlined in \emph{Subsection \ref{subsec preview}}.

\setcounter{section}{0}
\section*{Acknowledgements}
\addcontentsline{toc}{section}{Acknowledgements}

This paper derives from the research leading to my PhD thesis and it would not be conceivable without the guidance and support of Simon Donaldson.
I am also fundamentally indebted to Richard Thomas, Alexei Kovalev and Simon Salamon for many important contributions and to the referee, for a very thorough review of the first version. This project was fully funded by the Capes Scholarship (2327-03-1),
 from the Ministry of Education of the Federal Republic of Brazil.

The explicit construction of asymptotically stable bundles is joint work with Marcos Jardim, under a post-doctoral grant (2009/10067-0) by Fapesp - S\~ao Paulo State Research Council. I would finally like to thank Marcos Jardim and Amar Henni for some useful last-minute suggestions.

\setcounter{section}{0}
\section*{Summary}
\addcontentsline{toc}{section}{Summary}

I begin \emph{Section \ref{sec gauge theory on G2 mfds}} recalling basic properties of $G_{2}-$%
manifolds, then move on to adapt gauge-theoretical notions, such as self-duality of $2-$forms under  a $G_2-$structure (\emph{Definition \ref{def G2-structure}}),
topological energy bounds for the Yang-Mills functional and the relationship
between Hermitian Yang-Mills (HYM) connections over a Calabi-Yau $3-$fold $(W,\omega)$,  satisfying $\hat F:=F\cdot\omega=0$,  and so-called $G_{2}-$%
instantons over $W\times S^{1}$ (\emph{Subsections \ref{subsec backgroud G2}} to \emph{\ref{Subsection Hemitian YM}}).
Crucially, the question of finding $G_2-$instantons reduces to the existence of HYM metrics over $W$:

\newtheorem*{prop:G2-HYM}
        {Proposition \ref{prop HYM over W lifts to instanton on WxS1} }
\begin{prop:G2-HYM}
Given a holomorphic vector bundle  $\mathcal{E}\rightarrow W$ over  a Calabi-Yau $3-$fold, the canonical projection $p_1:W\times S^1\to W$ gives a one-to-one correspondence between HYM connections on $\mc{E}$ and $S^1-$invariant $G_{2}-$instantons
on the pull-back bundle $p_1^*\mathcal{E}$.
\end{prop:G2-HYM}

A quick introduction to Kovalev's manifolds \cite{kovalevzinho,kovalevzao} is provided in \emph{Subection \ref{Section ACCY}}, featuring the statement of his noncompact
Calabi-Yau-Tian theorem (\emph{Theorem \ref{Thm SU(3) mfds}}) and a discussion of its ingredients. In a nutshell, given a suitable compact Fano 3-fold $\bar W$, one deletes a K3 surface $D\in \left\vert -K_{\bar{W}}\right\vert$ to obtain a noncompact, asymptotically cylindrical Calabi-Yau manifold $W=\bar{W}\setminus D$. Topologically we interpret $W=W_0\cup W_\infty$ as a compact manifold with boundary $W_0$, with a tubular end $W_\infty$ attached along $\partial W_0$ and the divisor $D$ situated `at infinity'. This sets the scene for a natural evolution problem over  $W$, given by the `gradient flow' towards a HYM solution.

\emph{Sections \ref{sec HYM}} and \emph{\ref{Sec Time-uniform convergence}} form the technical core of the paper. We consider the HYM problem  on a
holomorphic bundle $\mathcal{E}\rightarrow W$ satisfying the
\emph{asymptotic stability} assumption that the restriction $\left.\mc{E}\right\vert_D$ is stable (\emph{Definition \ref{def bundle E->W}}), which allows the existence of a \emph{reference metric} $H_{0}$ with `finite energy' and suitable asymptotic behaviour (\emph{Definition \ref{def reference metric H0}}). This amounts to studying a parabolic equation on the space
of Hermitian metrics \cite{ASDYM,inf dets stab bdls...,Approx instantons} over the
asymptotically cylindrical base manifold (\emph{Subection \ref{Sec evolution
equation on W}}), and it follows a standard pattern 
\cite{Simpson,Guo,Buttler}, leading to the following existence result:
\newpage
\newtheorem*{thm:existence}
        {Theorem \ref{thm sols for all time and exponentially decaying} }
\begin{thm:existence}

Let $\mathcal{E}%
\rightarrow W$ be asymptotically stable, with reference metric $H_{0}
$, over an asymptotically cylindrical $SU\left( 3\right) -$manifold $W$; then,
for any $0<T<\infty $, $\mathcal{E}$ admits a $1-$parameter family $\left\{
H_{t}\right\} $ of smooth Hermitian metrics solving $%
\index{evolution equation!on Wx[0,T]@on $W\times \left[ 0,T\right] $}$%
\begin{equation*}
\left\{ 
\begin{array}{rcl}
H^{-1}%
\dfrac{\partial H}{\partial t}
&=&-2\mathbf{i}\hat{F}_{H} \\ 
H\left( 0\right) 
&=&H_{0}%
\end{array}%
\right. \quad \text{on}\quad W\times \left[ 0,T\right].
\end{equation*}%
Moreover, each $H_{t}$ approaches $H_{0}$ exponentially in all derivatives along the tubular end.
\end{thm:existence} One begins by solving the
associated Dirichlet problem on $W_{S}$, an arbitrary finite truncation of $W$ at `length $S$' along the tube,
first obtaining short-time solutions $H_{S}\left(
t\right) $ (\emph{Subsection \ref{subsec short time existence}}), then
extending them for all time. Fixing arbitrary finite time, one obtains $%
H_{S}\left( t\right) \underset{S\rightarrow \infty }{\longrightarrow }%
H\left( t\right) $ on compact subsets of $W$ (\emph{Subsection \ref{subsec
smooth sols for all time}}). Moreover, every metric in the $1-$parameter family $%
H\left( t\right) $ approaches exponentially the reference metric $H_{0}$, in
a suitable sense, along the cylindrical end (\emph{Subsection \ref{subsec asymp behaviour of sol}}). Hence one has solved the
original parabolic equation and its solution has convenient asymptotia.

From \emph{Section \ref{Sec Time-uniform convergence} }onwards I
tackle the issue of controlling $\lim\limits_{t\rightarrow \infty }H\left(
t\right) $. It is clear from the evolution problem that such a limit $H$, if it exists, will indeed be a HYM metric. Moreover, $H$ has the property of $C^{\infty}-$exponential decay, along cylindrical sections of fixed length, to the reference metric along the tubular end (\emph{Notation \ref{Not finite cylinder}}). The process will culminate in this article's main theorem: 
\newtheorem*{thm:MAIN}
        {Theorem \ref{Thm MAIN THEOREM} }
\begin{thm:MAIN}
In the terms of \emph{Theorem \ref{thm sols for all time
and exponentially decaying}}, the limit $H=\lim\limits_{t\rightarrow\infty}H_t$ exists and is a smooth Hermitian Yang-Mills metric on  $\mathcal{E}$, exponentially asymptotic  in all derivatives to $H_0$ along the tubular end: %
\begin{equation*}       \index{metric!Hermitian Yang-Mills}
\fbox{$\begin{array}{c}
        \hat{F}_{H}=0, \quad H \tublim H_0. \\
\end{array}$}
\end{equation*}
\end{thm:MAIN}

Adapting the `determinant line norm'
functionals introduced by Donaldson \cite{ASDYM,inf dets stab bdls...}%
, I predict in \emph{Claim \ref{claim (F hat)_L^2(Dz) > c  over large set}} a time-uniform lower bound  on the `energy density' $\hat{F}$ over
a finite piece down the tube, of size roughly proportional to $\left\Vert H(t) \right\Vert_{C^0(W)}$   (\emph{Subsection \ref{Subsec conjectured uniform Lp-bound}}). The proof of this fact is quite intricate and it is carried out in detail in \emph{Subsection \ref{Sec proof of conjecture}}. That, in turn,
is a sufficient condition for time-uniform $C^{0}-$convergence of $%
H\left( t\right) $ over the whole of $W$, in view of a negative  energy bound derived by the Chern-Weil method (\emph{Subsection \ref{sect time-uniform convergence}}). Such uniform bound then cascades back through the estimates behind \emph{Theorem \ref{thm sols for all time and exponentially decaying}} and yields a smooth solution to the HYM
problem in the $t\to\infty$ limit, as stated in \emph{Theorem \ref{Thm MAIN THEOREM}}.

Finally, \emph{Subsection \ref{Subsec Examples of ASB}} brings the illustrative example of a
setting $\mathcal{E}\rightarrow W$ satisfying the analytical assumptions of \emph{Theorem \ref{Thm MAIN THEOREM}}, based on a monad construction originally by Jardim \cite{Jardim Bull Braz} and further studied by Jardim and the author \cite{FAE}. 

In view of \emph{Proposition  \ref{prop HYM over W lifts to instanton on WxS1}}, the Chern connection $A_H$ of the HYM metric $H$ obtained in \emph{Theorem \ref{Thm MAIN THEOREM}} pulls back to a $G_2-$instanton over $W\times S^1$, which effectively completes step (1) of the programme outlined in the \emph{Introduction}.

For the reader's convenience, the schemes of proof leading to our key results 
\emph{Theorem \ref{thm sols for all time and exponentially decaying}} and \emph{Theorem \ref{Thm MAIN THEOREM}} are graphically sketched at the end of the paper, in \emph{Appendix \ref{App flowcharts}}.

This article is based on the author's thesis \cite{MyThesis}.
\newpage

\section{$G_{2}-$instantons and Kovalev's tubular construction}
\label{sec gauge theory on G2 mfds}
This \emph{Section} is devoted to the background language for the subsequent analytical investigation. The main references are \cite{Bryant,Salamon,Joyce} and of course \cite{kovalevzinho,kovalevzao}.

\subsection{Background on $G_2-$manifolds}
\label{subsec backgroud G2}
Ultimately we want to consider 7-dimensional Riemannian manifolds with holonomy group $%
G_{2}$. We adopt the conventions of \cite[p.155]{Salamon} for the
definition of $G_{2}$; denoting $\{e^i\}_{i=1,...,7}$ the standard basis of $\left(\mathbb{R}^{7}\right)^*$, $e^{ij}:= e^i\wedge e^j$ etc.:
\begin{definition}      \label{form phi}      
                        \index{G2@$G_{2}$!group}
The \emph{group }$G_{2}$ is the subgroup of $GL(7)$ preserving the $3-$form%
\begin{equation}        \label{eq form phi0}
        \varphi _{0}=\left( e^{12}-e^{34}\right) \wedge e^{5}
        +\left(e^{13}-e^{42}\right) \wedge e^{6}
        +\left( e^{14}-e^{23}\right) \wedge e^{7}+e^{567}
\end{equation}%
under the standard (pull-back) action on $\Lambda ^{3}\left( \mathbb{R}%
^{7}\right) ^{\ast }$, i.e., %
$
        G_{2}:= \left\{ g\in GL(7) \mid 
        g^{\ast}\varphi _{0}=\varphi _{0}\right\} .
$\end{definition}

This encodes various interesting geometrical facts, which are discussed in some detail in the author's thesis \cite{MyThesis} or this article's preprint version  \cite{G2I}, and are thoroughly explored in the above references. In particular, $\varphi _{0}$ defines an Euclidean metric:
\begin{equation}        \label{phi0 gives inner product}
        \left\langle a,b\right\rangle e^{1...7}
        =\left( a\lrcorner \varphi_{0}\right) 
        \wedge \left( b\lrcorner \varphi _{0}\right) 
        \wedge \varphi _{0},
\end{equation}%
and the group   $G_{2}$  has several distinctive properties \cite{Bryant}, which we recall for later use:%
\begin{theorem}         \label{thm properties of G2}
                        \index{G2@$G_{2}$!group}
The subgroup $G_{2}\subset SO(7)\subset GL(7)$
is compact, connected, simple and simply connected and $\dim(G_2)=14$.
Moreover, $G_{2}$ acts irreducibly on $\mathbb{R}^{7}$ and transitively on $%
S^{6}$.
\end{theorem}

\label{Sect G2-manifolds}

Let $M$ be an oriented simply-connected smooth $7-$manifold. 
\begin{definition}      \label{def G2-structure}
\index{G2@$G_{2}$!-structure@$-$structure}A $G_{2}-$\emph{structure} on the $%
7-$manifold $M$ is a $3-$form $\varphi \in \Omega ^{3}\left( M\right) $ such
that, at every point $p\in M$, $\varphi  _{p}=f_{p}^{\ast }\left(
\varphi _{0}\right) $ for some frame $f_{p}:T_{p}M\rightarrow \mathbb{R}^{7}$%
.
\end{definition}

Since $G_{2}\subset
SO\left( 7\right) $ [\emph{Theorem \ref{thm properties of G2}}], $\varphi $ fixes
the orientation given by some (and consequently any) such frame $f$ and also
the metric $g=g\left( \varphi \right) $ given pointwise by (\ref{phi0 gives inner product}).\emph{\ }We may refer indiscriminately to $\varphi$ or $g$ as the $G_{2}-$structure.
The \emph{torsion of }$\varphi $ is the covariant derivative $\nabla
\varphi $ by the induced Levi-Civita connection, and $\varphi $ is \emph{torsion-free }if $\nabla \varphi
=0$.

\begin{definition}
\index{G2@$G_{2}$!-manifold@$-$manifold}A $\emph{G}_{2}\emph{-}$\emph{%
manifold} is a pair $\left( M,\varphi \right) $ where $M$ is a $7-$manifold
and $\varphi $ is a torsion-free $G_{2}-$structure on $M$.
\end{definition}

The following theorem \cite[10.1.3]{Joyce}\cite{Fernandez&Gray} characterises the holonomy reduction on  $G_{2}-$manifolds:
\begin{theorem}         \label{G2 holonomy theorem}
Let $M$ be a $7-$manifold with $G_{2}-$structure $\varphi$ and associated metric $g$; then the following are equivalent:
\begin{enumerate}
        \item $\nabla \varphi =0$, i.e., $M$ is a $G_{2}-$manifold;
        \item $\func{Hol}\left( g\right) $ $\subset G_{2}$;
        \item denoting $\ast _{\varphi }$  the Hodge star from $g(\varphi)$, one has $d\varphi =0$ and $d\ast _{\varphi }\varphi =0$.
\end{enumerate}
\end{theorem}

Finally, on \emph{compact }$G_{2}-$%
manifolds, the holonomy group happens to be \emph{exactly} $G_{2}$ if and only if $\pi _{1}\left(M\right)$ is finite \cite[10.2.2]{Joyce}. The purpose of the `twisted' gluing in A. Kovalev's construction of asymptotically cylindrical $G_2-$manifolds is precisely to secure this topological condition, hence strict holonomy $G_2$.

\subsection{Yang-Mills theory in dimension $7$}

The $G_{2}-$structure allows for a $7-$dimensional analogue of conventional
Yang-Mills theory. The crucial fact is that $%
\varphi _{0}$ yields a notion of (anti-)self-duality for $2-$forms, as $\Lambda ^{2}=\Lambda ^{2}\left( \mathbb{R}%
^{7}\right) ^{\ast }$ splits into irreducible representations.

Since $G_{2}\subset SO\left( 7\right)$, we have $\mathfrak{g}%
_{2}\subset \mathfrak{so}\left( 7\right) \simeq \Lambda ^{2}$, under the
standard identification of $2-$forms with antisymmetric matrices. Denote $\Lambda _{-}^{2}:= \mathfrak{g}_{2}$
and $\Lambda _{+}^{2}$ its orthogonal complement in $\Lambda ^{2}$:%
\begin{equation}
\Lambda ^{2}=\Lambda _{+}^{2}\oplus \Lambda _{-}^{2}.  \label{split}
\end{equation}%
Then $\dim \Lambda _{+}^{2}=7$, and we identify $\Lambda _{+}^{2}\simeq \mathbb{R}^{7}$ as the linear span
of the contractions $\alpha _{i}:= v_{i}\lrcorner \varphi _{0}$. Indeed
the  $G_{2}-$action on $\Lambda _{+}^{2}$ translates to the
standard one on $\mathbb{R}^{7}$:%
\begin{eqnarray*}
\left( g.\alpha _{i}\right) (u_{1},u_{2}) &=&\alpha _{i}\left(
g.u_{1},g.u_{2}\right) =\varphi _{0}\left( v_{i},g.u_{1},g.u_{2}\right) \\
&=&\varphi _{0}\left( g^{-1}.v_{i},u_{1},u_{2}\right) =\left( \left(
g^{-1}.v_{i}\right) \lrcorner \varphi _{0}\right) \left( u_{1},u_{2}\right) .
\end{eqnarray*}%
One checks easily that $\alpha _{i}\in \left( \mathfrak{g}%
_{2}\right) ^{\perp }\subset \mathfrak{so}\left( 7\right) $;  moreover \cite[p. 541]{Bryant}:
\begin{claim}
\label{+2 -1 eigenspaces}The space $\Lambda _{\pm }^{2}$ has the following
properties:

\begin{enumerate}
\item $\Lambda _{\pm }^{2}$ is an irreducible representation of $G_{2}$;

\item $\Lambda _{\pm }^{2}$ is the $_{-1}^{+2}-$eigenspace of the  $%
G_{2}-$equivariant linear map%
        \begin{eqnarray*}
        T : \Lambda ^{2} &\rightarrow& \Lambda ^{2} \\
        \eta &\mapsto &T\eta =\ast \left( \eta \wedge \varphi _{0}\right) .
\end{eqnarray*}
\end{enumerate}
\end{claim}

By analogy with the $4-$dimensional case, we will call $\Lambda _{+}^{2}$
(resp. $\Lambda _{-}^{2}$) the space of \emph{self-dual }or \emph{SD} (resp. \emph{anti-self-dual} or \emph{ASD}) forms. 
Still in the light of \emph{Claim \ref{+2 -1 eigenspaces}}, there is
a convenient characterisation of the `positive' projection in $\left( \ref%
{split}\right) $. The $4-$form  dual to the $G_{2}-$structure $%
\left( \ref{form phi}\right) $ in our convention is%
\begin{equation}
\ast \varphi _{0}=\left( e^{34}-e^{12}\right) \wedge e^{67}+\left(
e^{42}-e^{13}\right) \wedge e^{75}+\left( e^{23}-e^{14}\right) \wedge
e^{56}+e^{1234}  \label{form *phi}
\end{equation}%
and we consider the  $G_{2}-$equivariant map (between representations of $G_2$):
\begin{eqnarray*}        \label{eq Lphi0}
        L_{\ast \varphi _{0}} : \Lambda ^{2} 
        &\rightarrow& \Lambda ^{6}\\ 
        \eta &\mapsto& \eta \wedge \ast \varphi _{0}.
\end{eqnarray*}
Since $\Lambda _{\pm }^{2}$ and $%
\Lambda ^{6}$ are irreducible representations and $\dim \Lambda
_{+}^{2}=\dim \Lambda ^{6}$, Schur's Lemma gives:

\begin{proposition}     \label{prop +-projection}
The above map restricts to $\Lambda^2_{\pm}$ as $$\left. L_{\ast \varphi _{0}} \right\vert _{\Lambda
_{+}^{2}}:\Lambda _{+}^{2} \: \tilde{\rightarrow} \: \Lambda^{6} 
\quad\text{and}\quad\left. L_{\ast \varphi _{0}}\right\vert _{\Lambda
_{-}^{2}}=0.$$
\begin{proof}
It only remains to check that the restriction $\left. L_{\ast \varphi _{0}} \right\vert _{\Lambda_{+}^{2}}$ is nonzero. Using $\left( \ref{form phi}\right) $ and $\left( \ref{form *phi}\right)$ we find, for instance, 
\begin{eqnarray*}
        L_{\ast \varphi _{0}}\alpha _{1} 
        &=&\left( v_{1}\lrcorner \varphi_{0}\right)\wedge\ast\varphi_{0}\\
        &=&\left( e^{25}+e^{36}+e^{47}\right) \wedge \ast\varphi_{0} \\
        &=&e^{234567}.
\end{eqnarray*}%
Not only does this prove the statement, but it also suggests carrying out the
full inspection of the elements $L_{\ast \varphi _{0}}\alpha _{i}$, which yields the somewhat aesthetical fact:
\begin{equation*}
        L_{\ast \varphi _{0}}\alpha _{i}=e^{1...\hat{\imath}...7}
        \qedhere
\end{equation*}
\end{proof}
\end{proposition}

Hence we may think of $L_{\ast \varphi _{0}}$, the `wedge product with $\ast\varphi _{0}$', as the orthogonal projection of $2-$forms into the subspace $\Lambda
_{+}^{2}\simeq \Lambda ^{6}$.

\label{Sub energy bounds}
Consider now a vector bundle $E\rightarrow M$ over a compact $G_{2}-$manifold $\left(
M,\varphi \right) $; the curvature $F_{A} $ of some connection $A$ conforms to
the splitting $\left( \ref{split}\right) $:%
\begin{equation*}
        F_{A}=F_{A}^{+}\oplus F_{A}^{-},
        \qquad  F_{A}^{\pm }\in \Omega _{\pm }^{2}\left(\func{End}E\right).
\end{equation*}%
The $L^2-$norm of $F_{A}$, if it is well-defined (e.g. if $M$ is compact), is the \emph{Yang-Mills functional}:%
\begin{equation}        \label{YM(A)}
\begin{array}{rcccl}
        YM\left( A\right) 
        &:=& \Vert F_{A}\Vert ^{2}
        &=&\Vert F_{A}^{+}\Vert ^{2}
        +\Vert F_{A}^{-}\Vert ^{2}.
\end{array}
\end{equation}

It is well-known that the values of $YM\left( A\right) $ can be related to a certain
characteristic class of the bundle $E$,%
\begin{equation*}
        \kappa \left( E\right) :=\int_{M}%
        \func{tr}\left( F_{A}^{2}\right) \wedge \varphi .
\end{equation*}%
Using the property $d\varphi =0$, a standard argument of Chern-Weil theory \cite{Milnor} shows that $\left[ \func{tr}\left(
F_{A}^{2}\right) \wedge \varphi \right]^{dR}$ is independent of $A$, thus the
integral is indeed a topological invariant. Using  the eigenspace decomposition from \emph{Claim \ref{+2 -1 eigenspaces}} we find:%
\begin{eqnarray*}        
        \kappa \left( E\right) &=&
        -\int_{M}\left\langle F_{A}\wedge 
        \left( F_{A}\wedge\varphi\right) \right\rangle _{\mathfrak{g}}
        \notag 
        \ =\ -\left(F_{A},2F_{A}^{+}-F_{A}^{-}\right)  \notag \\
        &=&\Vert F_{A}^{-}\Vert ^{2}
        -2\Vert F_{A}^{+}\Vert^{2}.
\end{eqnarray*}
Comparing with $\left( \ref{YM(A)}\right) $ we get:%
\begin{equation*}
        YM(A)=3\Vert F_{A}^{+}\Vert ^{2}
        +\kappa \left( E\right) 
        =\frac{1}{2}\left( 3\Vert F_{A}^{-}\Vert ^{2}
        -\kappa \left( E\right) \right)
\end{equation*}%
As expected,  $YM\left( A\right) $ attains its absolute minimum at a
connection whose curvature is either SD or ASD. Moreover, since $YM\geq0$, the sign of  $\kappa(E)$ obstructs the existence of one type or the other. Fixing $\kappa(E)\geq0$, these facts motivate our interest in the 
$G_{2}-$\emph{instanton equation}:%
\begin{equation}        \label{G2-instanton}
                        \index{G2@$G_{2}$!-instanton@$-$instanton}
        \boxed{F_{A}^{+}:=\left(\left.L_{\ast\varphi}
        \right\vert_{\Omega^2_+}\right)^{-1} \left(F_{A}\wedge*\varphi\right)=0}
\end{equation}
where $L_{\ast\varphi}$ is given pointwise by $L_{\ast\varphi_0}$, as above in page \pageref{eq Lphi0}.

\newpage
\subsection{Relation with $3-$dimensional Hermitian Yang-Mills}  
\label{Subsection Hemitian YM}

We now switch, for a moment, to the complex geometric picture and consider a holomorphic vector bundle $\mathcal{E}\rightarrow W$ over a K\"ahler manifold $\left( W,\omega \right)$. To every Hermitian metric $H$ on $\mathcal{E}$ there corresponds a unique  compatible (Chern) connection $A=A_H$, with 
$F_{A}\in\Omega^{1,1}\left( \func{End}\mc{E} \right).$ In this context, the \emph{Hermitian Yang-Mills (HYM)} condition is the vanishing of the  $\omega-$trace: 
\begin{equation}
\hat F_A:=F_A\cdot\omega=0\in\Omega^{0}\left( \func{End}\mc{E} \right).
\end{equation}
If $W$ is a Calabi-Yau $3-$fold, the Riemannian product $M=W\times S^{1}$ is
naturally a real $7-$dimensional $G_{2}-$manifold \cite[p.564]{Bryant}\cite[11.1.2]{Joyce}. In this \emph{Subsection}, we will check the corresponding gauge-theoretic fact that HYM 
\index{metric!Hermitian Yang-Mills}connections on $\mathcal{E}\rightarrow W$ pull back to $G_{2}-$%
instantons over the product $M$.

Indeed, given the K\"{a}hler form $\omega  $ and holomorphic volume form $\Omega $ on $W$, we obtain a natural $G_2-$structure [ibid.] as follows:%
\begin{equation}        \label{kahler -> G2}
\begin{array}{r c l}
        \varphi &:=&\omega \wedge d\theta +\func{Im}\Omega,\\
        \ast \varphi &=&\tfrac{1}{2}\omega\wedge\omega 
        -\func{Re}\Omega \wedge d\theta .
\end{array}
\end{equation}
Here $d\theta $ is the coordinate $1-$ form on $S^{1}$, and the Hodge
star on $M$ is given by the product of the K\"{a}hler metric on $W$ and the
standard flat metric on $S^{1}$. 

Now, a connection $A$ on $\mathcal{E}\rightarrow W$  pulls back to  $p_1^*\mathcal{E}\rightarrow M$
via the canonical projection 
$$
p_1:W\times S^{1}\to W,
$$and so do the forms $\omega $ and $\Omega $ (for simplicity I keep the same
notation for objects on $W$ and their pull-backs to $M$). In particular, under the isomorphism $\left. L_{\ast \varphi} \right\vert _{\Omega
_{+}^{2}}:\Omega_{+}^{2} \, \tilde{\rightarrow} \; \Omega^{6}$ [\emph{Proposition \ref{prop +-projection}}], the
SD part of curvature maps to 
\begin{equation}        \label{sd part kahler -> G2}
        L_{\ast\varphi}\left(F_{A}^{+}\right)
        =\ F_A\wedge\ast\varphi
        =\tfrac{1}{2}F_{A}\wedge \left( \omega \wedge \omega -2\func{Re}%
        \Omega \wedge d\theta \right) . 
\end{equation}

\begin{proposition}
\label{prop HYM over W lifts to instanton on WxS1}
Given a holomorphic vector bundle  $\mathcal{E}\rightarrow W$ over  a Calabi-Yau $3-$fold, the canonical projection $p_1:M=W\times S^1\to W$ gives a one-to-one correspondence between Hermitian Yang-Mills connections on $\mc{E}$ and $S^1-$invariant $G_{2}-$instantons
on the pull-back bundle $p_1^*\mathcal{E}$.
\begin{proof}
A HYM connection $A$ satisfies $F_{A}\in \Omega ^{1,1}\left( W\right) 
\text{and }\hat{F}_{A}=0$. Taking account of bidegree, the former implies $%
F_{A}\wedge \Omega =$ $F_{A}\wedge \overline{\Omega }=0$, hence
\begin{displaymath}
        F_{A}\wedge 2 \func{Re}\Omega 
        =F_{A}\wedge \left( \Omega +\overline{\Omega }\right) 
        =0.
\end{displaymath}
Replacing this in $\left( \ref{sd part kahler -> G2}\right)$, we check that $F_{A}^{+}$ maps isomorphically to the origin:
\begin{displaymath}
        \begin{array}{cclcc}
        F_{A}^{+} & \cong & \tfrac{1}{2} F_{A}\wedge \omega \wedge \omega 
                & \in & \Omega^{3,3}\left( W\right) \\
        &=& \left( cst.\right) \hat{F}_{A}\otimes d\vol\left( W\right)&&\\
        &=& 0 &&\\
        \end{array}
\end{displaymath}
using the HYM condition $\hat{F}_{A}=0$ and  $\omega
\wedge \omega =\frac{\left( cst.\right) }{\left\Vert \omega \right\Vert ^{2}}%
\ast \omega $.
\end{proof}
\end{proposition}

Thus, by solving the HYM
equation over $CY^{3}$, one obtains $G_{2}-$instantons over the product $CY^{3}\times S^1$, which is this article's motivation.
For some further discussion of $G_2-$manifolds of the form $CY^3\times S^1$, see \cite{Baraglia}. 

\subsection{Asymptotically cylindrical Calabi-Yau $3-$folds}
\label{Section ACCY}
I will give a brief account of the
building blocks in the construction of compact $G_{2}-$manifolds by A. Kovalev \cite{kovalevzinho,kovalevzao}.
These are achieved by
gluing together, in an ingenious way, a pair of noncompact asymptotically
cylindrical $7-$manifolds of holonomy $SU\left( 3\right) $ along their
tubular ends. Such components are of the form $W\times S^{1}$, where $\left(W,\omega\right)$ is $3-$%
fold given by a noncompact version of the Calabi conjecture, thus they carry $G_{2}-$structures as  in 
\emph{Subsection \ref{Subsection Hemitian YM}}. 

\begin{definition}      \label{def base manifold (W,w)}
                        \index{base manifold}%
A \emph{base manifold} for our purpose is a compact, simply-connected K\"{a}hler $3-$fold $\left( \bar{W},\bar{\omega}%
\right) $ admitting the following:%
\begin{itemize}         
        \item  
        a $K3-$surface 
        $D\in \left\vert -K_{\bar{W}}\right\vert $ (simply-connected, compact, $c_{1}\left( D\right)=0$) with holomorphically trivial         normal bundle $\mathcal{N}_{D/\bar{W}}$;
        \item
        the complement $W=\bar{W}\setminus D$ has finite fundamental group         $\pi _{1}\left( W\right) $.
\end{itemize}
\end{definition}

One wants to think of $W$ as a compact manifold $W_{0}$ with boundary $%
D\times S^{1}$ and a topologically cylindrical end attached there:%
\begin{equation}        \label{eq cylindrical picture}
\begin{array}{c}
        W=W_{0}\cup W_{\infty } \\ 
        W_{\infty }\simeq \left( D\times S^{1}
        \times \mathbb{R}_{+}\right) .%
\end{array}
\end{equation}
\begin{figure}[h]
        \begin{center}
        \includegraphics{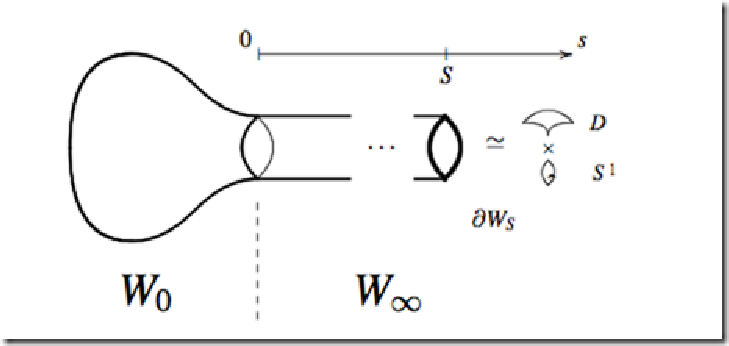}
        \end{center}
        \caption{Asymptotically cylindrical Calabi-Yau 3-fold $W$}
\end{figure}
\noindent Here   $W_S$ denotes the truncation of $W$ at `length $S$' on the $\mathbb{R}_+$ component.
\\\

Let $s_{0}\in H^{0}\left( \bar{W},K_{\bar{W}}^{-1}\right) $ be the defining
section of the \emph{divisor at infinity} $D$; then, $s_{0}$ defines a holomorphic coordinate $z$ on a
neighbourhood $U\subset\bar{W} $ of $D$. Since $%
\mathcal{N}_{D/\bar{W}}$ is trivial, we may assume $U$ is a \emph{tubular
neighbourhood of infinity}, i.e.,%
\begin{equation}        \index{infinity!neighbourhood of}
                        \label{eq neighbourhood of infinity}
        U\simeq D\times \left\{ \left\vert z\right\vert <1\right\}
\end{equation}%
as real manifolds. Denoting $s\in \mathbb{R}_{+}$ and $\alpha \in S^{1}$, we pass to the asymptotically cylindrical picture $\left( \ref{eq cylindrical
picture}\right) $ via $z=e^{-s-\mathbf{i}\alpha }$. For later reference,
let us establish a straightforward result on the decay of differential forms along $W_\infty$:
\begin{lemma}   \label{lemma decay of coordinates dz dz bar}
With respect to the model cylindrical metric, 
\begin{equation*}
        \left\vert dz\right\vert ,\left\vert d\bar{z}\right\vert 
        =O\left( \left\vert \e{-s}\right\vert \right) .
\end{equation*}
\begin{proof}
In holomorphic coordinates $\left( z,\xi ^{1},\xi ^{2}\right) $ around $%
D\subset U$ with $z=\e{-s-\mathbf{i}\alpha }$ and $%
D=\left\{ z=0\right\} $, we have $\left\vert z\right\vert =\left\vert \bar{z}\right\vert =\e{-s}$ and%
\begin{eqnarray*}
        dz &=&-z\left( ds-\mathbf{i}d\alpha \right) \\
        d\bar{z} &=&-\bar{z}\left( ds+\mathbf{i}d\alpha \right).
\end{eqnarray*}%
\vspace{-1.3cm}\[\qedhere\]
\end{proof}
\end{lemma}

Furthermore, the holomorphic coordinate $z$ on $U$ is the same as a local
holomorphic function $\tau $, say. The assumptions that $\bar{W}$ is simply connected, compact and K\"ahler imply the vanishing in Dolbeault cohomology \cite[Corollary 3.2.12, p.129]{Huybrechts} \begin{displaymath}
        H^{0,1}\left( \bar{W}\right)\oplus \overline{H^{0,1}\left( \bar{W}\right)}
        =H^1(\bar{W},\mathbb{C})=0,
\end{displaymath} 
in which case one can solve Mittag-Leffler's problem for $\dfrac{1}{\tau }$ on $\bar{W}$ \cite[pp.34-35]{Griffiths & Harris} and $\tau$ extends to
a global fibration%
\begin{equation}        \label{eq fibration tau}
        \tau :\bar{W}\overset{D}{\longrightarrow }\mathbb{C}P^{1}
\end{equation}%
with generic fibre a  $K3-$surface (diffeomorphic  to $D$) and some singular fibres. In fact, this
holomorphic coordinate can be seen as pulled-back from $\mathbb{C}P^{1}$, i.e.,  $K_{\bar{W}}^{-1}$ is the pull-back of a
degree-one line bundle $L\rightarrow \mathbb{C}P^{1}$ and $z=\tau ^{\ast
}s_{0}$ for some $s_{0}\in H^{0}\left( L\right) \cite[\S 3]{kovalevzao}$.

\label{Subsec Noncompact CY}
Finally, in order to state Kovalev's noncompact version of the Calabi
conjecture, notice that the $K3$ divisor $D$ has a complex
structure $I$ inherited from $\bar{W}$; by Yau's theorem, it admits a unique
Ricci-flat K\"{a}hler metric  $\kappa _{I}\in\left[\left. \bar{\omega}%
\right\vert _{D}\right]$. Ricci-flat  K\"{a}%
hler metrics on a complex surface%
\index{K3 surface!hyper-K\"{a}hler} are \emph{hyper-K\"{a}hler} \cite[pp.336-338]{BHPV}, which means $D$ admits
additional complex structures $J$ and $K=IJ$ satisfying the quaternionic
relations, and the metric is also K\"{a}%
hler with respect to any combination $aI+bJ+cK$ with $\left( a,b,c\right)
\in S^{2}$. Let us denote their K\"{a}hler forms by $%
\kappa _{J}$ and $\kappa _{K}$. In those terms we have  \cite[Theorem 2.2]{kovalevzinho}:
\begin{theorem}[Calabi-Yau-Tian-Kovalev]         \label{Thm SU(3) mfds}
For $W=\bar{W}\setminus D$ as in \emph{Definition \ref{def base manifold (W,w)}}: 
\begin{enumerate}
        \item
        $W$ admits a complete Ricci-flat K\"{a}hler metric $\omega $;
        \item
        along the cylindrical end $D\times S_{\alpha }^{1}\times 
        \left(\mathbb{R}_{+}\right) _{s}$, the K\"{a}hler form $\omega $
         and its holomorphic volume form $\Omega $ are 
        \emph{exponentially asymptotic} to those of the (product) cylindrical metric induced from $D$:%
        \begin{eqnarray}
        \begin{array}{r c l}      \label{asymptotic forms}
                \omega _{\infty } &=&\kappa _{I}+ds\wedge d\alpha \\
                \Omega _{\infty } &=&\left( ds+\mathbf{i}d\alpha \right)                 \wedge \left(\kappa _{J}+\mathbf{i}\kappa _{K}\right) \text{;}
        \end{array}
        \end{eqnarray}
        \item
        $\func{Hol}\left( \omega \right) =SU\left( 3\right) $, i.e., 
        $W$ is Calabi-Yau.
\end{enumerate}
\end{theorem}
By \emph{exponentially asymptotic} one means precisely that the forms can be written along the tubular end as
\begin{displaymath}
        \left.\omega \right\vert_{W_{\infty}}
        =\omega _{\infty }+d\psi, \qquad
        \left.\Omega \right \vert_{W_\infty}
        =\Omega _{\infty }+d\Psi 
\end{displaymath}%
where the $1-$form $\psi $ and the $2-$form $\Psi $ are smooth 
and decay as $O\left( \e{-\lambda s}\right) $ in all derivatives with respect
to $\omega _{\infty }$, for any $\lambda <\min \left\{ 1,
\sqrt{\lambda_{1}\left( D\right) }\right\} $, and 
$\lambda _{1}\left( D\right) $ is the first eigenvalue of the Laplacian
on differential forms on $D$, with the metric $\kappa _{I}$.

\newpage
\section{Hermitian Yang-Mills problem}
\label{sec HYM}

The interplay between the algebraic and differential geometry of vector bundles is one of the main aspects of  gauge theory, dating back at least to Narasimhan and Seshadri's famous correspondence between stability and flatness over Riemann surfaces.    In the general K\"ahler setting, given a bundle $\mc{E}$ over $(W,\omega)$, Hitchin and Kobayashi suggested a generalisation of the flatness condition for a Hermitian metric $H$, in terms of the natural contraction of $(1,1)-$forms, expressed by the \emph{Hermite-Einstein (HE)} condition:
$$
\hat{F}_{H}:=F_H\cdot \omega =\lambda.Id\in\Omega^0(\func{End}\mc{E}),
\qquad\lambda\in \mathbf{i}\R,
$$
where $\lambda$ is proportional to the slope of $\mc{E}$. The one-to-one correspondence between indecomposable HE connections and stable holomorphic structures was then established by Donaldson  over projective algebraic surfaces in \cite{ASDYM} and extended to compact K\"ahler manifolds of any dimension by Uhlenbeck and Yau, with a simpler proof in the projective case again by Donaldson in \cite{inf dets stab bdls...}.  The methods organised in those two articles rely on the interpretation of HE metrics as critical points of the Yang-Mills functional (henceforth \emph{Hermitian Yang-Mills} metrics), which allows for the application of PDE techniques, inspired by Eeels and Sampson \cite{Eels and Sampson}, to prove convergence and regularity of the associated gradient flow. That set of tools forms the theoretical backbone of the present investigation.

Several noncompact variants of the question unfold in the work of numerous authors, but two texts in particular have an immediate, pivotal relevance to my proposed  problem. Guo's study of instantons over certain cylindrical 4-manifolds \cite{Guo} shows that one may obtain uniform bounds on the `heat kernel' of the parabolic flow under suitable asymptotic conditions on the `bundle at infinity', hence prove regularity of solutions by bootstrapping, as in the compact case. On the other hand, Simpson explores the more general, higher-dimensional noncompact case    \cite{Simpson}, assuming an exhaustion by compact subsets and truncating `sufficiently far' with fixed Dirichlet conditions, then taking the limit in the family of solutions obtained over each compact manifold with boundary.
Both insights will be effectively combined in the technical argument to follow.
\subsection{\texorpdfstring
        {The evolution equation on $W$}
        {The evolution equation on W}}
\label{Sec evolution equation on W}

Let $\left( W,\omega \right) $ be an asymptotically cylindrical Calabi-Yau $3-$fold as given by \emph{Theorem \ref{Thm SU(3) mfds}} and $\mathcal{E}%
\rightarrow W$ the restriction of a holomorphic vector bundle on $\bar{W}$ under certain stability assumptions. Our guiding thread is the perspective of
obtaining a smooth Hermitian metric $H$ on $\mathcal{E}$ satisfying the 
\emph{Hermitian Yang-Mills} condition%
\begin{equation}        \label{Eq HYM}
        \hat{F}_{H}=0\in\Omega^0(\func{End}\mc{E}),
\end{equation}%
which would  solve the associated $G_2-$instanton equation on $W\times S^1$ [\emph{Proposition \ref{prop HYM over W lifts to instanton on WxS1}}]. Consider thus the following analytical problem.
 
Let $W_{S}$ be the compact manifold (with boundary) obtained by
truncating $W$ at length $S$ down the tubular end. On each $W_{S}$ we consider the nonlinear
`heat flow'%
\index{evolution equation!on WSx[0,T[@on $W_{S}\times \left[ 0,T\right[ $}%
\begin{equation}
\left\{ 
\begin{array}{r c l}
        H^{-1}\dfrac{\partial H}{\partial t}&=&-2\mathbf{i}\hat{F}_{H}\\
        H\left( 0\right) &=&H_{0}%
\end{array}%
\right. \text{\quad on\quad }W_{S}\times \left[ 0,T\right[   \label{Eq heat}
\end{equation}%
with smooth solution $H_{S}\left( t\right) $, defined for some (short) $T$, since $\left( \ref{Eq heat}\right) $ is parabolic. Here  $H_{0}$ is
a fixed metric on $\mathcal{E}\rightarrow W$ extending  the pullback `near infinity' of an ASD connection over $D$ [see \emph{Definition \ref{def reference metric H0}} below], and one imposes the \emph{Dirichlet} boundary condition
\begin{equation}        \label{Dirichlet}
        H\mid _{\partial W_{S}}=H_{0}\mid _{\partial W_{S}}.  
\end{equation}%
Taking suitable $t\rightarrow T$ and $S\rightarrow \infty $ limits of
solutions to $\left( \ref{Eq heat}\right) $ over compact subsets $W_{S_0}$, we obtain a solution $%
H\left( t\right) $ to the evolution equation  defined over $W$, for arbitrary  $T<\infty$, with two key properties:
\begin{itemize}
        \item 
        each metric $H\left( t\right) $ is exponentially asymptotic in all
derivatives to $H_{0}$ over typical finite cylinders along the tubular end;
        \item 
        if $H\left( t\right) $ converges as $t\leq T\rightarrow \infty$, then the limit is a HYM metric on $\mathcal{E}$.
\end{itemize}
 
Moreover, the infinite-time convergence of $H\left(
t\right) $ over $W$  can be
reduced to establishing a lower `energy bound' on $\hat{F}_{H(t)}$, over a
domain down the tube of length roughly proportional to $\left\Vert H(t)\right\Vert_{C^0(W)}$. That result  is  eventually proved, yielding a smooth HYM solution as the $t\to\infty$ limit of the family $H(t)$.

\subsection{The concept of asymptotic stability}
We start with a holomorphic vector bundle $\mathcal{E}$
over the original compact K\"{a}hler $3-$fold $\left( \bar{W},\bar{\omega}%
\right) $ [cf. \emph{Section \ref{Section ACCY}}] and we ask that its
restriction  $\left. \mathcal{E}\right\vert _{D}$ to the \emph{divisor at infinity} $D$ \index{infinity!divisor at}be (slope-)stable with respect to $\left[ \bar{\omega}\right] $. As stability is an open condition, it also holds over the nearby K3 `slices' in the neighbourhood of infinity $U\supset D$ [cf. (\ref{eq neighbourhood of infinity})], which we denote:
$$
D_z:=D\times\left\{ z \right\}\subset U.
$$

So there exists $\delta >0$
such that $\left. \mathcal{E}\right\vert _{D_{z}}$ is also stable for all $%
\left\vert z\right\vert <\delta $. In view of this hypothesis, we will say
colloquially that $\left. \mathcal{E}\right\vert _{\bar{W}\smallsetminus D}$, denoted simply $\mathcal{E}\rightarrow W$, is \emph{stable at
infinity }over the noncompact Calabi-Yau $\left( W,\omega \right) $ given by 
\emph{Theorem \ref{Thm SU(3) mfds}}; this is consistent since $\left[
\omega \right] =\left[ \left. \bar{\omega}\right\vert _{W}\right] $. Such asymptotic stability assumption will be crucial in establishing the time-uniform `energy' bounds on the solution to the evolution problem [cf. \emph{Lemma }\ref{Lemma NDz > xi_L^4/3}]. In
summary:

\begin{definition}
\label{def bundle E->W}A bundle $\mathcal{E}\rightarrow W$ will be called 
\emph{stable at infinity} (or \emph{asymptotically stable}) if it is the
restriction of a holomorphic vector bundle $\mathcal{E}\rightarrow \bar{W}$
such that $\left. \mathcal{E}\right\vert _{D}$ is stable, hence also $\left. \mathcal{E}\right\vert _{D_{z}}$ for $\left\vert z\right\vert <\delta $ (in particular, $\mathcal{E}$ is indecomposable as a direct sum.)%
\end{definition}
Explicit examples satisfying \emph{Definition \ref{def bundle E->W}} can be obtained for instance by monad techniques, following the ADHM paradigm [cf. \emph{Subsection \ref{Subsec Examples of ASB}}]. 

The last ingredient is a suitable metric for comparison. Fix a smooth trivialisation of $\mathcal{E} \vert_{U}$ `in the $z-$direction' over the neighbourhood of infinity $U\simeq D\times \left\{ \left\vert z\right\vert <1 \right\}\subset \bar{W}$, i.e., an isomorphism 
$$
\left. \mathcal{E}\right\vert _{U}\simeq\left\{ \left\vert z\right\vert <1 \right\}\times \left. \mathcal{E}\right\vert _{D}.
$$ 
Define $\left. H_{0}\right\vert _{D}$ as the Hermitian Yang-Mills metric on $\left. \mathcal{E}\right\vert _{D}$ and denote by $K$ its pull-back to $\left.\mathcal{E}\right\vert _{U}$ under the above identification. For definiteness, let us fix 
\begin{eqnarray}        \label{eq fixed det over D}
        \det \left. H_{0}\right\vert _{D}=\det K=1.
        \end{eqnarray}
Indeed, for the purpose of Yang-Mills theory we may always assume $\det H=1$, since the $%
L^{2}-$norm of $\func{tr}F_{H}$ is minimised independently by the
harmonic representative of $\left[ c_{1}\left( \mathcal{E}\right) \right]
^{dR}$.

Then, for each $0\leq \left\vert z\right\vert <\delta $, the stability assumption gives a self-adjoint element $h_z \in \func{End} \left( \left.\mathcal{E}\right\vert_{D_{z}} \right)$ such that $\left. H_{0}\right\vert _{D_{z}} := {K}.h_z$ is the HYM metric on $\left. \mathcal{E}\right\vert _{D_{z}}$. Supposing, for simplicity, that $\mathcal{E}$ is an $SL(n,\mathbb{C})-$bundle and fixing $\det h_z$, I  claim the family $h_z$ varies smoothly with $z$ across the fibres. This can be read from the HYM condition; writing  $F_{0\vert z}$ for the curvature of the Chern connection of $\left. H_{0}\right\vert _{D_{z}}$ and $\Lambda_z$ for the contraction with the K\"{a}hler form on $D_z$ (hence, by definition, $\hat{F} := \Lambda_z F$), we have:
\begin{eqnarray*}
        \widehat{F_{0\vert z}}=0 &\Leftrightarrow& P_z \left( h_z \right):=         \Delta_K h_{z}+\mathbf{i}\left( \hat{F}_K.h_z + h_z.\hat{F}_K \right)         + 2\mathbf{i}\Lambda_z\left(\bar{\partial}_K h_z.h_z^{-1}.\partial_Kh_z         \right) = 0
\end{eqnarray*}
where $\left\{P_z\right\}$ is a family of nonlinear partial differential operators depending smoothly on $z$. Linearising \cite[p.14]{ASDYM} and using the assumption $\widehat{F_{0\vert 0}}=0$ we find $\left.{\left(\delta_{z}P \right)} \right\vert_{z=0}=\Delta_0$, the Laplacian of the reference metric, which is invertible on metrics with fixed determinant. This proves the claim, by the implicit function theorem; namely, we obtain a smooth bundle metric $H_0$ over the neighbourhood $U$ which is `slicewise' HYM over each cross-section $D_z$.

Extending $H_0$ in any smooth way over the compact end $\bar{W} \setminus U$, we obtain
a
smooth Hermitian bundle metric on the whole of $\mathcal{E}$. For technical reasons,
I also require $H_{0}$ to have \emph{finite energy} [cf. equation $(\ref{Eq asymptotics of F_H0})$ in \emph{Subsection \ref{subsec asymp behaviour of sol}}]:%
\begin{equation*}
        \Vert \hat{F}_{H_{0}}\Vert _{L^{2}\left( W,\omega\right) }<\infty.
\end{equation*}

\begin{definition}      \label{def reference metric H0}
A \emph{reference metric} $H_{0}$ on the asymptotically stable bundle 
$\mathcal{E}\rightarrow W$ is (the restriction of) a smooth Hermitian metric on $\mathcal{E}\rightarrow \bar{W}$ such that:%
\begin{itemize}         \index{metric!reference $H_0$}
        \item 
        $\left. H_{0}\right\vert _{D_{z}}$ are the corresponding HYM metrics on $\left. \mathcal{E}\right\vert _{D_{z}}$, $0\leq
\left\vert z\right\vert <\delta $;
        \item
        $\det H_0\equiv1$;
        \item 
        $H_{0}$ has finite energy.
\end{itemize}
\end{definition}

\begin{remark}  \label{rem deformation complex}
Denote $A_{0}$ the Chern connection of $H_{0}$, relative to the fixed holomorphic structure of $\mathcal{E}$. Then by assumption each $\left. A_{0}\right\vert _{D_{z}}
$ is ASD \cite[p.47]{4-manifolds}. In particular, $\left.
A_{0}\right\vert _{D}$ induces an elliptic deformation complex%
\index{elliptic!complex}
\begin{equation*}
\Omega ^{0}\left( \mathfrak{g}\right) \overset{d_{A_{0}}}{\rightarrow }%
\Omega ^{1}\left( \mathfrak{g}\right) \overset{d_{A_{0}}^{+}}{\rightarrow }%
\Omega _{+}^{2}\left( \mathfrak{g}\right) 
\end{equation*}%
where $\mathfrak{g}=%
\func{Lie}\left( \left. \mathcal{G}\right\vert _{D}\right) $ generates
the gauge group $\mathcal{G}$ $=\func{End}\mathcal{E}$ over $D$. Thus,
the requirement that $\mathcal{E}\vert_D$ be indecomposable imposes a constraint on cohomology: 
\begin{equation*}
        \mathbf{H}_{A_{0}\vert_D}^{0}=0,
\end{equation*}%
as a non-zero horizontal section would otherwise split $\left. \mathcal{E}\right\vert _{D}$.

Furthermore, although this will not be essential here, it is worth observing that one might want to restrict
attention to \emph{acyclic }connections \cite[p.25]{floer},
i.e., whose gauge class $\left[ A_{0}\right] $ is \emph{isolated} in $%
\mathcal{M}_{\left. \mathcal{E}\right\vert _{D}}$. In other words, such requirement would prohibit infinitesimal deformations of $A_{0}$ across gauge
orbits, which translates  into the vanishing of the next cohomology
group:%
\begin{equation*}
\mathbf{H}_{A_{0}\vert_D}^{1}=0.
\end{equation*}
This nondegeneracy will be central in the discussion of the gluing theory \cite{G2II}.
\end{remark}

The underlying heuristic in our definitions is the analogy
between looking for finite energy solutions to partial differential
equations over a compact manifold, with fixed values on a hypersurface,
say, and over a space with  cylindrical ends, under exponential decay to a suitable condition at infinity.

\subsection{\texorpdfstring
        {Short time existence of solutions and $C^{0}-$bounds}
        {Short time existence of solutions and C0-bounds}}
\label{subsec short time existence}
The short-time existence of a solution to our evolution equation is
a standard result:

\begin{proposition}     \label{Prop - short time existence}
Equation $\left( \ref{Eq heat}\right) $
admits a smooth solution $H_{S}\left( t\right) $, $t\in \left[ 0,\varepsilon
\right[ $, for $\varepsilon $ sufficiently small.
\begin{proof}
By the K\"{a}hler identities, $\left( \ref{Eq heat}\right) $ is
equivalent to the parabolic equation%
\begin{eqnarray*}       \index{evolution equation!is nonlinear parabolic}%
        \frac{\partial h}{\partial t} 
        &=&-\left\{ \Delta _{0}h+\mathbf{i}\left( \hat{F}_{0}.h
        +h.\hat{F}_{0}\right) +2\mathbf{i}\Lambda 
        \left( \bar{\partial}_{0}h.h^{-1}.\partial _{0}h\right) \right\}\\
        h\left( 0\right)  &=&I,\qquad \left. h\right\vert _{\partial W_{S}}=I.
\end{eqnarray*}%
for a positive self-adjoint endomorphism $h\left( t\right)
=H_{0}^{-1}H_{S}\left( t\right) $ of the bundle, with Hermitian metric $H_{0}$. Then the claim is an instance of the general theory  \cite[Part IV, \S 11, p.122]{Hamilton}.
\end{proof}
\end{proposition}

The task of extending solutions for all time is left to  \emph{Subsection \ref{subsec smooth sols for all time}}; let us first collect some preliminary results. We begin by
recalling the parabolic maximum principle:
\begin{lemma}[Maximum principle]        
\label{Lemma max principle mfds w bdry}
Let $X$ be a compact Riemannian
manifold with boundary and suppose $f\in C^{\infty }\left( \mathbb{R}%
_{t}^{+}\times X\right) $ is a nonnegative function satisfying:%
\begin{equation*}       \index{maximum principle}%
        \left(\frac{d}{dt}+\Delta \right) f_{t}\left( x\right) \leq 0,
        \quad \forall \left(t,x\right) \in \mathbb{R}_{t}^{+}\times X
\end{equation*}
and the Dirichlet condition:%
\begin{equation*}
        \left. f_{t}\right\vert _{\partial X}=0.%
\end{equation*}%
Then either $\sup\limits_{X}f_{t}$ is a decreasing function of $t$ or $%
f\equiv 0$.
\end{lemma}

The crucial role of the K\"{a}hler
structure in this type of problem is that it often suffices to control $\sup
\vert \hat{F}_{H}\vert $ in order to obtain uniform bounds on $H$
and its derivatives, hence to take limits in one-parameter
families of solutions. Let us then establish such a bound; I denote generally $\hat{e}:= \vert \hat{F}_{H}\vert _{H}^{2}$ 
and, in the immediate sequel, $\hat{e}_{t}:= \vert \hat{F}%
_{H_{S}\left( t\right) }\vert _{H_{S}\left( t\right) }^{2}$.

\begin{corollary}       \label{Cor Bounded F hat}

Let $\left\{ H_{S}\left( t\right) \right\} _{0\leq t<T}$ be a smooth solution to $\left( \ref{Eq heat}\right)$ on $W_{S}$. Then $\sup\limits_{W_{S}}\hat{e}_{t}$ is non-increasing with $t$; in fact, there exists $B>0$, independent of $S$ and $T$, such that%
\begin{equation}        \label{Eq bound on F hat}
        \sup_{W_{S}}\vert\hat{F}_{H_{S}(t)}\vert^{2}\leq B.
\end{equation}
\begin{proof}
Using the Weitzenb\"{o}ck formula \cite[p. 221]{4-manifolds}, one finds%
\begin{equation}       \label{eq (d/dt+D)e^ < 0}
        \left( \frac{d}{dt}+\Delta \right) \hat{e_{t}}
        =-\vert d_{H_{S}\left(t\right) }^{\ast }F_{H_{S}
        \left( t\right) }\vert ^{2}\leq 0.
\end{equation}%
At the boundary $\partial W_{S}$, for $t>0$, the Dirichlet condition $\left( \ref{Dirichlet}\right) $ means precisely that $\left. H_{S}\right\vert
_{\partial W_{S}}$ is constant, hence the evolution equation gives $\left. \hat{e}_{t}\right\vert _{\partial W_{S}}=\left. \vert H_{S}^{-1}\dot{H}%
_{S}\vert ^{2}\right\vert _{\partial W_{S}}\equiv 0$. Then $%
B:=\sup\limits_{W_{S}}\hat{e}_{0}$.
\end{proof}
\end{corollary}

In order to obtain $C^{0}-$bounds and state our first convergence result,
let us digress briefly into two ways of measuring metrics, which will be
convenient at different stages. First, given two
metrics $H$ and $K$ of same determinant we write%
\begin{equation*}
H=K.\e{\xi }
\end{equation*}%
where $\xi \in \Gamma \left( \func{End}\mathcal{E}\right) $ is traceless
and self-adjoint with respect both to $H$ and $K$, and define 
\begin{equation}          \label{eq def lambda bar}
        \bar{\lambda}:\func{Dom}\left( \xi \right) 
        \subseteq W\rightarrow \mathbb{R}_{\geq 0}
\end{equation}%
as the highest pointwise eigenvalue of $\xi $.
 This is a nonnegative, Lipschitz function on $W$, and a transversality argument shows that it is, in fact, smooth away from a set of real codimension 3 \cite[pp.240, 244]{inf dets stab bdls...}. Now, clearly
\begin{equation}        \index{distance}
                        \label{Eq C0 norm and lambda bar}
        \left\vert H-K\right\vert \leq (cst).\left\vert K\right\vert 
        .\left( \e{\bar{\lambda}}-1\right) ,
\end{equation}%
so it is enough to control $\sup \bar{\lambda}$ to get a bound on $%
\left\Vert H\right\Vert _{C^{0}}$ relatively to $K$. Except where otherwise
stated, we will assume $K=H_{0}$ to be the reference metric.

\begin{remark}       \label{Rem Co norm is complete}
The space of continuous (bounded) bundle
metrics is \emph{complete} with respect to the $C^{0}-$norm \cite[Theorem
7.15]{Rudin}, so  the uniform limit of a family of metrics is itself a well-defined metric.
\end{remark}

Second, there is an alternative notion of `distance' \cite[Definition 12]{ASDYM},
which is more natural to our evolution equation, as will become clear in the next few results:
\begin{definition}
Given any two Hermitian metrics $H,K$ on a complex vector bundle $\mathcal{E}$,
let%
\begin{eqnarray*}       \index{distance}%
        \tau \left( H,K\right)   &=& \func{tr}H^{-1}K \\
        \sigma \left( H,K\right) &=&\tau \left( H,K\right) +\tau \left( K,H\right)
        -2\func{rk}\mathcal{E}.
\end{eqnarray*}%
The function $\sigma $ is symmetric, nonnegative (since $a+a^{-1}\geq 2,$ $\forall
a\geq 0$) and it vanishes if and only if $H=K$.
\end{definition}

Although $\sigma $ is not strictly speaking a distance, it does provide an equivalent criterion for $C^0-$convergence of metrics with fixed determinant, based on the following estimate:%
\begin{eqnarray*}
        \sigma(H,K) &=&\func{tr}\left( \e{\xi }+\e{-\xi }-2.\func{Id}\right) \\
        &\geq &\e{\bar{\lambda}}+\e{-\bar{\lambda}}-2
        =\e{-\bar{\lambda}}\left( \e{\bar{\lambda}}-1\right) ^{2}.
\end{eqnarray*}
\begin{remark}        \label{Rem conv sigma => conv lambda} 
A sequence $\left\{H_{i}\right\} $ of metrics with fixed determinant converges to $H$ in $C^{0}$ if, and only if, $\sup \sigma\left( H_{i},H\right) \rightarrow 0$. The former because $\sup \e{\bar{\lambda}}\rightarrow 1$ and the latter because obviously $$\lim\sup \left\vert H_i^{-1}H-I\right\vert=\lim \sup \left\vert H^{-1}H_i-I\right\vert=0.$$
\end{remark}
Indeed, $\sigma$ compares
to $\bar{\lambda}$ by increasing functions; from the previous inequality we deduce, in particular,
\begin{equation}        \label{eq sigma bounds lambda}
        \e{\bar{\lambda}}\leq \sigma +2.
\end{equation}%
Conversely, it is easy to see that%
\begin{equation*}
        \sigma \leq 2r.\e{\left( r-1\right) \bar{\lambda}}
        \qquad\text{with}\quad  r:=\func{rk}\mathcal{E}.
\end{equation*}%

Furthermore, in the context of our
evolution problem, $\sigma $ lends itself to applications of the maximum
principle \cite[Proposition 13]{ASDYM}:

\begin{lemma}         \label{Lemma max prin for sigma}
If $H_{1}(t) $ and $H_{2}(t)$ are solutions of evolution equation $\left( \ref{Eq heat}\right) $, then $\sigma (t) =\sigma \left( H_{1}(t),H_{2}(t) \right)$ satisfies%
\begin{equation*}
        \left( \frac{d}{dt}+\Delta \right) \sigma \leq 0.
\end{equation*}
\end{lemma}

Combining \emph{Lemma \ref{Lemma max principle mfds w bdry}} and \emph{Lemma \ref{Lemma max prin for sigma}}, we  obtain the following straightforward consequences:
\begin{corollary}       \label{cor uniqueness}
Any two solutions to (\ref{Eq heat}) on $W_{S}$, with Dirichlet boundary
conditions (\ref{Dirichlet}), which are defined for $t\in[0,T[$, coincide for all $t\in[0,T[$. 
\end{corollary}
\begin{corollary}       \label{Cor H(t) -> H(T) in C0}
If a smooth solution $H_{S}\left( t\right) $ to $\left( \ref{Eq heat}\right) $ on $W_{S}$, with Dirichlet boundary
conditions (\ref{Dirichlet}), is defined for $t\in[0,T[$, then 
$H_{S}\left( t\right) 
\underset{t\to T}{\overset{C^{0}}{\longrightarrow }}H_{S}\left( T\right)$ and $H_{S}\left( T\right) $ is continuous.
\begin{proof}
The argument is analogous to \cite[Corollary 15]{ASDYM}. As discussed in \emph{Remarks \ref{Rem Co norm is complete} }and \emph{\ref{Rem conv sigma => conv lambda}}%
, it suffices to show that $\sup \sigma \left( H_{S}\left( t\right)
,H_{S}\left( t^{\prime }\right) \right)\to0 $ when $t^{\prime
}>t\rightarrow T$. Clearly 
$$
f_{t}:=\sup\limits_{W_{S}}
        \sigma \left(H_{S}\left( t\right) ,H_{S}\left( t+\tau \right) \right)
$$ 
satisfies the (boundary) conditions of \emph{Lemma \ref{Lemma max principle mfds w bdry}},
 so it is decreasing and 
\begin{equation*}
        \sup_{W_{S}}\sigma\left(H_{S}\left(t\right),
        H_{S}\left(t+\tau\right)\right)
        <\sup_{W_{S}}\sigma\left( H_{S}\left( 0\right) ,
        H_{S}\left( \tau \right)\right)
\end{equation*}%
for all $t,\tau ,\delta >0$ such that $0<T-\delta <t<t+\tau <T$. Taking $%
\delta <\varepsilon $ in \emph{Proposition \ref{Prop - short time existence}}
we ensure continuity of $H_{S}\left( t\right) $ at $t=0$, so the right-hand
side is arbitrarily small for all $t$ sufficiently close to $T$. Hence the family $\left\{H_{S}\left( t\right)\right\}_{t\in[0,T[} $ is uniformly Cauchy as $t\rightarrow T$ and the limit is continuous.%
\end{proof}
\end{corollary}
 
Looking back at $\left( \ref{Eq heat}\right) $, we may interpret $\hat{F}%
_{H} $ intuitively as a velocity vector along $1-$parameter families $%
H\left( t\right) $ in the space of Hermitian metrics. In this case, \emph{%
Corollary \ref{Cor Bounded F hat}} suggests an absolute bound on the
variation of $H$ for finite (possibly small) time intervals $0\leq t\leq T$ where solutions exist. A straightforward calculation yields
\begin{displaymath}
\begin{array}{cclcl}
        \dot{\sigma}&=&\displaystyle\frac{d}{dt}\sigma\left(H_S(t),H_0 \right)
        &\leq& \func{tr}\left[ \left( H_S^{-1}\dot{H}_S\right)
        .\left( \e{\xi}-\e{-\xi }\right) \right]\\
        &\leq& 2\left\vert \func{tr} \left( \mathbf{i}\hat{F}_{H_{S}}\right).
        \left( \e{\xi }-\e{-\xi }\right) \right\vert 
        &\leq& (cst.)\vert \hat{F}_{H_{S}}\vert.\e{\bar{\lambda}}\\
        &\leq &\tilde{B}.\e{\bar{\lambda}}&&\\
\end{array}
\end{displaymath}
with $\tilde{B}:=(cst.)\sqrt{B}$, using the evolution equation and \emph{Corollary \ref{Cor Bounded F hat}}.
Combining with $\left( \ref{eq sigma bounds lambda}\right) $ and integrating,%
\begin{equation}        \label{eq bound on sigma}
\fbox{$%
\begin{array}{c}
        \e{\bar{\lambda}}\leq \sigma 
        +2\leq 2\e{\tilde{B}T}:= C_{T},\quad \forall t\leq T. 
\end{array}%
$}  
\end{equation}%
Consequently, for any fixed $S_{0}>0$, the restriction of $H_{S}\left( t\right) $ to $W_{S_{0}}$ lies in a $C^{0}-$ball of radius $\log C_{T}$ about $H_{0}$ in the space of Hermitian metrics, for all $S\geq S_{0}$ and $t\leq T$. Since $C_{T}$ doesn't depend on $S$, the next \emph{Lemma} shows that the $H_{S}$ converge uniformly on compact subsets $W_{S_{0}}\subset W$ for any fixed interval $\left[ 0,T\right] $ (possibly trivial) where solutions exist for all $S>S_{0}$:        
\newpage
\begin{lemma}           \label{Lemma smooth sols on Ws for finite time}
If there exist $C_{T}>0$ and $S_{0}>0$ such that, for all $S^{\prime }>S\geq S_{0}$, the evolution equation
\begin{equation*}
\left\{ 
\begin{array}{l}
        H^{-1}\dfrac{\partial H}{\partial t}=-2\mathbf{i}\hat{F}_{H} \\ 
        \multicolumn{1}{c}{H\left( 0\right) =H_{0},
        \quad \left. H\right\vert_{\partial W_{S}}
        =\left. H_{0}\right\vert _{\partial W_{S}}}%
\end{array}%
\right. \quad \text{on } \quad W_{S}\times \left[ 0,T\right]
\end{equation*}%
admits a smooth solution $H_{S}$ satisfying%
\begin{equation*}
        \left.\sigma \left( H_{S},H_{S^{\prime }}\right) \right\vert _{W_{S}}
        \leq C_{T},
\end{equation*}%
then the $H_{S}$ converge uniformly to a (continuous) family $H$ defined on $W_{S_{0}}\times \left[ 0,T\right] $.
\begin{proof}
It is of course possible to find a function $\phi :W\rightarrow \mathbb{R}$
such that%
\begin{equation*}
\left\{ 
\begin{tabular}{l}
        $\phi \equiv 0$, on $W_{0}$ \\ 
        $\phi \left( y,\alpha ,s\right) =s$, for $s\geq 1$ \\ 
        $\left\vert \Delta \phi \right\vert \leq L$%
\end{tabular}%
\right. 
\end{equation*}%
thus giving an exhaustion of $W$ by our compact manifolds with boundary $W_{S}\simeq
\left\{ p\in W\mid \phi \left( p\right) \leq S\right\} $, $S\geq S_{0}$.
Taking $S_{0}<S<S^{\prime }$, I claim 
\begin{equation*}
        \sigma \left( \left. H_{S}\left( t\right) \right\vert _{W_{S}},
        \left. H_{S^{\prime }}\left(t\right) \right\vert _{W_{S}}\right)
         \left( p\right) \leq \frac{C_{T}}{S}
         \left(\phi \left( p\right) +Lt\right) ,\quad 
         \forall \left( p,t\right) \in W_{S}\times \left[ 0,T\right] ,
\end{equation*}%
which yields the statement, since its restriction to $W_{S_{0}}$ gives 
\begin{equation*}
        \left.\sigma\left(H_{S},H_{S^{\prime}}\right)\right\vert_{W_{S_{0}}}
        \leq \frac{C_{T}\left( S_{0}+LT\right) }{S}
        \underset{S\rightarrow \infty }{\longrightarrow }0.
\end{equation*}%
The inequality holds trivially at $t=0$ and on $\partial W_{S}$ by $\left( \ref{eq bound on sigma}\right) $, hence on the whole of $W_{S}\times \left[ 0,T%
\right] $, by the maximum principle [\emph{Lemma \ref{Lemma max principle
mfds w bdry}}]:%
\begin{equation*}
\left( \frac{d}{dt}+\Delta \right) \left( \sigma \left( H_{S},H_{S^{\prime
}}\right) -\frac{C_{T}}{S}\left( \phi +Lt\right) \right) \leq -\frac{C_{T}}{S%
}\left( \Delta \phi +L\right) \leq 0,
\end{equation*}%
using $\left\vert \Delta \phi \right\vert \leq L$ and \emph{Lemma \ref{Lemma
max prin for sigma}}.
\end{proof}
\end{lemma}

This defines a 1-parameter family of \emph{continuous} Hermitian bundle metrics over  $W$:%
\begin{equation}        \label{Eq limit metric H(t)}
\fbox{$%
\begin{array}{c}
        H\left( t\right) := \lim\limits_{S\rightarrow \infty }
        H_{S}\left(t\right) ,\quad t\leq T. 
\end{array}%
$}  
\end{equation}%
It remains to show that any $H_{S}$ can be smoothly extended for all $t\in \lbrack 0,\infty \lbrack $ and that the limit $H\left( t\right) $ is itself a smooth solution of the evolution equation on $W$, with satisfactory asymptotic properties along the tubular end.

\subsection{Smooth solutions for all time}
\label{subsec smooth sols for all time}

Following a standard procedure, we will now check that the bound $\left( \ref{Eq bound on F hat}\right) $ allows us to smoothly extend
solutions $H_{S}$ up to  $t=T$, hence past $T$, for all time \cite[Corollary 6.5]{Simpson}. More precisely, we can exploit the features of our problem to control a Sobolev norm $\left\Vert \Delta H\right\Vert $ by $\sup
\vert \hat{F}\vert $ and weaker norms of $H$.
The fact that one still has `Gaussian' bounds on the norm of the heat kernel in our noncompact case will be central to the argument [\emph{Theorem \ref{Thm Gaussian bound on heat kernel}}, Appendix \ref{App Gaussian upper bound}].
\begin{lemma}
Let $H$ and $K$  be smooth Hermitian metrics on a holomorphic bundle over a K\"ahler manifold with K\"ahler form $\omega$; then, for any submultiplicative pointwise norm $\left\Vert .\right\Vert$,
\begin{equation}        \label{eq bound on Laplacian}
\fbox{$%
\begin{array}{c}
        \left\Vert\Delta_K H \right\Vert 
        \leq (cst.)\left[ \left( \Vert\hat{F}_H\Vert +1 \right)
        \left\Vert H \right\Vert + \left\Vert \nabla_{K}H\right\Vert ^{2}
        \Vert H^{-1}\Vert \right]  
\end{array}%
$}
\end{equation}
where $\Delta_K:= 2\mathbf{i}\Lambda_\omega\bar{\partial}\partial_K$ is the K\"ahler Laplacian and $(cst.)$ depends on $K$ and $\left\Vert .\right\Vert$ only.
\begin{proof}
Write  $h=K^{-1}H$ and $\nabla_K$ for the Chern connection of $K$. Since $\nabla_{K}K=0$,
\begin{displaymath}
        \Delta_K H=K.\Delta_K h.
\end{displaymath}
On the other hand, the Laplacian satisfies \cite[p.46]{4-manifolds}\cite[p.15]{ASDYM}
\begin{displaymath}
        \Delta_Kh =  h \left( \hat{F}_{H}-\hat{F}_{K} \right)+\mathbf{i}\Lambda_{\omega}         \left( \bar{\partial}h.h^{-1}\wedge\partial _{K}h\right)
\end{displaymath}
so the triangular inequality and again $\nabla_{K}K=0$ yield the result.
\end{proof}
\end{lemma}
That will be the key to the recurrence argument behind \emph{Corollary \ref%
{Cor limit metric is smooth}}, establishing smoothness of $H_{S}$ as $%
t\rightarrow T$. We will need the following
technical facts, adapted from   \cite[Lemma 6.4]{Simpson}.
\begin{lemma}           \label{Lemma Hi is in Lp2}
Let $\left\{ H_{i}\right\} _{0\leq i< I}$ be a
one-parameter family of Hermitian metrics on a bundle $\mathcal{E}%
\rightarrow X$ over a compact K\"{a}hler manifold with boundary such that
\begin{enumerate}
        \item 
        $H_{i}\underset{i\to I}{\overset{C^{0}(X)}{\longrightarrow }}H_{I}$, where $H_{I}$ is a
continuous metric,
        \item 
        $\sup\limits_{X}\vert \hat{F}_{H_{i}}\vert $ is bounded
uniformly in $i$,
        \item 
        $\left. H_{i}\right\vert _{\partial X}=H_{0}$;
\end{enumerate}
then the family $\left\{ H_{i}\right\} $ is bounded in $L_{2}^{p}\left( X\right) $,
 for all $1\leq p<\infty $, so $H_{I}$ is of class $C^{1}$.
\end{lemma}

\begin{corollary}       \label{Cor Hs(T) is in Lp2}
If $\left\{ H_{S}\left( t\right)\right\} _{0\leq t<T}$ is a solution of $\left( \ref{Eq heat}\right) $ with Dirichlet condition (\ref{Dirichlet}) on $\partial W_{S}$, then the $H_{S}\left( t\right) $ are bounded in $L_{2}^{p}\left( W_{S}\right) $ uniformly in $t$, for all $1\leq p<\infty $, and $H_{S}\left( T\right)$ is of class $C^{1}$.
\begin{proof}
By \emph{Corollary \ref{Cor H(t) -> H(T) in C0} }and $\left( \ref{Eq bound
on F hat}\right) $, $\left\{ H_{S}\left( t\right) \right\} _{0\leq t<T}$
satisfies respectively (1) and (2) in  \emph{Lemma \ref{Lemma Hi is in Lp2}}.
\end{proof}
\end{corollary}

The \emph{Corollary} gives, in particular, a time-uniform bound on $\left\Vert
F_{H_{S}}\right\Vert _{L^{p}\left( W_{S}\right) }$. This can actually be
improved to a uniform bound on all derivatives of curvature:

\begin{lemma}   \label{Lemma all derivatives of F are unif bounded}
$F_{H_{S}}$ is bounded in $C^k\left( W_{S}\right) $, uniformly in $t\in[0,T] $, for each $k\geq 0$.
\begin{proof}
By induction on $k$:

\underline{$k=0$}: following \cite[Lemma 18]{ASDYM}, we obtain a uniform
bound on $e_S(t):= \left\vert F_{H_{S}\left( t\right)
}\right\vert ^{2}$, using the fact that 
\begin{equation*}
\left( \frac{d}{dt}+\Delta \right) e_{S}\leq \left( cst.\right) \left( (e_{S})^{\frac{3%
}{2}}+e_{S}\right)
\end{equation*}%
\cite[Proposition 16, (ii)]{ASDYM}, and consequently%
\begin{equation}
e_{S}\left( t\right) \leq \left( cst.\right) \left( 1+\int_{0}^{t}\left\Vert
K_{t-\tau }\right\Vert _{L^{p}\left( W_{S}\right) }\Vert \left(e_S\right)^{\frac{3}{2}%
}+e_S\Vert _{L^{q}\left( W_{S}\right) }\right) ,
\label{Eq heat bound on e}
\end{equation}%
where $K_{t}$ is the heat kernel associated to $\frac{d}{dt}+\Delta $ and $%
\frac{1}{p}+\frac{1}{q}=1$. On the complete [\emph{Theorem \ref{Thm
SU(3) mfds}}] $6-$dimensional Riemannian manifold $W_{S}$, $K_{t}$ satisfies 
\cite[\S 9]{Eels and Sampson} the diagonal condition%
\begin{equation*}
K_{t}\left( x,x\right) \leq \frac{\left( cst.\right) }{t^{3}},\qquad \forall
x\in W_{S}
\end{equation*}%
of \emph{Theorem \ref{Thm Gaussian bound on heat kernel}} in Appendix \ref{App Gaussian upper bound}, which gives a `Gaussian' bound on the heat kernel. So, fixing $C>4$ and denoting $%
r\left( .,.\right) $ the geodesic distance, we have 
\begin{equation*}
K_{t}\left( x,y\right) \leq \frac{\left( cst.\right) }{t^{3}}\exp \left\{ -%
\frac{r\left( x,y\right) ^{2}}{Ct}\right\} ,\qquad \forall x,y\in W_{S}.
\end{equation*}%
Hence, for each $x\in W_{S}$, we obtain the bound%
\begin{eqnarray*}
\left\Vert K_{t}\left( x,.\right) \right\Vert _{L^{p}\left( W_{S}\right) }
&\leq &\frac{\left( cst.\right) }{t^{3}}\left( \int_{W_{S}}\exp \left\{ -p%
\frac{r\left( x,y\right) ^{2}}{Ct}\right\} dy\right) ^{\frac{1}{p}} \\
&\leq &\frac{\left( cst.\right) }{t^{3}}\left( \int_{0}^{\infty }\left( 
\tfrac{Ct}{p}\right) ^{3}u^{5}\e{-u^{2}}du\right) ^{\frac{1}{p}} \\
&\leq &\tilde{c}_{p}\;t^{\frac{3}{p}\left( 1-p\right) }.
\end{eqnarray*}%
Now, $p<\frac{3}{2}\Leftrightarrow \frac{3}{p}\left( 1-p\right) >-1$, in
which case%
\begin{equation*}
\int_{0}^{T}\left\Vert K_{t}\left( x,.\right) \right\Vert _{L^{p}\left(
W_{S}\right) }dt\leq c_{p}\left( T\right) .
\end{equation*}%
Inequality $\left( \ref{Eq heat bound on e}\right) $ gives the desired
result provided $\left(e_S\right)^{\frac{3}{2}}\in L^{q}\left( W_{S}\right) $ for some $q>3$%
; this means $F_{H_{S}\left( t\right) }\in L^{\tilde{q}}\left( W_{S}\right) $
for some $\tilde{q}>9$, which is guaranteed by \emph{Corollary \ref{Cor
Hs(T) is in Lp2}}.

\underline{$k\Rightarrow k+1$}: The general recurrence step is identical to 
\cite[Corollary 17 (ii)]{ASDYM}, using the maximum principle [\emph{Lemma \ref%
{Lemma max prin for sigma}}] with boundary conditions.
\end{proof}
\end{lemma}

We are now in shape to put into use the K\"{a}hler setting, combining the $C^k-$bounds on $F$ (hence on $\hat{F}$) with inequality (\ref{eq bound on Laplacian}), via elliptic regularity:

\begin{lemma}        \label{Lemma smoothness of limit}
For $0\leq I\leq \infty$, let $\left\{ H_{i}\right\}_{0\leq i< I} $ be a
one-parameter family of Hermitian metrics on a holomorphic vector bundle $\mathcal{E}\rightarrow X$ over a compact (real) $2n-$dimensional K\"{a}hler manifold with boundary
such that
\begin{enumerate}
        \item
        $\left\{ H_{i}\right\} $ is bounded in $L_{2}^{p}\left( X\right), $ for all $1\leq p<\infty $, and $\displaystyle H_{i}\underset{i\rightarrow I}{\overset{L_{2}^{p}(X)}{\longrightarrow }}H$, where $H$ is a continuous metric,        
        \item
        $\{\hat{F}_{H_{i}}\}$is bounded in $L_{k}^{p}\left( X\right) $, for each $(p,k)\in [1,\infty[\times\mathbb{N}$, 
        \item
        $\{H_{i}\}$is bounded in $L_{k}^{p}\left( \partial X\right) $, for each $(p,k)\in [1,\infty[\times\mathbb{N}$; 
        
\end{enumerate}
then $\displaystyle H_{i}\underset{i\rightarrow I}{\overset{C^\infty(X)}{\longrightarrow }}H$ and $H$ is smooth.
\begin{proof}
Fixing $p>2n$, I will prove the following statement by induction in $k$:%
\begin{equation*}
\begin{tabular}{l}
$\left\Vert H_i\right\Vert _{L_{k+2}^{p}(X)}$ and $\Vert H_i^{-1}\Vert
_{L_{k}^{p}(X)}$ are bounded uniformly in $i$, $\forall k\geq0$.%
\end{tabular}%
\end{equation*}%
The first hypothesis gives step $k=0$, as well as $\Vert H_i^{-1}\Vert
_{L_{1}^{p}(X)}<\infty$, since the Sobolev embedding implies $H_i$ is in $C^1$. Now, assuming the statement up
to step $k-1\geq0$ implies in particular $k\geq1>\frac{2n}{p}$, which authorises the multiplication $L^p_k\times L^p_{k-1}\to L^p_{k-1}$ and so 
\begin{eqnarray*}
        \Vert H_i^{-1}\Vert _{L_{k}^{p}}
        &=&\Vert H_i^{-1}\Vert_{L^{p}}
        +\Vert \nabla \left( H_i^{-1}\right)\Vert_{L_{k-1}^{p}}\\
        &\leq &\Vert H_i^{-1}\Vert _{L^{p}}
        +\Vert H_i^{-1}\Vert_{L_{k-1}^{p}}^{2}
        \Vert \nabla H_i\Vert _{L_{k}^{p}} \\
        &\leq &\Vert H_i^{-1}\Vert _{L_{k-1}^{p}}
        \left( 1+\Vert H_i^{-1}\Vert _{L_{k-1}^{p}}
        \left\Vert H_i\right\Vert _{L_{k+1}^{p}}\right)
\end{eqnarray*}%
so $\Vert H_i^{-1}\Vert _{L_{k}^{p}}$ is bounded. On the other hand,
elliptic regularity on manifolds with boundary and $\left( \ref{eq
bound on Laplacian}\right) $, with $K=H_0$, give%
$$
\begin{array}{rl}
\left\Vert H_i\right\Vert _{L_{k+2}^{p}}^{2} \leq &c_0
        \left(\left\Vert \Delta H_i\right\Vert _{L_{k}^{p}}^{2}
        +\left\Vert H_i\right\Vert_{L_{k+1}^{p}}^{2}
        +\left\Vert H_i\right\Vert _{L_{k+\frac{3}{2}}^{p}
        \left(\partial X\right) }^{2}\right) \\
        \leq &c_0
        \left[ \left\Vert H_i\right\Vert _{L_{k+1}^{p}}^{2}
        \left(1+\left( 1+\Vert \hat{F}_{H_i}\Vert _{L_{k}^{p}}
        +\left\Vert H_i\right\Vert _{L_{k+1}^{p}}
        \Vert H_i^{-1}\Vert _{L_{k}^{p}}\right)^{2}\right)+
        \left\Vert H_i\right\Vert 
        _{L_{k+\frac{3}{2}}^{p}\left( \partial X\right) }^{2} \right]\end{array}
$$
where $c_0$ depends on $H_0$ and $X$ only, and all those terms are bounded by assumption. 

Since $p$ could be chosen arbitrarily big, the family $\{H_{i}\}$ is uniformly bounded in each $C^r$. But we know it converges to $H$ in $C^1$, hence in fact it converges in $C^\infty$ and the limit is
smooth.
\end{proof}
\end{lemma}

\begin{corollary}
\label{Cor limit metric is smooth}Under the Dirichlet conditions (\ref{Dirichlet}), the limit
metric $H_{S}\left( T\right) $ is smooth.

\begin{proof}
We apply \emph{Lemma \ref{Lemma smoothness of limit}} to the one-parameter family $\{H_S(t)\}_{0\leq t\leq T}$ of Hermitian metrics on the restriction $\mc{E}\to W_S$ given by \emph{Lemma \ref{Lemma smooth sols on Ws for finite time}}. Then 
\emph{Corollary \ref{Cor Hs(T) is in Lp2}} gives hypothesis (1), \emph{%
Lemma \ref{Lemma all derivatives of F are unif bounded}} gives (2) and the
Dirichlet condition on $\partial W_{S}$ gives (3), as $H_{0}$ is smooth.
\end{proof}
\end{corollary}

Since $H_{S}\left( t\right) \overset{C^{\infty }(W_S)}{\underset{t\rightarrow T}{\longrightarrow} }%
H_{S}\left( T\right) $, the solution can be smoothly extended beyond $T$, by short-time existence, 
hence for all time \cite[Proposition 6.6]{Simpson}:

\begin{proposition}
\label{Prop smooth solution for arbitrary finite time}Given any $T>0$, the
family of Hermitian metrics $H\left( t\right) $ on $\mathcal{E}\rightarrow W$
defined by $\left( \ref{Eq limit metric H(t)}\right) $ is the unique, smooth
solution of the evolution equation$%
\index{evolution equation!on Wx[0,T]@on $W\times \left[ 0,T\right] $}$%
\begin{equation*}
\left\{ 
\begin{array}{c}
H^{-1}%
\dfrac{\partial H}{\partial t}=-2\mathbf{i}\hat{F}_{H} \\ 
H\left( 0\right) =H_{0}%
\end{array}%
\right. \quad \text{on }W\times \left[ 0,T\right] 
\end{equation*}%
with $\det H=\det H_0$ and $\sup_{W}\left\vert H\right\vert <\infty $. Furthermore, 
$$
\sup_{W}\vert \hat{F}_{H\left( t\right) }\vert \leq
B=\sup_{W}\vert \hat{F}_{H_{0}}\vert. 
$$
\begin{proof}
Using \emph{Lemma \ref{Lemma smoothness of limit}} on any compact subset $\Omega_{S_0}:=W_{S_{0}}\times \left[ 0,T\right] $, the $H_{S}$ are $%
C^{\infty }-$bounded (uniformly in $S$). By the evolution equation, 
$$
\frac{\partial H_{S}}{\partial t}
\underset{S\rightarrow \infty}{\overset{C^{\infty }\left( \Omega_{S_0} \right)}{\longrightarrow}}
\frac{\partial H}{\partial t}
$$
so $H$ is a solution on $\Omega_{S_0} $ satisfying the
same bounds, and this is independent of the choice of $S_{0}$.
The bound on $\sup_{W}\vert \hat{F}_{H\left( t\right) }\vert $ is immediate from \emph{Corollary \ref{Cor Bounded F hat}} and uniqueness is the statement of \emph{Corollary \ref{cor uniqueness}}.
\end{proof}
\end{proposition}

\subsection{Asymptotic behaviour of the solution}
\label{subsec asymp behaviour of sol} 
We have a solution $\{ H\left( t\right) \} $ of the flow on $W$ [\emph{Proposition \ref{Prop smooth solution for arbitrary
finite time}}], giving a Hermitian metric on $\mathcal{E}\rightarrow W$ for
each $t\in \left[ 0,T\right] $. Let us study the asymptotic properties of $%
H\left( t\right) $ along the non-compact end. Set
\begin{equation*}
        \hat{e}_{t}=\vert \hat{F}_{H\left( t\right) }\vert ^{2}.
\end{equation*}%
First of all, as a direct consequence of \emph{Lemma \ref{lemma decay of
coordinates dz dz bar},} I claim%
\begin{equation}        \label{Eq (e hat)0 decays exp}
                        \index{asymptotic decay} 
\hat{e}_{0}\leq B\epsilon ,\qquad \epsilon := \left\{ 
\begin{tabular}{ll}
$1$ & on $W_{0}$ \\ 
$\e{-s}$ & on $\partial W_{s}$,$\quad s\geq 0$%
\end{tabular}%
\right.
\end{equation}%
where $B=\sup_{W}\hat{e}_{0}$ [\emph{Corollary \ref{Cor Bounded
F hat}}]. In the trivialisation $\left. \mathcal{E}\right\vert _{U}\simeq\left\{ \left\vert z\right\vert <1 \right\}\times \left. \mathcal{E}\right\vert _{D}$
 over the
neighbourhood of infinity $U$, with coordinates $\left( z,\xi ^{1},\xi ^{2}\right)$ such that
  $D=\left\{ z=0\right\} \subset \bar{W}$, the curvature $F_{H_{0}}$ is represented by the endomorphism-valued $\left(
1,1\right) -$form:%

\begin{equation}       \label{Eq asymptotics of F_H0}
\left.F_{H_{0}}\right\vert_{U\smallsetminus D}=F_{zz}\:\underset{O( \left\vert
z\right\vert ^{2}) }{\underbrace{dz\wedge d\bar{z}}}+\sum_{i} (F_{zi}\;%
\underset{O\left( \left\vert z\right\vert \right) }{\underbrace{dz\wedge d\bar{\xi}^{i}}}+F_{iz}\;\underset{O\left( \left\vert z\right\vert \right) }{%
\underbrace{d\xi ^{i}\wedge d\bar{z}}})+\sum_{i,j} F_{ij}\; d\xi ^{i}\wedge d\bar{\xi%
}^{j}.
\end{equation}%
The terms involving $dz$ or $d\bar{z}$ decay at least as $O\left( \left\vert
z\right\vert \right) $ along the tubular end [\emph{Lemma \ref%
{lemma decay of coordinates dz dz bar}}], and all the coefficients of ${F}_{H_{0}}$ are bounded, so $\left. F_{H_{0}}\underset{%
\left\vert z\right\vert \rightarrow 0}{\longrightarrow }\sum F_{ij}d\xi
^{i}\wedge d\xi ^{j}\right. $. Consequently,
\begin{equation*}
        \hat{F}_{H_{0}}\left( z,\xi ^{1},\xi ^{2}\right) 
        \underset{\left\vert z\right\vert \rightarrow 0}{\longrightarrow}
        \left.\hat{F}_{H_{0}}\right\vert _{D}
        \left(\xi ^{1},\xi ^{2}\right) 
        =0,
\end{equation*}%
i.e., $\hat{F}_{H_{0}}$ decays exponentially to zero as $s\rightarrow \infty 
$. From $\left( \ref{Eq (e hat)0 decays exp}\right) $ we now obtain the
 exponential decay of each $\hat{e}_{t}$ along the cylindrical end:

\begin{proposition}     \label{Prop e hat decays exponentially along S}
Take $B$ and $\epsilon$ as in (\ref{Eq (e hat)0 decays exp}); then
\begin{equation*}
        \hat{e}_t \leq \left(B\e{t}\right)\epsilon 
        \quad \text{on}\quad W.
\end{equation*}
\begin{proof}
The statement is obvious on $W_{0}$. For any $s_0,t_0\geq 0$, take $T=S>\max
\left\{ s_0,t_0\right\} $, let $\Sigma _{S}:= W_{S}\smallsetminus W_{0}$ and
consider on $\Sigma _{S}\times \left[ 0,T\right] $ the comparison function $%
g\left( t,s\right) := B\e{t-s}$. Using the Weitzenb\"ock formula one shows that  $\left( \frac{d}{dt}+\Delta \right) 
\hat{e}_{S}\leq 0$ [cf. $(\ref{eq (d/dt+D)e^ < 0})$], where $\hat{e}_{S}=\vert \hat{F}_{H_S} \vert^2$ and $H_S$ is a solution of our flow on $W_S$ as in $\emph{Lemma}$\emph{\ \ref{Lemma smooth sols on Ws for finite time}}. For $\psi := $ $\hat{e}%
_{S}-g$, one clearly has $\left( \frac{d}{dt}+\Delta \right) \psi \leq 0\ $(recall that our sign convention for the Laplacian is $\Delta
\overset{loc}{=}\mathbf{-}\sum\frac{\partial ^{2}}{\partial x_{i}^{2}}$) and, by the maximum principle [\emph{Lemma \ref{Lemma max principle mfds w bdry}}],%
\begin{equation*}
\psi \leq \max\limits_{\partial \left( \left[ 0,T\right] \times \Sigma
_{S}\right) }\left\{ \hat{e}_{S}-B\e{t-s}\right\} \leq 0.
\end{equation*}%
To see that the r.h.s. is zero, there are four boundary terms to check:
\begin{description}
        \item[\quad$s=S:$]
the Dirichlet condition means $\hat{e}_{S}\left(
t,S\right) =0$, $\forall t>0$, so $\psi \left( t,S\right) \leq 0$;

        \item[\quad$s=0:$]
$\psi \left( t,0\right) \leq B\left( 1-\e{t}\right)\leq 0$;

        \item[\quad$t=0:$]
$\left( \ref{Eq (e hat)0 decays exp}\right) $ gives $\psi
\left( 0,s\right) \leq 0$;

        \item[\quad$t=T:$]
again by \emph{Corollary \ref{Cor Bounded F hat} }we have 
$\psi \left( T,s\right) \leq B\left( 1-\e{T-s}\right) \leq 0$.
\end{description}
This shows that $\hat{e}_{S}\left( t,s\right) \leq B\e{t-s}$ on $\Sigma_S\times[0,T]$. Take $T=S\rightarrow \infty $.
\end{proof}
\end{proposition}

To conclude exponential $C^{0}-$%
convergence of $H\left( t\right) $ along the cylindrical end, recall from  $%
\left( \ref{eq bound on sigma}\right) $  that the constant  $\tilde{B}=(cst.)\sqrt{B}$ is obtained from the uniform bound $\hat{e} \leq B$. Now, in the context of \emph{Proposition \ref{Prop e hat decays exponentially along S}}, this control is improved to an exponentially decaying pointwise bound along the tube, thus we may replace $\tilde{B}\e{\frac{1}{2}(T-s)}$ for $\tilde{B}$ in that expression:%
\index{asymptotic decay!of $\hat{F}$}%
\begin{equation}        \label{eq exp decay C0}
        \left.\sigma\left(H\left( t\right) ,H_{0}\right) \right\vert 
        _{\partial W_{S}}
        \leq 2\left( \e{\tilde{B}T\e{\frac{1}{2}(T-S)}}-1 \right)
        =
        O\left(\e{-S}\right) .  
\end{equation}

The next result establishes exponential decay of $H\left(
t\right) $ in $C^{1}$, emulating the proof of \cite[Proposition 8]{Approx instantons}. I state it in rather general terms to highlight
the fact that essentially all one needs to control is the Laplacian, hence $\hat{F}$ in view of (\ref{eq bound on Laplacian}). 
\begin{proposition}     \label{Prop regularity on smaller open set}
Let $V$ be an open set of a
Riemannian manifold $X$, $V^{\prime }\Subset\ V$ an interior domain and $Q\rightarrow X$ some bundle with connection $\nabla $ and a continuous
fibrewise metric. There exist constants $\varepsilon ,A>0$ such that, if a smooth section $\phi \in \Gamma (Q) $ satisfies:

\begin{list}{\labelitemi}{\leftmargin=1em}
        \item[(1)] 
$\Vert \phi \Vert _{C^{0}\left( V\right) }\leq \varepsilon $;        
        \item[(2)] 
$\vert \Delta \phi \vert \leq f\left( \left\vert \nabla
\phi \right\vert \right) $ on $V$, for some non-decreasing function $f:%
\mathbb{R}^{+}\rightarrow \mathbb{R}^{+}$;        
        \item[(2')] 
assumption (2) remains valid under local rescalings, in the sense that, on every ball $B_{r}\subset V$, it still holds
for some function $\tilde{f}$ after the radial rescaling 
$\tilde{\phi}\left( \tilde{x}\right) := \phi
\left( mx\right) $, $m>0$;
\end{list}
then 
\begin{eqnarray*}
        \left\Vert \phi \right\Vert _{C^{1}\left( V^{\prime }\right) }
        &\leq& 
        A\left\Vert \phi \right\Vert _{C^{0}\left( V\right) }.  
\end{eqnarray*}

\begin{proof}
I first contend that $\phi $ obeys an \emph{a priori} bound%
\begin{equation*}
        \left\vert \left( \nabla \phi \right) _{x}\right\vert r\left( x\right)         \leq
        1,\quad \forall x\in V
\end{equation*}%
where $r\left( x\right) :V\rightarrow \mathbb{R}$ is the distance to $\partial V$. Since the term on the left-hand side is zero on $\partial V$,
its supremum is attained at some $\hat{x}\in V$ (possibly not unique).
Write \begin{displaymath}
        m:= \left\vert \left( \nabla \phi \right) _{\hat{x}}\right\vert,
        \quad 
        R=r\left( \hat{x}\right) 
\end{displaymath}
and suppose, for contradiction, that $R>\frac{1}{m}$. If that is the case, then we rescale the ball $B_{R}\left( \hat{x}\right) $ by the factor $m$, obtaining a rescaled local section $\tilde{\phi}$ defined in $\tilde{B}_{mR}\supset \tilde{B}_{1}$. In
this picture, any point in $\tilde{B}:= \tilde{B}_{\frac{1}{2}}$ is
further from $\partial V$ than $\frac{R}{2}$, hence, by definition of $\hat{x}$, $\Vert \nabla \tilde{\phi}\Vert _{C^{0}\left( \tilde{B}\right)}\leq 2$. 
By (2) and (2'), there exists $L>0$ such
that 
$\Vert \Delta \tilde{\phi} \Vert _{C^{0}\left( \tilde{B}\right)}\leq L$,
and elliptic regularity gives%
\begin{equation*}
        \Vert \nabla \tilde{\phi}\Vert _{C^{0,\alpha }
        \left( \tilde{B}\right) }\leq c_{\alpha }
        .\left( L+\varepsilon \right) 
        := \tilde{c}_{\alpha }
\end{equation*}%
using assumption (1). Now, the rescaled gradient at $\hat{x}$ has norm $\vert ( \nabla \tilde{\phi}) _{\hat{x}}\vert =1$ so,
taking $\alpha =\frac{1}{2}$ (say) in a smaller ball of radius 
$\rho =(\frac{1}{2\tilde{c}_{\frac{1}{2}}})^{^2}$, 
\begin{equation*}
        \vert ( \nabla \tilde{\phi}) _{x}\vert 
        \geq 
        1-\tilde{c}_{\frac{1}{2}}.\rho ^{\frac{1}{2}}
        \geq 
        \frac{1}{2},\quad \forall x\in \tilde{B}_{\rho }.
\end{equation*}%
This means $\vert \tilde{\phi}\vert $ varies by some definite $%
\delta >0$ inside $\tilde{B}_{\rho }$ and we reach a contradiction
choosing $\varepsilon <\delta $. So%
\begin{equation*}
        \left\vert \left( \nabla \phi \right) _{x}\right\vert 
        \leq 
        \left(\inf\limits_{\partial U}r\right) ^{-1},
        \quad \forall x\in U\Subset V
\end{equation*}%
for some open set $U\Supset V^{\prime }$. To conclude the proof, it
suffices to control the $L^{2}-$norm of $\nabla \phi $ on $U$: 
\begin{eqnarray*}
        \left( cst.\right) \left\Vert \nabla \phi \right\Vert^2 _{C^{0}\left(
V^{\prime }\right) } &\leq &\left\Vert \nabla \phi \right\Vert ^2_{L^{2}\left(
U\right) }=\int_{U}\left\langle \nabla \phi ,\nabla \phi \right\rangle 
=\int_{U}\left\langle
\phi ,\Delta \phi \right\rangle \\
&\leq &\left[f\left(\left\Vert \nabla \phi \right\Vert _{C^{0}\left( U\right)
}\right)\right]^{2}\left\Vert \phi \right\Vert _{L^{2}\left( V^{\prime }\right) }^{2}
\end{eqnarray*}%
and the last term is obviously bounded by $\left( cst.\right) \left\Vert
\phi \right\Vert _{C^{0}\left( V\right) }^{2}$.
\end{proof}
\end{proposition}

Now let $\func{End}\mathcal{E}=Q$ in \emph{Proposition \ref{Prop regularity on smaller open set}}, with connection $\nabla_0$ induced by $H_0$.

\begin{notation}        \label{Not finite cylinder}%
Given $S>r>0$, write $\Sigma _{r}\left( S\right) $ for the
interior of the cylinder $\left( W_{S+r}\smallsetminus
W_{S-r}\right) $ of `length' $2r$. We denote the \emph{$C^k-$exponential tubular limit} of an element in  $C^k\left(\Gamma(Q)\right)$ by:
\begin{eqnarray*}
        \phi \overset{C^k}{\underset{S\rightarrow\infty}{\longrightarrow}}\phi_{0}         & \dot{\Leftrightarrow}& \left\Vert \phi -\phi_{0} \right\Vert_{C^k\left( \Sigma_1(S),\omega \right)}=O\left( \e{-S} \right).
\end{eqnarray*}
\end{notation}

For $S\geq 3$, let $V=\Sigma _{3}\left( S\right) $ and $%
V^{\prime }=\Sigma _{2}\left( S\right) $ so that the distance of $V^{\prime }
$ to $\partial V$ is always $1$. In view of (\ref{eq exp decay C0}), for whatever $\varepsilon>0$  given by the statement, it is possible to choose $S\gg0$ so that $\phi=\left(\left.H\left( t\right) -H_{0}\right)\right\vert_{\Sigma_3(S)}$
satisfies the first condition (for arbitrary fixed $t$), hence also the
second one by (\ref{eq bound on Laplacian}), with  $f(x)=(cst.)\left[ \left(B+1\right)\varepsilon +x^2 \right]$ and $(cst.)$ depending only on $H_0$ and $\varepsilon$. We conclude, in particular, that $H(t)$ is  $C^1-$exponentially asymptotic to $H_0$ in the tubular limit:
\begin{equation}        \label{eq H decays in C1 as S -> infinity}
        H(t) \overset{C^1}{\underset{S\rightarrow\infty}{\longrightarrow}}H_{0}.
\end{equation}%

Furthermore, in our case the bound on the Laplacian $\left( \ref{eq bound on
Laplacian}\right) $ holds for any $L_{k}^{p}-$norm, given our control over
all derivatives of the curvature [\emph{Lemma \ref{Lemma all derivatives of F are unif bounded}}], so the argument above lends itself to the obvious iteration over shrinking tubular segments $\Sigma_{1+\frac{1}{k}}(S)$:
\begin{corollary}
\label{Cor H -> H0 exponentially in all derivatives}Let $\left\{ H\left(
t\right) \mid t\in \left[ 0,T\right] \right\} $ be the solution to the
evolution equation on $\mathcal{E}\rightarrow W$ given by \emph{Proposition %
\ref{Prop smooth solution for arbitrary finite time}}; then%
\index{asymptotic decay!of all derivatives of $H(t)$} 
\begin{equation*}
\boxed{H(t) \overset{C^k}{\underset{S\rightarrow\infty}{\longrightarrow}}H_{0},\quad \forall k\in \mathbb{N}.}
\end{equation*}
\end{corollary}
Combining existence and uniqueness of the solution for arbitrary time [\emph{Proposition \ref{Prop smooth solution for arbitrary finite
time}}] and $C^\infty-$exponential decay [\emph{Corollary \ref{Cor H -> H0 exponentially in all derivatives}}], one has the main statement:

\begin{theorem}
\label{thm sols for all time and exponentially decaying}Let $\mathcal{E}%
\rightarrow W$ be stable at infinity, with reference metric $H_{0}
$, over an asymptotically cylindrical $SU\left( 3\right) -$manifold $W$ as given by the \emph{Calabi-Yau-Tian-Kovalev Theorem \ref{Thm SU(3) mfds}}; then,
for any $0<T<\infty $, $\mathcal{E}$ admits a $1-$parameter family $\left\{
H_{t}\right\} $ of smooth Hermitian metrics solving $%
\index{evolution equation!on Wx[0,T]@on $W\times \left[ 0,T\right] $}$%
\begin{equation*}
\left\{ 
\begin{array}{rcl}
H^{-1}%
\dfrac{\partial H}{\partial t}
&=&-2\mathbf{i}\hat{F}_{H} \\ 
H\left( 0\right) 
&=&
H_{0}%
\end{array}%
\right. \quad \text{on}\quad W\times \left[ 0,T\right].
\end{equation*}%
Moreover, each $H_{t}$ approaches $H_{0}$ exponentially in all derivatives
over tubular segments $\Sigma_1(S)$ along the noncompact end.
\end{theorem}

\newpage
\section{Time-uniform convergence}
\label{Sec Time-uniform convergence}

There is a standard way \cite[\S 1.2]{ASDYM} to build a
functional on the space of Hermitian bundle metrics over a compact K\"{a}%
hler manifold the critical points of which, if any, are precisely the
Hermitian Yang-Mills metrics. This procedure is analogous
to the Chern-Simons construction, in that it amounts to integrating along paths a prescribed
first-order variation, expressed by a closed $1-$form. 
I will adapt this prescription to  $W$, restricting attention to metrics with suitable
asymptotic behaviour, and to the $K3$ divisors $D_{z}=\tau ^{-1}\left(
z\right) $ along the tubular end. 
\index{norm!functional}On one hand, the resulting functional $\mathcal{N}_{W}
$ will illustrate the fact that our evolution equation converges
to a HYM metric. On the other hand, crucially, the family $\mathcal{N}_{D_{z}}$ will mediate the role of stability in the time-uniform control of $%
\left\{ H_{t}\right\} $ over $W$.

\subsection{Variational formalism of the functional $\mathcal{N}$}
 
I will set up this analogous framework in some generality at first, defining an \emph{a priori} path-dependent functional $\mathcal{N}_{W}$ on a
suitable set of Hermitian metrics on $\mathcal{E}$. When restricted to
the specific $1-$parameter family $\left\{ H_{t}\right\} $ from our
evolution equation, we will see that $\mathcal{N}_{W}\left( H_{t}\right)$
is in fact decreasing and the study of its derivative will reveal that the $t\rightarrow \infty $
limit metric, if it exists, must be HYM on $\mathcal{E}$. Let 
\begin{displaymath}
        \mathcal{I}_0 := \left\{ h\in\func{End}\mathcal{E}  \left\vert\; h\text{ is Hermitian}, \; h \tublim 0\right. \right\}
\end{displaymath}
 denote the space of fibrewise Hermitian matrices
which decay exponentially along the tube.
\begin{lemma}   \label{Lemma 1st order var of curvature}
Let $H$ be a Hermitian bundle
metric and  $h$ an element of $\mathcal{I}_0$ and denote $%
\tau := H^{-1}h$; then the curvature of the
Chern connection of $H$ varies, to first order, by%
\begin{equation*}
        F_{H+h}=F_{H}
        +\bar{\partial}\partial _{H}\tau 
        +O( \left\vert\tau\right\vert ^{2}).
\end{equation*}

\begin{proof}
Set $g=H^{-1}\left( H+h\right) =1+\tau $, so that \cite[p.46]{4-manifolds}\cite[p.15]{ASDYM}%
\begin{equation*}
F_{H+h}=F_{H}+\bar{\partial}\left( g^{-1}\partial _{H}g\right) .
\end{equation*}%
Observing that $g^{-1}=1-\tau +O( \left\vert\tau\right\vert ^{2}) $, we expand the
variation of curvature:%
\begin{eqnarray*}
\bar{\partial}\left( g^{-1}\partial _{H}g\right) &=&-\left( g^{-1}.\bar{%
\partial}g.g^{-1}\right) \partial _{H}g+g^{-1}\bar{\partial}\partial _{H}g \\
&=&-\left( 1-\tau \right) \bar{\partial}\tau \left( 1-\tau \right) \partial
_{H}\tau +\left( 1-\tau \right) \bar{\partial}\partial _{H}\tau +O(
\left\vert\tau\right\vert ^{2}) \\
&=&\bar{\partial}\partial _{H}\tau +O( \left\vert\tau\right\vert ^{2}).
\end{eqnarray*}
\vspace{-1.3cm}\[\qedhere\]
\end{proof}
\end{lemma}

\begin{definition}      \label{Def space of Hermitian bundle metrics}
Let $\mathcal{H}_{0}$ be the set of smooth Hermitian metrics $H$ on $\mathcal{E}\rightarrow W$ such that:
\begin{displaymath}
\fbox{$
\begin{array}{c}
        H \tublim H_0. \\
\end{array}
$} 
\end{displaymath}

\begin{remark}       \label{Rem All about H0}
About the \emph{Definition}:
\begin{list}{\labelitemi}{\leftmargin=0em}
        \item[(1)] 
        The exponential decay $\left(\ref{Eq (e hat)0 decays exp}\right)$ implies $H_{0}\in \mathcal{H}_{0}$. Indeed, $\mathcal{H}_{0}$ is a \emph{star domain} in the affine space $H_0+\mathcal{I}_0$, in the sense that $H_{0}+\ell\left( H-H_{0} \right)\in \mathcal{H}_0$, $\forall\left(\ell,H\right)\in[0,1]\times \mathcal{H}_0$, with $H-H_{0}\in \mathcal{I}_0$. Thus $\mathcal{H}_{0}$ is contractible, hence \emph{connected} and \emph{simply connected}.

        \item[(2)] 
There is a well-defined notion of `infinitesimal variation' of a metric $H$, as an object in the `tangent space'
        \begin{equation*}
                T_{H}\mathcal{H}_{0}\simeq \mathcal{I}_{0}.
        \end{equation*}

        \item[(3)]
We know from $(\ref{Eq (e hat)0 decays exp})$ that $\hat{F}_{H_0} \overset{C^0}{\underset{S\rightarrow\infty}{\longrightarrow}} 0$, hence \emph{Lemma \ref{Lemma 1st order var of curvature}} implies
\begin{displaymath}
        \left\Vert \hat{F}_{H}\right\Vert _{L^{1}\left( W,\omega\right)}
        <\infty, \quad\forall H \in\mathcal{H}_0.
\end{displaymath}
        \item[(4)] 
Any `nearby' $H\in \mathcal{H}_{0}$, for which $\xi =\log H_{0}^{-1}H$ is well-defined (i.e.,  $\Vert H_{0}^{-1}H-I\Vert_{C^{0}(W,\omega)}<1$), is joined to $H_{0}$ by 
\begin{equation}        \label{eq def path in H0}
\begin{array}{c}
        \gamma :\left[ 0,1\right] \rightarrow \mathcal{H}_{0} \\ 
        \gamma \left( \ell \right) =H_{0}\e{\ell \xi }%
\end{array}.
\end{equation}%
Clearly $\ell\xi \tublim 0$, so $\gamma(\ell)\in \mathcal{H}_0$, $\forall\ell \in [0,1]$.
        \item[(5)] 
        Given any $T>0$, the solutions $\left\{ H_{t}\right\} _{t\in \left[0,T\right] }$ of our flow form a path in $\mathcal{H}_{0}$%
, since $\hat{F}_{H_{t}}$ decays exponentially along
the tube for each  $t$ [\emph{Theorem \ref{thm sols for all time and exponentially decaying}}].
\end{list}
\end{remark}
\end{definition}

Following \cite[pp.8-11]{ASDYM}, let $\theta \in \Omega ^{1}\left( \mathcal{H}_{0},\Omega ^{1,1}\left(
W\right) \right)$ be given by 
\begin{equation}        \label{eq def theta_H}      
\begin{array}{c}
        \theta _{H} : T_{H}\mathcal{H}_{0}\rightarrow \Omega ^{1,1}(W)\\
        \theta _{H}\left( k\right) 
        = 2\mathbf{i}\func{tr}\left(H^{-1}.k.F_{H}\right)
\end{array}.
\end{equation}
Then we may, at first formally, write
\begin{equation}        \label{eq def rhoW}
        \left( \rho _{W}\right) _{H}\left( k\right) 
        =\int_{W}\theta _{H}\left(k\right) \wedge \omega ^{2},
\end{equation}%
which will define a smooth $1-$form on any domain $H_0\in\mathcal{U}\subset 
\mathcal{H}_{0}$ where the integral converges, for all $H\in \mathcal{U}$ and all $k\in T_{H}\mathcal{H}_{0}$. The
crucial fact is that $\rho $ is identically zero precisely at the HYM metrics:
\begin{eqnarray*}
        \left. \left( \rho _{W}\right) _{H}=0\right.  &\Leftrightarrow &
        \int_{W}\func{tr}\left(H^{-1}k.F_{H}\right)\wedge \omega ^{2}=0,
        \quad \forall k\in T_{H}\mathcal{H}_{0} \\
        &\Leftrightarrow &\hat{F}_{H}=\left( F_{H},\omega \right) =0.
\end{eqnarray*}%
Following the analogy with Chern-Simons formalism, this suggests
integrating $\rho _{W}$ over a path to obtain a function having the HYM
metrics as critical points. Given $H\in \mathcal{H}_{0}$, let $%
\gamma \left( \ell \right) =H_{\ell }$ be a path in $\mathcal{H}_{0}$ connecting $H$ to the reference
metric $H_{0}$, and form the evaluation of $\theta $ along $\gamma $:%
\begin{equation}        \label{eq def Phizao of l}
        \Phi ^{\gamma }\left( \ell \right) := 
        \left[ \theta _{\gamma }\left(\dot{\gamma}\right) \right] 
        \left( \ell \right) =2\mathbf{i}\func{tr}%
        \left( H_{\ell }^{-1}.\dot{H}_{\ell }.F_{H_{\ell }}\right) 
        \in \Omega^{1,1}\left( W\right) .
\end{equation}%
For instance, with $\gamma $ as in $\left( \ref{eq def path in H0}\right)$,
we have $H_\ell^{-1}.\dot{H}_{\ell }=\underset{\e{-\ell \xi}}{\underbrace{H_\ell^{-1}.H_0}}. \left(\frac{\partial}{\partial\ell} \e{\ell\xi}\right)=\xi$ and
\begin{eqnarray*}
        \left( \rho _{W}\right) _{H_{\ell }}\left( \dot{H}_{\ell }\right)
        &=&\int_{W}\Phi ^{\gamma }\left( \ell \right) \wedge \omega ^{2}
        \quad= \quad2\mathbf{i}\int_{W}\func{tr}\xi.F_{H_\ell}\wedge\omega^{2}\\
        &=&2\mathbf{i}\int_{W}\func{tr}\xi.\hat{F}_{H_\ell}\dvol
\end{eqnarray*}%
is well-defined near $H_0$, since $\xi=\func{log}H_0^{-1}H$ is bounded and $\hat{F}_{H_{\ell }}$ is integrable [\emph{Remark \ref{Rem All about H0}}]. Thus, in this setting at least, the integral is rigorously defined:%
\begin{equation}        \label{Eq func N first def}
                        \index{norm!functional!$\mathcal{N}_{W}$}
        \mathcal{N}_{W}^{\gamma}\left( H\right):=\int_{\gamma }\rho _{W}.
\end{equation}%
There is a convenient relation between $\Phi ^{\gamma }\left( \ell \right) $ and the rate of change of the `topological' charge density $\func{tr}F^2$ along $\gamma$, which will be useful later:

\begin{lemma}       \label{lemma: trF^2 and Phizao}
Let $\left\{ \gamma \left( \ell \right) =H_{\ell }\right\} \subset 
\mathcal{H}_{0}$ be a $1-$parameter family of metrics on $\mathcal{E}$; then the evaluation 
$\Phi ^{\gamma }$ from $\left( \ref{eq def Phizao of l}\right) $ satisfies 
\begin{equation*}
        -\mathbf{i}\bar{\partial}\partial \Phi ^{\gamma }\left( \ell \right)
        =\frac{d}{d\ell }\func{tr}F_{H_{\ell }}^{2}.
\end{equation*}
\begin{proof}
From the first variation of $F$ [\emph{Lemma \ref{Lemma 1st order var of curvature}}] and the Bianchi identity:%
\begin{eqnarray*}
        \frac{d}{d\ell }\func{tr}F_{H_{\ell }}^{2} &=&
        2\func{tr}\left( \frac{d}{d\ell }F_{H_{\ell }}\right)\wedge F_{H_{\ell}}
        \quad=\quad2\func{tr}\bar{\partial}\partial _{H_{\ell }}
        \left( H_{\ell }^{-1}.%
        \dot{H}_{\ell }\right) \wedge F_{H_{\ell }} \\
        &=&-\mathbf{i}\bar{\partial}\partial \Phi ^{\gamma }\left( \ell \right).
\end{eqnarray*}
\vspace{-1.3cm}\[\qedhere\]
\end{proof}
\end{lemma}
By the same token, if we restrict attention to our family $\left\{\gamma(t)=
H_{t}\right\} \subset \mathcal{H}_{0}$ satisfying the evolution equation%
\begin{equation}        \label{eq evolution again}
\left\{ 
\begin{array}{c}
H^{-1}\dfrac{\partial H}{\partial t}=-2\mathbf{i}\hat{F}_{H} \\ 
H\left( 0\right) =H_{0}%
\end{array}%
\right. ,  
\end{equation}%
set $\mathcal{N}_{W}\left( H_{0}\right) =0$ and write for short $\Phi _{t}:=\Phi^{\gamma}(t)$, we obtain a real smooth
function%
\index{norm!functional!$\mathcal{N}_{W}$} 
\begin{equation}        \label{Eq func N using Phizao}
        \mathcal{N}_{W}\left( H_{T}\right) 
        =\int_{0}^{T}\left( \rho _{W}\right)_{H_{t}}
        \left(\dot{H}_{t}\right) dt
        =\int_{0}^{T}\left( 
        \int_{W}\Phi _{t}\wedge \omega^{2}\right) dt.  
\end{equation}

\begin{proposition}
\label{Prop Nw is defined and bounded above}
The function $\mathcal{N}_{W}\left( H_{t}\right) $ 
is well-defined, $\forall t\in \left[ 0,\infty\right[ $, and
\begin{equation*}
        \frac{d}{dt}\mathcal{N}_{W}\left( H_{t}\right) 
        =-\frac{2}{3}\Vert \hat{F}_{H_{t}}\Vert _
        {L^{2}\left( W\right) }^{2}.
\end{equation*}
\begin{proof}
Using the evolution equation $\left( \ref{eq evolution again}\right)$:
\begin{eqnarray*}
        \frac{d}{dt}\mathcal{N}_{W}\left( H_{t}\right) 
        &=& \left( \rho _{W}\right)_{H_{t}}\left( \dot{H}_{t}\right) 
        = \int_{W}\Phi _{t}\wedge \omega ^{2} \\
        &=& 2\int_{W}\func{tr}\underset{2\hat{F}_{H_{t}}}{\underbrace{\mathbf{i}%
        H_{t}^{-1}.\dot{H}_{t}}}.\underset{\frac{1}{6}\hat{F}_{H_{t}}\dvol}{\underbrace{F_{H_{t}}
        \wedge \omega ^{2}}} 
        = \frac{2}{3}\int_{W}\func{tr}\hat{F}_{H_{t}}^{2}\dvol \\
        &=& -\frac{2}{3}\Vert \hat{F}_{H_{t}}\Vert
        _{L^{2}\left( W\right)}^{2}
\end{eqnarray*}%
and this is finite, as $\hat{F}_{H_{t}}$ decays exponentially along $W$ [\emph{Proposition \ref{Prop e hat decays
exponentially along S}}].
\end{proof}
\end{proposition}

The above \emph{Proposition} confirms that we are on the right track: if the $\left\{ H_{t}\right\} $ converge to a smooth metric $H=H_{\infty }$ at all, then $H$ must be HYM.

Finally, our definition of $\mathcal{N}_{W}$ by integration of $\rho _{W}$ is \emph{a
priori }path-dependent and we have briefly examined two examples [(\ref{Eq func N first def}) and 
(\ref{Eq func N using Phizao})] which will be relevant in the ensuing analysis. Let us now eliminate the dependence, so that these settings
are, in fact,  equivalent.
\begin{lemma}   \label{Lemma antisymmetric theta}
Let $H\in\mathcal{H}_0$ and $h,k\in T_H\mathcal{H}_0 \cong\mathcal{I}_0 $; in the terms of (\ref{eq def theta_H}), the difference 
\begin{equation*}
\eta _{H}\left( h,k\right) := \frac{1}{2\mathbf{i}}\left( \theta
_{H+h}\left( k\right) -\theta _{H}\left( k\right) \right)
\end{equation*}
 is anti-symmetric to first order, modulo $\func{img}\partial +\func{img}\bar{\partial}$.
\begin{proof}
In the notation of \emph{Lemma \ref{Lemma 1st order var of curvature}} and setting $\sigma:= hK^{-1}$, the
anti-symmetrisation of $\eta _{H}$ is 
\begin{equation}     \label{antisymmetrisation form}
\begin{array}{r c l}
        \xi _{H}\left( h,k\right) &:= &\eta _{H}\left( h,k\right) 
        -\eta_{H}\left( k,h\right)   \\
        &=&\func{tr}\left( \left( \sigma .\tau -\tau .\sigma \right)
        .F_{H}+\sigma .\bar{\partial}\partial _{H}\tau 
        -\tau .\bar{\partial}\partial _{H}\sigma \right)+ \\
        &&+ \; O( \left\vert\sigma\right\vert.\left\vert\tau\right\vert ^{2})
        +O( \left\vert\tau\right\vert.\left\vert\sigma\right\vert ^{2}).\\
\end{array}
\end{equation}
The curvature of the Chern connection of $H$ obeys $F_{H}=\bar{\partial}%
\partial _{H}+\partial _{H}\bar{\partial}$, so%
\begin{eqnarray*}
\sigma .\bar{\partial}\partial _{H}\tau &=&\sigma .F_{H}.\tau -\sigma
.\partial _{H}\bar{\partial}\tau \\
&=&\sigma .F_{H}.\tau -\partial _{H}\left( \sigma .\bar{\partial}\tau
\right) -\bar{\partial}\left( \partial _{H}\sigma .\tau \right) +\tau .\bar{%
\partial}\partial _{H}\sigma ;
\end{eqnarray*}%
\emph{mutatis mutandis}, 
\begin{equation*}
\tau \bar{\partial}\partial _{H}\sigma =\tau .F_{H}.\sigma -\partial
_{H}\left( \tau .\bar{\partial}\sigma \right) -\bar{\partial}\left( \partial
_{H}\tau .\sigma \right) +\sigma .\bar{\partial}\partial _{H}\tau .
\end{equation*}%
Substituting these and using the cyclic property of trace in $\left( \ref%
{antisymmetrisation form}\right) $ we find%
\begin{eqnarray*}
        \xi _{H}\left( h,k\right) &=&\tfrac{1}{2}\func{tr}
        \left( \partial_{H}\left( \tau .\bar{\partial}\sigma 
        -\sigma .\bar{\partial}\tau \right) 
        +\bar{\partial}\left( \partial _{H}\tau .\sigma 
        -\partial _{H}\sigma .\tau\right) \right)\\
        &&+ \; O(\left\vert\sigma\right\vert.\left\vert\tau\right\vert^{2})
        +O( \left\vert\tau\right\vert.\left\vert\sigma\right\vert ^{2}) \\
        &\in &\func{img}\partial +\func{img}\bar{\partial} 
        \quad\text{ \emph{modulo} }\quad
        O(\left\vert\sigma\right\vert.\left\vert\tau\right\vert ^{2})
        +O( \left\vert\tau\right\vert.\left\vert\sigma\right\vert ^{2}).
\end{eqnarray*}
\vspace{-1.3cm}\[\qedhere\]
\end{proof}
\end{lemma}

\begin{corollary}       \label{Cor rho derivative of R is closed}
Let $\mathcal{U}\subset \mathcal{H}_{0}$ be a subset where the integral defining the $1-$form $\rho_{W}$ in (\ref{eq def rhoW}) converges, for all $H\in \mathcal{U}$ and all $k\in T_{H}\mathcal{H}_{0}$; then $\left.\rho_{W}\right\vert_\mathcal{U}$ is closed.

\begin{proof}
Recall that a $1-$form is closed precisely when its infinitesimal variation is symmetric to first order. In view of the previous \emph{Lemma}, it remains to check that $\xi
_{H}\left( h,k\right) \wedge \omega ^{2}$ integrates to zero modulo terms of higher order:%
\begin{equation*}
        \lim\limits_{S\rightarrow \infty }\int_{W_{S}}\xi _{H}\left( h,k\right)
        \wedge \omega ^{2}=0.
\end{equation*}%
Taking account of bi-degree and de Rham's theorem, modulo terms with decay  $O( \left\vert\sigma\right\vert.\left\vert\tau\right\vert ^{2})
        +O( \left\vert\tau\right\vert.\left\vert\sigma\right\vert ^{2})$, one has %
\begin{eqnarray*}
        \int_{W_{S}}\xi _{H}\left( h,k\right) \wedge \omega ^{2} &=&\tfrac{1}{2}%
        \int_{\partial W_{S}}\func{tr}\left[ \left( \tau .\bar{\partial}\sigma
        -\sigma .\bar{\partial}\tau \right) +\left( \partial _{H}\tau .\sigma
        -\partial _{H}\sigma .\tau \right)\right] \wedge \omega ^{2}\\ 
        &=&\tfrac{1}{2}\int_{\partial W_{S}}\func{tr}\left[ 2\tau .\bar{\partial}
        \sigma +2\partial _{H}\tau .\sigma -\nabla _{H}\left( \sigma .\tau\right)\right]
        \wedge \omega ^{2}\\
        &=&\int_{\partial W_{S}}\func{tr}\left( \tau .\bar{\partial}\sigma
        +\partial _{H}\tau .\sigma \right) \wedge \omega ^{2}
        \underset{S\rightarrow\infty}{\longrightarrow}0. 
\end{eqnarray*}%
On the other hand, from (the Calabi-Yau-Tian-Kovalev) \emph{Theorem \ref{Thm SU(3) mfds}} [cf. $\left( \ref{asymptotic forms}%
\right) $],
\begin{equation*}
        \left. \omega ^{2}\right\vert _{\partial W_{S}}
        =\kappa _{I}^{2}+O\left(\e{-S}\right).
\end{equation*}%
Consequently, as $S\rightarrow\infty$, the operation `$.\,\wedge\, \omega ^{2}$' annihilates all components of the $%
1-$form $\func{tr}\left( \tau .\bar{\partial}\sigma +\partial _{H}\tau
.\sigma \right)$ except those transversal to $D_{z}$ ($\left\vert
z\right\vert =\e{-S}$) in $\partial W_{S}\simeq D_{z}\times S^{1}$:%
\begin{eqnarray*}
\int_{\partial W_{S}}\func{tr}\left( \tau .\bar{\partial}\sigma +\partial
_{H}\tau .\sigma \right) \wedge \omega ^{2} &=&\int_{\partial W_{S}}\func{%
tr}\left( \tau .\frac{\partial \sigma }{\partial \bar{z}}d\bar{z}+\left( \left(\partial_{H}\tau \right)_{z} .\sigma dz\right)\wedge \left( \kappa _{I}^{2}+O\left( \e{-S}\right) \right) \right) \\
&= &\int_{\left\vert z\right\vert =e^{-S}}O\left( \left\vert
z\right\vert \right) \wedge \left( \kappa _{I}^{2}+O\left( \left\vert
z\right\vert \right) \right) \underset{\left\vert z\right\vert \rightarrow 0}%
{\longrightarrow }0.
\end{eqnarray*}%
Here we used that $\left\vert dz\right\vert ,\left\vert d\bar{z}\right\vert =
O\left( \e{-S}\right) =O\left( \left\vert z\right\vert \right) $ [\emph{%
Lemma \ref{lemma decay of coordinates dz dz bar}}], while $%
\tau $ and $\sigma $ also decay exponentially in all derivatives [\emph{Definition \ref{Def space of Hermitian bundle metrics}}].
\end{proof}
\end{corollary}

Since $\mathcal{H}_{0}$ is simply connected, we conclude that%
\index{norm!functional!$\mathcal{N}_{W}$} 
\begin{equation}        \label{eq func N path-indep}
        \boxed{\; \mathcal{N}_{W}
        =\int_{\gamma }\rho _{W} \;}  
\end{equation}%
doesn't depend on the choice of path $\gamma $.

\subsection{\texorpdfstring
        {A lower bound on `energy density' via $\mathcal{N}_{D_{z}}$}
        {A lower bound on `energy density' via NDz}}    
\label{Subsec conjectured uniform Lp-bound}

It is easy to adapt this prescription to the $K3$ divisors $D_{z}=\tau ^{-1}\left(
z\right) $ along the tubular end. 
\index{norm!functional}Such a family $\mathcal{N}_{D_{z}}$ will mediate the role of stability in the time-uniform control of $%
\left\{ H_{t}\right\} $ over $W$.
Still following \cite[pp.8-11]{ASDYM}, by analogy with (\ref{eq def theta_H}) and (\ref{eq def rhoW}), define $\theta_{z} \in \Omega ^{1}\left( \mathcal{H}_{0},\Omega ^{1,1}\left( D_{z}\right) \right) $  by 
\begin{equation}
\begin{array}{c}        \label{eq def (theta_z)_H}
\left(\theta_{z}\right) _{H}:T_{H}\mathcal{H}_{0}\rightarrow \Omega ^{1,1}\left( D_{z}\right) \\ 
\left(\theta_{z}\right) _{H}\left( k\right) =2\mathbf{i}\func{tr}
\left(\left.H\right\vert _{D_{z}}^{-1}.k.
\left.F_{H}\right\vert _{D_{z}}\right)%
\end{array}.  
\end{equation}%
and accordingly 
\begin{equation}        \label{eq def rhoD_z}
        \left( \rho _{D_z}\right) _{H}\left( k\right) 
        =\int_{D_{z}}\theta _{H}\left(k\right) \wedge \omega,
\end{equation}%
Setting each $\mathcal{N}_{D_{z}}\left( H_{0}\right) =0$, and choosing, at first, a curve   $\gamma
\left( \ell \right) =\left. H_{\ell }\right\vert _{D_{z}}$ of
Hermitian metrics on $\left. \mathcal{E}\right\vert _{D_{z}}$, we write [cf. $%
\left( \ref{eq def rhoW}\right)$]%
\index{norm!functional!$\mathcal{N}_{D_z}$} 
\begin{equation}
\mathcal{N}_{D_{z}}:= \int_{\gamma }\rho _{D_{z}}=\int_{\gamma
}\int_{D_{z}}\theta \wedge \omega.   \label{eq def NDz and NE}
\end{equation}%
Each $D_{z}$
being a compact complex surface, it is known  [ibid.] that this definition of $\mathcal{N}%
_{D_{z}}$ is in fact path-independent. 

In view of \emph{Proposition \ref{Prop Uniform upper bound}} below, which underlies this article's main result [\emph{Theorem \ref{Thm MAIN THEOREM}}], one would like to
derive, for  small enough $\left\vert z\right\vert $, a time-uniform \emph{lower} bound on the `energy density' given by the $\omega-$trace of the restriction of curvature  $\left.F_{H_t}\right\vert _{D_{z}}$:
\begin{displaymath}
\widehat{F_{t\vert z}} := \widehat{\left.F_{H_t}\right\vert_{D_z}} = \left( \left.F_{H_t}\right\vert_{D_z} , \left.\omega \right\vert_{D_z}\right).
\end{displaymath}
Recalling that $\xi _{t}\in \Gamma \left( 
\func{End} \mathcal{E} \right) $ is defined by $H_{t}=H_{0}e^{\xi _{t}}$ (hence is self-adjoint with respect to both metrics), write $\bar{\lambda}_{t}$ for its highest eigenvalue as in $\left( \ref{eq def lambda bar}\right)$ and set 
\begin{displaymath}
        L_{t}:=\sup_{W}\bar{\lambda}_{t}.
\end{displaymath}

\begin{claim}
\label{claim (F hat)_L^2(Dz) > c  over large set} 
There are constants $c,c'>0$  independent of $t$ and $z$ such that, for every $t \in \left]0,\infty \right[$, there exists an open set $A_t\subset\tau(W_\infty)\subset\mathbb{C}P^{1}$   of parameters [cf. (\ref{eq fibration tau})]  satisfying:
\begin{enumerate}[(1)]
        \item
                $\forall z\in A_t$, $\left\vert z\right\vert <\delta $ as                 in \emph{Definition \ref{def bundle E->W}}, i.e., $\left.\mathcal{E}\right\vert_{D_z}$                 is stable;
        \item 
                if $L_t\gg1$ then, in the cylindrical measure $\mu_{\infty}$ induced on $\tau(W_{\infty})$ by $ds^2+d\alpha^2$, with $z=e^{-s+\mathbf{i}\alpha}$  [cf. (\ref{asymptotic forms})], the following estimate holds:               \begin{displaymath}
                        \boxed{\int_{A_t}\Vert \widehat{F_{t\vert z}}\Vert^{2}_{L^2(D_z)}
                        ds\wedge d\alpha\ 
                        \geq {\frac{c}{2}\mu_{\infty}\left( A_t \right) \cdot
                        \left( 1+O\left(\frac{1}{L_{t}}\right) \right)}}\quad;
                \end{displaymath}
        \item
moreover, when $L_t\to \infty$, one has \;
$\displaystyle
\frac{\mu_{\infty}(A_t) }{ \sqrt{L_t}}\to c'.
$
\end{enumerate}
\end{claim}
Together with a uniform upper `energy bound' on $F_{H_t}$ over $W$,   \emph{Claim \ref{claim (F hat)_L^2(Dz) > c  over large set}} will suffice to establish the time-uniform $%
C^{0}-$bound on $\bar\lambda_t$. In the course of its proof, at last, the asymptotic stability assumption on $\mathcal{E}$ intervenes, by an instance of the following result:
\begin{lemma}   \label{Lemma NDz > xi_L^4/3}
        Suppose $\left.\mathcal{E}\right\vert_{D_z}$ is stable with Hermitian         Yang-Mills metric $\left.H_0\right\vert_{D_z}$; then there exists         a  constant $c_z>0$ such that 
        \begin{displaymath}
                \mathcal{N}_{D_{z}}\left( H_{t}\right)\geq c_z \left( \left\Vert\xi_t\right\Vert_{L^\frac{4}{3}(D_z)}-1                 \right).
        \end{displaymath}
        Moreover, one can choose $\delta>0$ in \emph{Definition \ref{def bundle E->W}} small enough to obtain a definite infimum $$
\left(\inf_{\left\vert z \right\vert<\delta}{c_{z}}\right)>0.
$$ 
\begin{proof}
The estimate for any given $D_z$ is the content of \cite[Lemma 24]{inf dets stab bdls...}. In particular, each constant $c_z$ is determined by the Riemannian geometry of $D_z$, as it comes essentially from an application of Sobolev's embedding theorem, hence it varies continuously with $z$. As $\left.\mc{E}\right\vert_{D}$ is itself an instance of our hypotheses, with metric $\left.H_0\right\vert_{D}$ over $D=D_0$, we have in this case $c_0>0$. Now one can certainly pick smaller $\delta>0,$ in \emph{Definition \ref{def bundle E->W}}, so that $\left\vert z\right\vert < \delta \Rightarrow c_z>\frac{c_0}{2}$, say.   
\end{proof}
\end{lemma}

Intuitively, our argument for \emph{Claim \ref{claim (F hat)_L^2(Dz) > c  over large set}} will proceed as follows: on a fixed $D_{z}$ far enough
along the tube, the quantity $\mathcal{N}_{D_{z}}\left( H_{t}\right) $ is controlled, in a certain sense, by the $\omega-$traced restriction of curvature   $\widehat{F_{t\vert z}}$ [\emph{Lemma \ref{Lemma Xsi bounds NDz}}, below]. On the other hand, the stability assumption implies that the same $\mathcal{N}_{D_{z}}\left(
H_{t}\right) $  controls above the slicewise norm  $\left\Vert \xi _{t}\right\Vert _{L^{\frac{4}{3}}(D_z)}$ [\emph{Lemma \ref{Lemma NDz > xi_L^4/3}}], so $\left.\xi _{t}\right\vert_{D_z}$ arbitrarily `big' would imply on $\widehat{F_{t\vert z}}$ being at least `somewhat big'. Moreover, if this happens at some $z_{0}$ then it must still hold over a `large' set $A_{t}\subset\mathbb{C}P^{1}\ $of parameters $z$, roughly proportional to the supremum $L_{t}=\left\Vert \bar\lambda_{t}\right\Vert _{C^{0}\left( D_{z_{0}}\right) }$ [\emph{Claim \ref{claim (F hat)_L^2(Dz) > c  over large set}}] at some $z_0$. 
Adapting the archetypical Chern-Weil technique [cf. \emph{Subsection \ref{Sub energy bounds}}], I  establish an absolute bound on the $L^2-$norm of $F_{H_{t}}$ over $W$ [estimate $\left( \ref{eq uniform energy bound} \right)$, below], so that the set $A_{t}$, carrying a `proportional amount of energy', cannot be too large in the measure $\mu_{\infty}$. Hence the supremum $L_t$, roughly of magnitude $\mu_{\infty}(A_t)$, can only grow up to a definite, time-uniform value. 

I start by proving essentially `half' of \emph{Claim \ref{claim (F hat)_L^2(Dz) > c  over large set}}:

\begin{lemma}       \label{Lemma Xsi bounds NDz}
There is a constant  $c_1>0$, independent of $t$ and $z$, such that 
\begin{equation*}
\mathcal{N}_{D_{z}}\left( H_{t}\right) \leq c_{1} \: L_{t}\:\Vert \widehat{F_{t\vert z}}\Vert_{L^2(D_z)} \qquad \forall t\in \left]0,\infty \right[.
\end{equation*}
\begin{proof}
Fixing $t>0$ and $\left\vert z\right\vert $ $<\delta $, let us simplify notation by $\xi =\xi _{t}$ and $\left\Vert .\right\Vert =\left\Vert
.\right\Vert _{L^{2}\left( D_{z}\right) }$, and consider the curve
\begin{equation*}
\begin{array}{c}
\gamma :\left[ 0,1\right] \rightarrow \left. \mathcal{H}_{0}\right\vert
_{D_{z}} \\ 
\gamma \left( \ell \right) =H_{0}\e{\ell \xi }%
\end{array}%
\end{equation*}%
with $\gamma \left( 1\right) =H_{t}$ and $\gamma _{\ell }^{-1}.\dot{\gamma}%
_{\ell }=\xi $. Using the first variation of curvature [\emph{Lemma \ref%
{Lemma 1st order var of curvature}}] we obtain $\frac{d}{d\ell }F_{\ell}=\bar{\partial}\partial _{\ell }\xi $, with $\partial_\ell:=\partial_{\gamma_{\ell}}$ and $F_\ell := F_{\gamma \left( \ell \right) }$.
Form [cf. (\ref{eq def Phizao of l})]%
\begin{equation*}
m\left( \ell \right) := \int_{0}^{\ell }\int_{D_{z}}\Phi ^{\gamma
}\left( \ell \right) \wedge \omega =2\mathbf{i}\int_{0}^{\ell }\left(
\int_{D_{z}}\xi .F_{\ell}\wedge \omega \right) d\ell
\end{equation*}%
so that $m\left( 1\right) =\mathcal{N}_{D_{z}}\left( H_{t}\right)$ and $%
m\left( 0\right) =0$; differentiating along $\gamma $ we have%
\begin{equation*}
m^{\prime }\left( \ell \right) =2\mathbf{i}\int_{D_{z}}\xi .F_{\ell}\wedge \omega .
\end{equation*}%
The function $m\left( \ell \right) $ is in fact convex:
\begin{displaymath}
        m^{\prime \prime }\left( \ell \right) 
        =2\mathbf{i}\int_{D_{z}}\func{tr}\xi .\left( \bar{\partial}_{\ell}
        \partial _{\ell }\xi \right) \wedge \omega 
        = 2\left\Vert \partial _{\ell }\xi \right\Vert ^{2}
        =\left\Vert \nabla_{\ell }\xi \right\Vert ^{2} \geq 0
\end{displaymath}
since $\xi $ is real, and so $\left\vert \partial _{\ell }\xi
\right\vert ^{2}=\left\vert \bar{\partial}_{\ell }\xi \right\vert ^{2}=\frac{%
1}{2}\left\vert \nabla _{\ell }\xi \right\vert ^{2}$. Now, by the mean value theorem, there exists some $\ell \in\left[ 0,1 \right]$ such that
\begin{eqnarray*}
        \mathcal{N}_{D_{z}}\left( H_{t}\right)
        &=&m\left( 1\right) \;=\; m(0) + m^{\prime }\left( \ell \right) 
        \;\leq\; m^{\prime}(1)
        \;\leq\; 2\int_{D_{z}}\left\vert \xi.F_{1}\wedge \omega \right\vert \\
        &\leq & c_{1}\;L_t \;\Vert \widehat{\left.F_{1}\right\vert_{D_z}}\Vert
\end{eqnarray*}%
using convexity and Cauchy-Schwarz, with $c_1:= 2\sup\limits_{\left\vert z \right\vert < \delta}\sqrt{\func{Vol}(D_z)}$ and the notation $F_1=F_{H_t}$.
\end{proof}
\end{lemma}
\begin{remark}          \label{rem N functional is convex}
                        \index{norm!functional!is convex}
Convexity implies that $\mathcal{N}_{D_{z}}\left( H_{t}\right) =m\left( 1\right) $ is \emph{positive} for all $t$, because $m^{\prime }\left( 0\right) =0$ gives an absolute minimum at $\ell =0$, so $m\left( 1\right) \geq m\left( 0\right) =0$.
\end{remark}
 
Finally, it can be shown that $\Delta\bar{\lambda}$ is (weakly) uniformly bounded \cite[p.246]{inf dets stab bdls...}, hence the maximum principle suggests that $\left\vert\xi_t\right\vert$ cannot `decrease faster' than a certain concave parabola along the cylindrical end. With that in mind, let us assume, for the sake of argument, that the following can be made rigorous:
\begin{quote}
In the terms of \emph{Claim \ref{claim (F hat)_L^2(Dz) > c  over large set}}, there exists a set $A_t$, `proportional' to $c'\sqrt{L_t}$, such that%
\begin{equation*}
        \left\Vert \xi_{t}\right\Vert _{L^{\frac{4}{3}}\left( D_{z}\right)}
        \geq c_{2} \; L_{t} \qquad \forall t\in\left]0,\infty\right[,
\end{equation*}
where $c_2$ is independent of $t$ and $z$.
\end{quote}
Then, together with \emph{Lemma \ref{Lemma NDz > xi_L^4/3}} and \emph{Lemma \ref{Lemma Xsi bounds NDz}}, this would prove \emph{Claim \ref{claim (F hat)_L^2(Dz) > c  over large set}}, with%
\begin{equation*}
c=\frac{2c_2}{c_1}\left(\inf_{\left\vert z \right\vert<\delta}{c_{z}}\right).
\end{equation*}
This is the heuristics underlying the proof, in the next \emph{Section}.

\subsection{Proof of \emph{Claim \ref{claim (F hat)_L^2(Dz) > c  over large set}}}
\label{Sec proof of conjecture}

Recall that one would like to
establish, for  small enough $\left\vert z\right\vert $, a time-uniform \emph{lower} bound on the `energy density' given by the $\omega-$trace of the restriction  $\left.F_{H_t}\right\vert _{D_{z}}$:
\begin{displaymath}
\widehat{F_{t\vert z}} := \widehat{\left.F_{H_t}\right\vert_{D_z}} = \left( \left.F_{H_t}\right\vert_{D_z} , \left.\omega \right\vert_{D_z}\right),
\end{displaymath}
 in the \emph{weak} sense that its $L^2-$norm over a cylindrical segment $\Sigma$ far enough down the tubular end is bounded below by a scalar multiple of $\vol \Sigma$. Moreover, the length of such a cylinder $\Sigma$ can be assumed roughly proportional to $L_{t}:=\sup_{W}\bl_{t}$, so that $L_t\gg0$ implies a large `energy' contribution. This  will yield the $C^0-$bound in \emph{Proposition \ref{Prop Uniform upper bound}}.

The strategy here consists, on one hand, of using the weak control over the Laplacian from \emph{Lemma \ref{Lemma Delta Lambda weakly bounded}} below to show that, around the furthest point down the tube where $L_{t}=\max\limits_W \bl_t$ is attained, the slicewise supremum of $\bl_{t}$ is always on top of a certain concave parabola $P_t$. 
On the other hand, the integral along the tube of the slicewise norms $\Vert \bl\Vert_{L^\frac{4}{3}(D_z)}$, which bounds below $\mathcal{N}_{D_z}(H_t)$ [\emph{Lemma \ref{Lemma NDz > xi_L^4/3}}], can be shown to be itself bounded below by those slicewise suprema, using again the weak bound on the Laplacian to apply a Harnack estimate on `balls' of a standard shape, which fill essentially `half' of the corresponding tubular volume. 

Since $\mathcal{N}_{D_z}(H_t)$ is controlled above by $\hat F_{H_t}$, in the sense of \emph{Lemma \ref{Lemma Xsi bounds NDz}}, this leads to the desired minimal `energy' contribution, `proportional' to the length (roughly $\sqrt{L_t}$) of the tubular segment underneath the parabola $P_t$.

\subsubsection{Preliminary analysis}
\label{subsubsec analytical lemmas}
Here I establish some basic wording and technical facts underlying the proof.
  
\begin{definition}      \label{def weak bound on Laplacian}
Given functions $f,g:W\rightarrow \R$, we denote
\begin{eqnarray*}
        \Delta f \wleq g 
        &\dot\Leftrightarrow&
        \int_W f\Delta\varphi \leq \int_W g\varphi ,
        \quad \forall \varphi\in C_c^\infty (W,\R^+).
\end{eqnarray*} 
In particular, a constant $\beta\in\R$ is called a \emph{weak bound on the Laplacian of $f$} if $\Delta f \wleq\beta$
\end{definition}
The following \emph{Lemma} is intuitively assumed in  \cite[p.246]{inf dets stab bdls...}, so the proof may also interest careful readers of that reference.
\begin{lemma}[weak maximum principle]   
\label{lemma weak max principle}
Let $f:W\to\R$ be a (nonconstant) nonnegative Lipschitz function, smooth away from a singular set $N_0$ with $\func{codim} N_0\geq3$; if $\Delta f\wleq 0$ [cf. \emph{Definition \ref{def weak bound on Laplacian}}] over a (bounded open) domain $U\subset W$, then $$\left. f \right\vert_{U}\leq \max\limits_{\partial U} f.$$

\begin{proof}
Outside of an  $\varepsilon-$neighbourhood $N_\varepsilon$ of the singular set $N_0$, the function $f$ is smooth, so the strong maximum principle on $U\smallsetminus N_\varepsilon$ implies, at the limit $\varepsilon\to0$, that any interior maximum point $q\in U$ must be in $N_0$.
Take then a (small) coordinate neighbourhood $q\in V\Subset U$ such that the gradient field $\nabla f$ `points inwards' to $q$ almost everywhere (a.e.) along $\partial V$. 

I claim one can choose  test functions  $\varphi\in C^\infty_c(W,\R^+)$,  with $\func{supp}(\varphi)\subset V$, such that  $\left(\nabla f,\nabla\varphi\right)$ is  practically nonnegative, in a suitable sense. Since $f$ is Lipschitz and a.e. differentiable, $f\in L^2_{1,loc}(W)\hookrightarrow L^2_1(V)$; by density, for every $\delta>0$ there exists $\varphi_\delta\in C^\infty_c(V)\hookrightarrow C^\infty_c(W)$ such that 
$$
\left\Vert f-\varphi_\delta\right\Vert^2_{L^2_1(V)}<\delta.
$$
Moreover, since $f$ is (nonconstant and) nonnegative and $q$ is a maximum, we may assume without loss that $\varphi_\delta\geq0$, so $\varphi_\delta$ is indeed a test function on $V$. Then 
\begin{eqnarray*}
        0\geq \int_W f\Delta\varphi_\delta
        =\int_V f\Delta\varphi_\delta
        =\lim_{\varepsilon\to0}\int_{V\smallsetminus N_\varepsilon}f\Delta\varphi_\delta
        \overset{(*)}{=} \int_V \left( \nabla f,\nabla\varphi_\delta \right)
        >-\delta
        \underset{\delta\to0}{\longrightarrow}0
\end{eqnarray*}
which either iterates all over $U$ to imply $f\equiv \max\limits_{\partial U} f$ or contradicts the assumption that $q$ is an interior point. Step $(*)$ is rigorous because $N_0$ has large enough codimension so that $\left.\dvol\right\vert_{\partial N_\varepsilon}=O(\varepsilon^2)$. Thus boundary terms in the integration by parts over $V\smallsetminus N_\varepsilon$ vanish when $\varepsilon\rightarrow 0$, again since $f$ is Lipschitz \cite[p.244]{inf dets stab bdls...}.
\end{proof}
\end{lemma}

Recalling that $\xi _{t}\in \Gamma \left( 
\func{End} \mathcal{E} \right) $ is defined by $H_{t}=H_{0}e^{\xi _{t}}$ (hence is self-adjoint with respect to both metrics), write $\bar{\lambda}_{t}$ for its highest eigenvalue; then:

\begin{lemma}           \label{Lemma Delta Lambda weakly bounded}
The Laplacian of $\bl$ admits a weak bound $\beta>0$:
\begin{displaymath}
        \Delta\bl_t \wleq \beta, \qquad\forall t>0. 
\end{displaymath}
\begin{proof}
We know from \cite[Lemma 25]{inf dets stab bdls...} that\begin{displaymath}
        \Delta\bl \wleq 
        2 \left( \left\Vert \hat{F}_{H_t} \right\Vert_{H_t} 
        + \left\Vert \hat{F}_{H_0} \right\Vert_{H_0}\right),
\end{displaymath}
and the right-hand side is controlled by the time-uniform bound on $\hat{F}_{H_t}$ [\emph{Corollary \ref{Cor Bounded F hat}}].
\end{proof}
\end{lemma}

\begin{lemma}   \label{lemma p-norm and 1-norm}
Let $\left(D,g\right)$ be a compact Riemannian manifold, $f\in L^\infty\left(D,\R^+\right)$,  $p>1$ and $x>0$; then there exists a constant $k_p=k_p\left(D,g\right)>0$ such that
\begin{displaymath}
        \left\Vert f\right\Vert_p \geq \frac{k_p}{F^x}\left\Vert f^{1+x}         \right\Vert_1
\end{displaymath}
with $\left\Vert.\right\Vert_q:=\left\Vert . \right\Vert_{L^q\left(D,g\right)},1<q\leq\infty,$ and $F:= \left\Vert f \right\Vert_\infty$.
\begin{proof}
It suffices to write
\begin{displaymath}
        \left\Vert f \right\Vert_p := \left(\int_D f^p \; d\func{Vol}_g\right)^{\frac{1}{p}}
        \geq \left(\int_D f^p\left( \frac{f}{F} \right)^{xp}
        d\func{Vol}_g\right)^{\frac{1}{p}} \\
        = F^{-x} \left\Vert f^{1+x} \right\Vert_p
\end{displaymath}
then apply H\"older's inequality, finding $k_p=\vol_g(D)^{\frac{1}{p}-1}$.
\end{proof}
\end{lemma}

\subsubsection{A concave parabola as lower bound}
In tubular coordinates $\left\vert z\right\vert=\e{-s}$, the supremum of $\bl_t$ on a transversal slice along $W_\infty$ defines a smooth function
\begin{displaymath}
\begin{array}{rcrcl}
        \ell_t&:&\R^+&\rightarrow&\R^+\\
        &&\ell_t(s)&:=& \sup\limits_{\partial W_s} \bl_t.\\
\end{array}
\end{displaymath}
Moreover, for each $t>0$, denote $S_t$ the `furthest length' down the tube at which $L_t$ is attained, i.e., $S_t=\max \left\{s\geq0\left\vert L_t=\ell_t\left( s \right)\right.\right\}$, and set
\begin{displaymath}
        I_t := \left[ S_t,S_t + \delta^{+}_t\right]\subset \R^+
\end{displaymath}
with $\delta^{+}_t:= \frac{1}{2}\left(\sqrt{1+\frac{8}{\beta}L_t}-1\right)$ and $\beta$ as in \emph{Lemma \ref{Lemma Delta Lambda weakly bounded}}.
\begin{lemma}           \label{lemma parabola} 
For each $t>0$, the transversal supremum $\ell_t$ is bounded below over $I_t$ by the concave parabola $P_t(s-S_t):= L_t -\frac{\beta}{2}(s-S_t)(s-S_t+1)$, i.e.,
\begin{displaymath}
        \ell_t(s)\geq P_t(s), \quad \forall s\in I_t.
\end{displaymath}
\begin{proof}
Fix $t>0$, $\delta^+_t\geq\delta>0$ and set $J_t\left( \delta \right):= \left]S_{t}-1,S_{t}+\delta\right[ $ (I suppress henceforth the $t$ subscript everywhere, for clarity). The parabola takes the value $P(-1)=P(0)=L$ at the points $S-1$ and $S$, and its concavity $-P''_t=\beta$ is precisely the weak bound on $\Delta\bl$ [\emph{Lemma \ref{Lemma Delta Lambda weakly bounded}}], so we have $\Delta\left( \bl - P\right)\wleq0$ on $\overline{J(\delta)}$ and  \emph{Lemma \ref{lemma weak max principle}} gives
\begin{eqnarray*}
         \left.\left( \bl - P \right)\right\vert_{J(\delta)}
        &<& \max \left\{ 
        \sup\limits_{s=S-1} \left( \bl - P \right),  
        \sup\limits_{s=S+\delta}\left( \bl - P \right)\right\}\\
        \Rightarrow \quad0=\ell(S)-P(S) &<& \ell(S+\delta)-P(S+\delta),
        \quad \forall 0<\delta<\delta^+ 
\end{eqnarray*}
since by assumption $\ell(S-1)\leq L=P(-1)$. Hence $\ell\geq P,\quad \forall s\in I$.
\end{proof}
\end{lemma}

\begin{remark}  \label{rem delta+,epsilon,t}
Fixing $0<\epsilon<1$ and setting $\delta^+_{\epsilon,t}:=\frac{1}{2}\left(\sqrt{1+\frac{8}{\beta}(1-\epsilon)L_t}-1\right)$, one has\begin{equation}        \label{eq ell bounded below}
        \ell_t(s) \geq \epsilon L_t, \quad \forall s\in I_{\epsilon,t}
        :=\left[ S_t,S_t+\delta^+_{\epsilon,t} \right].
\end{equation}
Clearly, if $L_t\to\infty$, the interval length grows quadratically as fast:
$$\frac{\delta^+_{\epsilon,t}}{\sqrt{L_t}}\to c_\epsilon':=\sqrt{\frac{2(1-\epsilon)}{\beta}}.
$$
\end{remark}


Moreover, we can use a Harnack-type estimate over transversal slices to establish the following inequality for the finite cylinder under the parabola:
\begin{lemma}   \label{lemma Moser}
Given $t>0$ and $0<\epsilon<1$, let $\Sigma_{t}(\epsilon):= I_{\epsilon,t}\times D \simeq \overline{W_{S_t+\delta^+_{\epsilon,t}}\setminus W_{S_t}}$ be the finite cylinder along $W$, under the parabola $P_t$ of \emph{Lemma \ref{lemma parabola}}, determined by the interval of length $\delta^+_{\epsilon,t}$ on which $(\ref{eq ell bounded below})$ holds, and suppose ${2\mathbf{\pi}}\delta^+_{\epsilon,t}\in\mathbb{N}$; then, for each $x>0$, there exist time-uniform constants $ a_{x,\epsilon},b_{x,\epsilon}>0$ such that, if $L_t>b_{x,\epsilon}$, the following estimate holds:
\begin{equation}        \label{eq lower bound over 'half' the cylinder}
        \int_{\Sigma_{t}(\epsilon)}\bl_t^{1+x} \; 
        d\vol_{\omega_\infty} 
        \geq a_{x,\epsilon} 
        \cdot \delta^+_{\epsilon,t}
        \cdot \left(L_t-b_{x, \epsilon} \right)^{1+x}. 
\end{equation}
\begin{proof}  
Again let me suppress the $t$ subscript, for tidiness, and work all along in the cylindrical metric $\omega_\infty$. For each $s\in I_\epsilon$, let $p_s\in \partial W_s \simeq D_s\times S^1$ be a point on the corresponding transversal slice such that the maximum $\ell(s)=\bl(p_s)$ is attained, and form the `unit' open cylinder $B_s\subset \Sigma_\epsilon$ of length $\frac{1}{2\mathbf{\pi}}$ (so that the volume integral over $B_s$ along the $S^1\times I_\epsilon$ directions is 1), centered on $p_s$, such that 
\begin{displaymath}
        \vol B_s = \vol \left( B_s\cap D_s \right) = \frac{1}{2} \vol D,
\end{displaymath}
where $\vol D\equiv \vol D_s$ denotes the (same) four-dimensional volume of (every) $D_s$.

By \emph{Lemma \ref{Lemma Delta Lambda weakly bounded}}, $\bl$ is a weak subsolution of the elliptic problem $\Delta u=\beta$, so we may apply Moser's iteration method   \cite[Theorem 8.25]{Gilb-Trud}\cite[Theorem 2]{Moser} over an open set $V$ such that $\partial W_s\subset V\subset B_s$ to obtain local boundedness of $\bl$ in terms of its $L^{1+x}(B_s,\omega_\infty)-$norm. Indeed we can choose $V$ such that:   
$$
\ell(s)
\leq 
C_x\left[
        \left( 
        \frac{1}{\vol B_s}
        \right)^{\frac{1}{1+x}}
        \left\Vert\bl\right\Vert_{1+x} +1
\right]
$$
for some constant $C_x=C_x(B_s,\beta,x)>0$ which in fact is uniform in $s$, as all $B_s$ are congruent by translation. Setting $\displaystyle b_{x,\epsilon}:=\frac{C_x}{\epsilon}$ and using (\ref{eq ell bounded below}) we have
\begin{equation}        \label{eq Moser}
        \left[
        \epsilon \left(L-b_{x,\epsilon}\right)
        \right]^{1+x}
\leq
\frac{1}{\vol B_s}\int_{B_s}\bl^{1+x}.
\end{equation}
In particular, one can choose at most $2\mathbf{\pi}\delta^+_\epsilon\in \mathbb{N}$ values $s_j\in I_\epsilon$ such that the corresponding $B_{s_j}$ are necessarily disjoint, and form their union 
\begin{displaymath}
        B\left( \epsilon \right):=\coprod\limits_{j=1}^{2\mathbf{\pi}\delta^+_\epsilon}
        B_{s_j}.
\end{displaymath}
Clearly $\vol B(\epsilon)\geq \frac{1}{2}\left(2\mathbf{\pi}\delta^+_\epsilon\right)\vol D$. Now, the statement about averages $(\ref{eq Moser})$ goes over to the disjoint union, hence
\begin{displaymath}
        \int_{\Sigma(\epsilon)}\bl^{1+x} 
        \geq \int_{B(\epsilon)}\bl^{1+x}
        \geq \vol B(\epsilon) 
        \cdot 
        \left[
        \epsilon \left(L-b_{x,\epsilon}\right)
        \right]^{1+x}
\end{displaymath}
 which proves the \emph{Lemma}, with 
 $a_{x, \epsilon} := \left(\mathbf{\pi} \vol D\right)\epsilon^{1+x}$.
\end{proof}
\end{lemma}

\subsubsection{End of proof}

It is now just a matter of putting together the previous results. Recalling the final remarks of \emph{Subsection \ref{Subsec conjectured uniform Lp-bound}}, we know from  \emph{Lemma \ref{Lemma NDz > xi_L^4/3}} and \emph{Lemma \ref{Lemma Xsi bounds NDz}} that
\begin{displaymath}
        L\Vert \hat F \Vert_2 \geq k'
        \left( \left\Vert \bl\right\Vert_{\frac{4}{3}} -1\right)
\end{displaymath}  
over each $D_z$ sufficiently far down the tube, for a uniform constant $k'>0$. Choosing $0<\epsilon<1$ and $x>0$, integrating over $A_\epsilon:= I_\epsilon\times S^1$ and applying \emph{Lemma \ref{lemma p-norm and 1-norm}} we have
\begin{displaymath}
        \int_{A_\epsilon} \Vert \hat F \Vert_2  \;ds\wedge d\alpha
        \geq  \frac{k''}{L^{1+x}}
        \left(\int_{A_\epsilon}\left\Vert \bl^{1+x} \right\Vert_1 
        ds\wedge d\alpha \right)
        -k'\frac{\delta_\epsilon^+}{L}
\end{displaymath}
where $k'':= k' \cdot k_\frac{4}{3}$ is still a uniform constant. Moreover, by \emph{Lemma \ref{lemma Moser}}, the integral term is bounded below by 
$k'' \cdot a_{x,\epsilon} \cdot\delta^+_{\epsilon } \cdot
        \left(1-\frac{b_{x, \epsilon}}{L} \right)^{1+x}$, so H\"older's inequality gives
\begin{eqnarray*}
        \left(\vol A_\epsilon \right)^{\frac{1}{2}}
        \left(\int_{A_\epsilon}\Vert \hat F \Vert_2^2 \; ds\wedge d\alpha\right)^\frac{1}{2}
        &\geq& \delta^+_{\epsilon}
        \left(
                k'' \cdot a_{x,\epsilon}  \cdot 
                \left(1-\frac{b_{x, \epsilon}}{L} \right)^{1+x} 
                -\frac{k'}{L}
        \right).
\end{eqnarray*}
Since the interval $I_\epsilon$ has length precisely $\delta^+_\epsilon$, we have
$\mu_{\infty}\left( A_\epsilon \right)=\vol A_\epsilon =2\mathbb{\pi}\delta^+_\epsilon$ and so
\begin{eqnarray*}
         \int_{A_\epsilon}\Vert \hat F \Vert_2^2 \; ds\wedge d\alpha         &\geq& \frac{\mu_{\infty}\left( A_\epsilon \right)}{4\mathbb{\pi}^2}
    \left(
                k'' \cdot a_{x,\epsilon}  \cdot 
                \left(1-\frac{b_{x, \epsilon}}{L} \right)^{1+x} 
                -\frac{k'}{L}
        \right)^2
\end{eqnarray*}
which yields the \emph{Claim}, choosing e.g. $\epsilon=\frac{1}{2}$, $x=1$ and [cf. \emph{Remark \ref{rem delta+,epsilon,t}}]
\begin{displaymath}
        c=
        \left( 
        \frac{k'' \cdot a_{1,\frac{1}{2}}}{\sqrt{2}\mathbf{\pi}}
        \right)^2 
        =
        \frac{\left(k'\right)^2\left(\vol D\right)^{\frac{3}{2}}}{32}
        \quad \text{and} \quad   
        c'=2\mathbb{\pi}\sqrt{\frac{1}{\beta}}.\\
\end{displaymath}

\newpage
\section{Conclusion}
 
\subsection{Solution of the Hermitian Yang-Mills problem}
\label{sect time-uniform convergence}

Let $\left\{ H_{t}\right\} $ be the family of smooth Hermitian metrics on $%
\mathcal{E}\rightarrow W$ given for arbitrary finite time by \emph{%
Proposition \ref{Prop smooth solution for arbitrary finite time}}. In order
to obtain a HYM metric as $H=\lim\limits_{t\rightarrow \infty }H_{t}$ it would suffice to show that $\left\{ H_{t}\right\} $ is $C^{0}-$bounded, for then it
is actually $C^{\infty }-$bounded on any compact subset of $W$ and the limit $H$ is
smooth [\emph{Lemmas \ref{Lemma smooth sols on Ws for finite time}}, \emph{%
\ref{Lemma Hi is in Lp2} }and \emph{\ref{Lemma smoothness of limit}}].
Concretely, this would mean improving the constant $C_{T}$ in $\left( \ref%
{eq bound on sigma}\right) $ to a time-uniform bound $C_{\infty }$ or, what is the same, controlling the sequence $\bar{\lambda}_{t}$ of highest
eigenvalues of $\xi _{t}=\log H_{0}^{-1}H_{t}$ [cf. $\left( \ref{eq def
lambda bar}\right) $]:%
\begin{equation}
\boxed{\begin{array}{cr}
\left\Vert \bar{\lambda}_{t}\right\Vert _{C^{0}\left( W\right) }\leq
C_{\infty },& \forall t>0. \\
\end{array}  }\label{eq time-unif bound on lambda bar}
\end{equation}%
I will show that this task reduces essentially to \emph{Claim \ref{claim (F hat)_L^2(Dz) > c  over large set}}, as the problem $\left( \ref{eq time-unif bound on
lambda bar}\right) $ amounts in fact to controlling the size of the set $A_t$ where the `energy density' $\widehat{F_{t\vert z}}$ is bigger than a definite constant.
I begin by stating the announced upper bound:
\begin{proposition}     \label{prop negative energy} 
A solution $\left\{H_t\right\}$ to the evolution equation
given by \emph{Theorem \ref{thm sols for all time and exponentially decaying}} satisfies
the \emph{negative energy} condition\begin{eqnarray}        \label{eq uniform energy bound}
        E(t):= \int_{W}
        \left( 
        \left\vert F_{H_{t}}\right\vert^{2}
        -
        \left\vert F_{H_{0}}\right\vert ^{2}
        \right) \dvol%
        &\leq& 
        0,\qquad \forall t\in \left]0,\infty \right[.
\end{eqnarray}%
\begin{proof}
The curvature of a Chern connection splits orthogonally as $F=\hat{F}.\omega  \oplus  F^{\perp }$ in $\Omega ^{1,1}\left( \func{End}\mathcal{E}\right)$, so
$\left\vert F\right\vert ^{2}=\vert F^{\perp }\vert
^{2}+\vert \hat{F}\vert ^{2}$
(setting $\left\vert \omega \right\vert =1$). On the other hand, the Hodge-Riemann  equation $\left( \ref{Eq Hodge-Riemann}\right) $ [see \emph{Appendix A}] reads
 \begin{equation*}
        \func{tr}F^{2}\wedge \omega 
        =\left( \vert F^{\perp }\vert^{2}
        -\vert \hat{F}\vert ^{2}\right) \omega ^{3}.
\end{equation*}%
Comparing we find $\left\vert F\right\vert ^{2}\omega ^{3}=\func{tr}F^{2}\wedge \omega +2\vert \hat{F}\vert ^{2}\omega ^{3}$ [cf. \emph{Subsection \ref{Sub energy bounds}}], so
\begin{eqnarray*}
        E(t)&=&\int_{W}\left(\func{tr}F_{H_{t}}^{2}
        -\func{tr}F_{H_{0}}^{2}\right)\wedge \omega
        +2 \int_{W}\left(\hat{e}_t - \hat{e}_0\right)\omega ^{3}\\
        \overset{\frac{d}{dt}}{\Longrightarrow}\quad
        \dot{E}(t)&\leq&\int_W \left(-\mathbf{i}\bar{\partial}\partial           \Phi_t\right)\wedge\omega +2\int_W\left(-\Delta\hat{e}_t\right)\omega^{3}\\
        &\leq&\lim_{S\rightarrow\infty}\int_{\partial W_S}
        2\left[\partial \func{tr} \left(\hat{F}_{H_t}.F_{H_t}\right)\wedge\omega
        +\frac{\partial\hat{e}_t}{\partial\nu} \left.\dvol\right\vert_{\partial         W_S}\right]\\
        &=&0
\end{eqnarray*}
using \emph{Lemma \ref{lemma: trF^2 and Phizao}}
along with its \emph{Proof} and $\left(\frac{d}{dt}+\Delta\right)\hat{e}_t \leq 0$ as in $(\ref{eq (d/dt+D)e^ < 0})$, then complex integration by parts [\emph{Lemma \ref{Lemma integration by parts on W}} in \emph{Appendix \ref{App exterior algebra}}] and the Gauss-Ostrogradsky theorem, and finally the exponential decay $\hat{F}_{H_t} \tublim 0$, a direct consequence from \emph{Proposition \ref{Prop e hat decays exponentially along S}} and \emph{Corollary \ref{Cor H -> H0 exponentially in all derivatives}}. This shows that $E(t)$ is non-increasing, while obviously $E(0)=0$.
\end{proof}
\end{proposition}
From now on I will write, in cylindrical coordinates,  $D_s$ for $D_z$ when $\left\vert z\right\vert = e^{-s}$. Reasoning as above we find, for a Chern connection on $\left.\mathcal{E}\right\vert_{D_s}$,%
\begin{equation}       \label{Eq energy bound on surface Ds}
\left\vert F\right\vert ^{2}\omega ^{2}=\func{tr}F^{2}
+2\vert \hat{F}\vert ^{2}\omega ^{2}.
\end{equation}
Denoting, for convenience, the component of $\dvol$ containing the `tubular' factor $ds\wedge d\alpha$ by   
$$
\dvtub:=\; \tfrac{1}{2} ds\wedge d\alpha
        \wedge \left(\kappa_I +d\psi \right)^2,
$$
the remainder certainly satisfies  $\widetilde{d\psi}:=\dvol-\dvtub=O\left(e^{-s}\right)$. Moreover, around each $D_s$ we have
$\displaystyle
        {\left.\dvol\right\vert_{D_s}} 
        = \tfrac{1}{2}\left(\left.\omega\right\vert_{D_s}\right)^2 
        = \left.\tfrac{1}{2}\left(\kappa_I +d\psi \right)^2 \right\vert_{D_s}
$,
so that 
\begin{eqnarray}        \label{eq tubular decomp volume forms}
        \dvtub&=&ds\wedge d\alpha\wedge{\left.\dvol\right\vert_{D_s}}.
\end{eqnarray}

\begin{lemma}   \label{lemma Energy on cylindrical domain}
Let $\Sigma=\tau^{-1}(A)\subset W_\infty$ be any cylindrical domain along the tubular end, parametrised by $A\subset\mathbb{C}P^{1} $; then
\begin{eqnarray*}
        \int_{\Sigma} \left(\left\vert F_t\right\vert^2
        - \left\vert F_0 \right\vert^2 +R_0\right) 
        \dvtub         
        &\geq&
        2\int_{A} \Vert \widehat{F_{t \vert s}} \Vert^2_{L^2(D_s)}
        ds\wedge d\alpha 
\end{eqnarray*}
where $R_0$ decays exponentially along the tube. 
\begin{proof}
Isolating the component of curvature along each transversal slice $D_s$ in the asymptotia (\ref{Eq asymptotics of F_H0}) of $F_{0}$ , we have $\left\vert F_0 \right\vert^2 =\vert F_{0\vert s} \vert^2+R_{0}$ over $W_\infty$, where the remainder indeed satisfies $R_0 = O\left(e^{-s}\right)$. The Hodge-Riemann property  (\ref{Eq energy bound on surface Ds}) and the decomposition (\ref{eq tubular decomp volume forms}) then give, by Fubini's theorem,
\begin{eqnarray*}
        \int_{\Sigma} \left(\left\vert F_0 \right\vert^2 -R_0 \right)
        \dvtub
        &=& \int_{A} \left\{\int_{D_{s}}\left\vert F_{0\vert s} \right\vert^2
        \left.\dvol\right\vert_{D_s}\right\} ds\wedge d\alpha \\
        &=& \int_{A} \left\langle c_2\left( \left.\mc{E}\right\vert_{D_s}\right),\left[ D_s \right]\right\rangle
        ds\wedge d\alpha.        
\end{eqnarray*}
On the other hand, again by (\ref{Eq energy bound on surface Ds}), 
\begin{eqnarray*}
        \int_{\Sigma } \left\vert F_t \right\vert^2 
        \dvtub
        & \geq& \int_{\Sigma } \left\vert F_{t\vert s} \right\vert^2
        \dvtub\\
        &= &\int_{A} \left\{\int_{D_s} \left\vert F_{t\vert s}\right\vert^2
        {\left.\dvol\right\vert_{D_s}}\right\} ds\wedge d\alpha \\
        &\geq &\int_{A}\left\{ \left\langle c_2 \left( \left.\mathcal{E}\right\vert_{D_s}         \right),\left[ D_s \right]\right\rangle 
        +2 \int_{D_s}\vert \widehat{F_{t\vert s}} \vert^2 \left.\dvol\right\vert_{D_s}         \right\} ds\wedge d\alpha,
\end{eqnarray*}
using the normalisation $c_1(\left.\mathcal{E}\right\vert_{D_s})=0$ [cf. discussion following (\ref{eq fixed det over D})] so that only $c_2 \left( \left.\mathcal{E}\right\vert_{D_s}\right)$ contributes to the integral of $\func{tr}F_{t\vert s}^{2}$. Comparing and cancelling the topological terms yields the result.\end{proof}
\end{lemma}
In the terms of \emph{Claim \ref{claim (F hat)_L^2(Dz) > c  over large set}}, denote by $\Sigma_t := A_t \times D_s$ the finite cylinder $\tau^{-1}\left( A_t \right)$ along the tubular end $W_\infty$. Then the \emph{Lemma}  gives a lower  estimate on the curvature over $\Sigma_t$:
\begin{eqnarray}       \label{eq energy over finite cylinders}
        \quad\int_{\Sigma_t} \left(\left\vert F_t\right\vert^2
        - \left\vert F_0 \right\vert^2 +R_0\right) 
        \dvtub         &\geq&2 \int_{A_t} \Vert \widehat{F_{t \vert s}} \Vert^2_{L^2(D_s)}
        ds\wedge d\alpha.
\end{eqnarray} 

This discussion culminates at the following result:
\begin{proposition}     \label{Prop Uniform upper bound}
Under the hypotheses of  \emph{Claim \ref{claim (F hat)_L^2(Dz) > c  over large set}}, there exists a constant $C_\infty$, independent of $t$ and $z$, such that
\begin{displaymath}
L_t \leq C_\infty , \quad \forall t\in \left]0,\infty \right[.
\end{displaymath}

\begin{proof}

Either $L_t$ is uniformly bounded in $t$ and there is nothing to prove, or  $\displaystyle\frac{\mu_{\infty}\left( A_t \right)}{\sqrt{L_{t}}}\underset{t\to\infty}{\longrightarrow} c'$ and so (a subsequence of) the family of sets $\left\{A_t\right\}_{0<t<\infty}$ gets arbitrarily $\mu_\infty-$large as $t\rightarrow \infty$.
In the latter case, the term $O\left( \frac{1}{L_t} \right)$ in the \emph{Claim} becomes negligible, so
\begin{eqnarray*}
        c.\mu_{\infty}\left( A_t \right) 
        &\leq& \int_{\Sigma_t} \left(\left\vert F_t\right\vert^2
        - \left\vert F_0 \right\vert^2 +R_0\right) 
        \dvtub.
\end{eqnarray*}

 Let us examine the energy contributions of each component of $W=W_0\cup W_\infty$ [cf. (\ref{eq cylindrical picture})]. On one hand, the complement $W_\infty\smallsetminus\Sigma_t$ is a cylindrical domain, which by \emph{Lemma \ref{lemma Energy on cylindrical domain}}  has nonnegative energy in the model metric $\dvtub$; so we may extend the above integral over $W_\infty$. On the other hand, in the actual metric $\dvol$ we have:
\begin{eqnarray*}
        E(t)&=&\left(\int_{W_0}+\int_{W_\infty} \right)
        \left(\left\vert F_t\right\vert^2
        - \left\vert F_0 \right\vert^2 \right) 
        \dvol\\
        &\geq&
       -\left\Vert F_0\right\Vert^2 _{L^2(W_0,\omega)} 
        + \int_{W_\infty}
        \left(\left\vert F_t\right\vert^2
        - \left\vert F_0 \right\vert^2\right) 
        \dvol.
\end{eqnarray*}
To connect both facts, consider the ratio of volume forms
over $W_\infty$
$$
\dvtub=f.\dvol
$$
and extend it to a  bounded positive function $f:W\to\R^+$. Then 
\begin{eqnarray*}
        c.\mu_{\infty}\left( A_t \right) 
        &\leq& \left(\sup_W f\right)\left(E(t)
        +\left\Vert F_0\right\Vert^2 _{L^2(W_0,\omega)} 
        +\left\Vert R_{0}\right\Vert_{L^1(W,\omega)}\right)
\end{eqnarray*}
 and we know from the negative energy condition (\ref{eq uniform energy bound}) that this is bounded above, uniformly in $t$. Hence the cylinder $\Sigma_t$ cannot stretch indefinitely and, by contradiction, there must exist a uniform upper bound$$ 
C_\infty \simeq 
\left[
\frac
{\sup f} 
{c.c''}
\left(\left\Vert F_0\right\Vert^2 _{L^2(W_0)} 
+\left\Vert R_{0}\right\Vert_{L^1(W)}
\right)
\right]^2.\qedhere
$$
\end{proof}
\end{proposition}

Finally, replacing the uniform bound for $C_T$     in $\left(\ref{eq bound on sigma} \right)$,  this control cascades into the exponential $C^0-$decay in $\left( \ref{eq exp decay C0} \right)$ and hence the $C^\infty-$decay in \emph{Corollary \ref{Cor H -> H0 exponentially in all derivatives}}. By \emph{Proposition \ref{Prop Nw is defined and bounded above}}, the limit metric $H=\lim\limits_{t\rightarrow\infty}H_t$ satisfies $\hat F_H=0$. We have thus proved the following instance of the HYM problem:

\begin{theorem}         \label{Thm MAIN THEOREM}
Let $\mc{E}\rightarrow W$ be stable at
infinity [cf. \emph{Definition \ref{def bundle E->W}}], equipped with a reference metric $H_{0}$ [cf. Definition \ref{def reference metric H0}], over an asymptotically
cylindrical $SU\left( 3\right) -$manifold $W$ as given by the
Calabi-Yau-Tian-Kovalev \emph{Theorem \ref{Thm SU(3) mfds}}, and let $\left\{
H_{t}=H_{0}e^{\xi _{t}}\right\} $ be the $1-$parameter family of Hermitian
metrics on $\mathcal{E}$ given by \emph{Theorem \ref{thm sols for all time
and exponentially decaying}}; the limit $H=\lim\limits_{t\rightarrow\infty}H_t$ exists and is a smooth Hermitian Yang-Mills metric on  $\mathcal{E}$, exponentially asymptotic  in all derivatives [cf. \emph{Notation \ref{Not finite cylinder}}] to $H_0$ along the tubular end of $W$: %
\index{asymptotic decay!of HYM metric}%
\begin{equation*}       \index{metric!Hermitian Yang-Mills}
\fbox{$\begin{array}{c}
        \hat{F}_{H}=0, \quad H \tublim H_0. \\
\end{array}$}
\end{equation*}
\end{theorem}

\subsection{Examples of asymptotically stable bundles}
\label{Subsec Examples of ASB}

It is fair to ask whether there are any holomorphic bundles at all satisfying
the asymptotic stability conditions of \emph{Definition \ref{def bundle E->W}%
}, thus providing concrete\emph{\ }instances for  \emph{Theorem \ref{Thm MAIN THEOREM}}. 

We know, on one hand, that Kovalev's base manifolds [cf. \emph{Definition \ref{def base manifold (W,w)}}] are K\"ahler 3-folds
coming from blow-ups $\bar W=\func{Bl}_C X$, where $X$ is Fano and the curve $C=D\cdot D$ represents the self-intersection of a $K3$ divisor $D\in\left\vert -K_X \right\vert$ \cite[6.43]{kovalevzao}. In addition, several of the examples provided (e.g. $X=\mathbb{C}P^3$, complete intersections etc.) are nonsingular, projective and satisfy the \emph{cyclic conditions}: $\pic(X)=\mathbb{Z}$ and $\pic(D)=\mathbb{Z}$. 
On the other hand, certain linear monads over projective varieties of the above kind yield stable bundles as their middle cohomology \cite{Jardim Bull Braz}. Those are called \emph{instanton monads} and have the form
\begin{equation}        \label{eq instanton monad}
\begin{array}{ccccccccc}
0 & \rightarrow & \oo\left(-1\right)^{\oplus c} & \overset{\alpha}{\longrightarrow} & \oo^{\oplus2+2c} & \overset{\beta}{\longrightarrow} & \oo\left(1\right)^{\oplus c} & \rightarrow & 0.\\
\end{array}
\end{equation} Then,  denoting $K:=\ker\beta$,  the middle cohomology $\mc{E}:= K/\img \alpha$ is always a stable bundle with
        \begin{displaymath}
                \rank(\mc{E})=2,\quad c_1(\mc{E})=0, 
                \quad c_2(\mc{E})=c\cdot h^2
        \end{displaymath} 
where  $h:= c_1(\oo(1))$ is the hyperplane class and $c\geq1$ is an integer.

In those terms, twist the monad by $\oo(-d)$, $d=\deg D$, so the relevant data fit in the following \emph{canonical diagram}:
\begin{displaymath}
\begin{array}{ccccccccc}
  &  & 0 &  &  &  &  &  &  \\
  &  & \downarrow &  &  &  &  &  &  \\
  &  & \oo\left(-(d+1)\right)^{\oplus c} &  &  &  &  &  &  \\
  &  & \downarrow &  &  &  &  &  &  \\
0 & \rightarrow & K(-d) & \longrightarrow & \oo(-d)^{\oplus 2+2c} & \longrightarrow & \oo(-\left(d-1\right))^{\oplus c} & \rightarrow & 0 \\
 &  & \downarrow &  &  &  &  &  &  \\
0 & \rightarrow & \mc{E}(-d) & \longrightarrow & \mc{E} & \longrightarrow & \left.\mc{E}\right\vert_D & \rightarrow & 0 \\
 &  & \downarrow &  &  &  &  &  &  \\
 &  & 0 &  &  &  &  &  &  \\
\end{array}
\end{displaymath}
Computing cohomologies, one checks by Hoppe's criterion \cite[pp.165-166]{Okonek-Spindler-Schneider} that indeed $\left.\mc{E}\right\vert_D$ is stable, i.e., the bundle $\mc{E}$ is asymptotically stable. 

\noindent NB.: Fixing  $\bar{W}=\func{Bl}_{C}\left(\mathbb{C}P^{3}\right)$ and $c=1, d=4$ as monad parameters, this specialises to the well-known null-correlation bundle \cite{Okonek-Spindler-Schneider}\cite{Barth}.

A detailed study with many more examples of asymptotically stable bundles over base manifolds admissible by Kovalev's construction  will be published separately in \cite{FAE}.

\subsection{Future developments}
\label{subsec preview}
Kovalev's construction culminates at producing new examples of \emph{ compact} 7-manifolds with holonomy strictly equal to $G_2$. Loosely speaking, the asymptotically cylindrical  $G_{2}-$structures on certain `matching' pairs $W'\times S^1$ and $W''\times S^1$ are superposed along a truncated gluing region down the tubular ends, using cut-off functions to yield a global 
$G_{2}-$structure $\varphi $ on the compact manifold 
\begin{displaymath}
        M_S:= W'\tsum W'',
\end{displaymath}
defined by a certain `twisted' gluing procedure. This structure can be chosen to be torsion-free,
by a `stretch the neck' argument on the length parameter $S$. Step (2) in our broad programme, as proposed in the \emph{Introduction}, is the corresponding problem of gluing the $G_2-$instantons obtained in the present article. 

First of all, one needs to extend the twisted sum operation to bundles $\mc{E}\ii\to W\ii$ in some natural way, respecting the technical `matching' conditions.
More seriously, however, one should expect transversality to 
play an important role. As mentioned in \emph{Remark \ref{rem deformation complex}}, the easiest way to proceed  is probably to restrict
attention at first to \emph{acyclic }instantons \cite[p.25]{floer},
i.e., whose gauge class `at infinity' is \emph{isolated} in its moduli space $%
\mathcal{M}_{\left. \mathcal{E}\right\vert _{D}}$. 
This would require, of course,
the existence of asymptotically stable bundles which are also \emph{asymptotically rigid}, in the obvious sense. Fortunately, Kovalev's construction admits (at least) a certain class of prime  Fano 3-folds  $X_{22}\hookrightarrow\pp^{13}$, thoroughly studied by Mukai   \cite[\S3]{MukaiFano3-folds}\cite[\S3]{MukaiBundlesOfK3}, which admit bundles exhibiting precisely  those properties, thus such investigation does not seem  void from the outset.

Furthermore, if an asymptotically rigid gluing result does not present any further meanders, one might then consider the full question of transversality towards a gluing theory for families of instantons.
A detailed study of this matter will appear in the sequel, provisionally cited as \cite{G2II}.

\newpage
\appendix

\section{Facts from geometric analysis}
\label{App exterior algebra}

I collect in this Appendix three analytical results from both complex and Riemannian geometry which, although well-known to specialists, should be stated in precise terms and in compatible notation due to their importance in the text. 

\subsection{Integration by parts on complex manifolds with boundary}
\begin{lemma}[Integration by parts]\label{Lemma integration by parts on W}

Let $X^{n}\subseteq W$ be a compact complex (sub)manifold (possibly $n=3$), $\Phi $ a $\left( 1,1\right)-$form, $\Omega $ a closed $\left( n-2,n-2\right) -$form and $f$ a
meromorphic function on $X$; then 
\begin{equation*}
        \int_{X}
        \Phi \wedge dd^{c}f\wedge \Omega 
        =\int_{X}
        f.\left( -\mathbf{i}\bar{\partial}\partial \Phi \right) \wedge \Omega         +\mathbf{i}\int_{\partial X}
        \left( \Phi \wedge \bar{\partial}f+f.\partial \Phi \right)\wedge\Omega.
\end{equation*}
\begin{proof}
By the Leibniz rule and Stokes' theorem, using $d=\partial +\bar{\partial}$
and taking account of bi-degree, we have%
\begin{eqnarray*}
\int_{X}\Phi \wedge dd^{c}f\wedge \Omega &=&\int_{X}\Phi \wedge \mathbf{i}%
\partial \bar{\partial}f\wedge \Omega \\
&=&\mathbf{i}\int_{\partial X}\Phi \wedge \bar{\partial}f\wedge \Omega 
\underset{\left( \ast \right) }{\underbrace{-\int_{X}\mathbf{i}\partial \Phi
\wedge \bar{\partial}f\wedge \Omega }}
\end{eqnarray*}%
and again%
\begin{eqnarray*}
\left( \ast \right) &=&\mathbf{i}\int_{\partial X}f.\partial \Phi \wedge
\Omega +\int_{X}f.\left( -\mathbf{i}\bar{\partial}\partial \Phi \right)
\wedge \Omega .
\end{eqnarray*}
\vspace{-1.3cm}\[\qedhere\]
\end{proof}
\end{lemma}

\subsection{The Hodge-Riemann bilinear relation}
 
The curvature on a K\"{a}hler $n-$fold splits as $F=\hat{F}.\omega \oplus
F^{\perp }\in \Omega ^{1,1}\left( \func{End}\mathcal{E}\right) $, so: 
\begin{equation*}
F^{2}\wedge \omega ^{n-2}=\hat{F}^{2}.\omega ^{n}+( F^{\perp })
^{2}\wedge \omega ^{n-2}.
\end{equation*}%
The Hodge-Riemann pairing $\left( \alpha ,\beta \right) \longmapsto \alpha
\wedge \beta \wedge \omega^{n-2} $ on $\Omega ^{1,1}\left( W\right) $ is
positive-definite along $\omega $ and negative-definite on the primitive
forms in $\left\langle \omega \right\rangle ^{\perp }$ \cite[pp.39-40]%
{Huybrechts} (with respect to the reference Hermitian bundle metric); since
the curvature $F$ is real as a bundle-valued $2-$form, we have:%
\begin{equation}        \label{Eq Hodge-Riemann}
        \func{tr}F^{2}\wedge \omega ^{n-2}
        =\left( \vert F^{\perp}\vert ^{2}
        -\vert \hat{F}\vert ^{2}\right) \omega ^{n}
\end{equation}%
using that $\func{tr}\xi ^{2}=-\left\vert \xi \right\vert ^{2}$ on the
Lie algebra part.

\newpage
\subsection{Gaussian upper bounds for the heat kernel} 
\label{App Gaussian upper bound}
The following instance of \cite[Theorem 1.1]{Heat kernel bounds} stems from a long series, going back to J. Nash (1958) and D. Aronson (1971) [op.cit.], of generalised `Gaussian' upper bounds (i.e., given by a Gaussian exponential on the geodesic distance $r$) for the heat kernel $K_t$ of a Riemannian manifold: 
\begin{theorem} \label{Thm Gaussian bound on heat kernel}
        Let $M$ be an arbitrary connected Riemannian $n-$manifold, $x,y\in M$         and $0\leq T\leq\infty$; if there exist suitable [see below] real         functions $f$ and $g$ satisfying the `diagonal' conditions 
        \begin{eqnarray*}
                K_t(x,x)\leq\frac{1}{f(t)}\quad \text{and} \quad K_t(y,y)\leq\frac{1}{g(t)},\quad         \forall t\in \left(0,T\right),
        \end{eqnarray*}
        then, for any $C>4$, there exists $\delta=\delta(C)>0$ such that
        \begin{eqnarray*}
                \boxed{ K_t(x,y) \leq \frac{(cst.)}{\sqrt{f(\delta t)g(\delta t)}}\exp\left\{-\frac{r(x,y)^2}{Ct}\right\},                 \quad \forall t\in \left(0,T\right)}
        \end{eqnarray*}
        where $(cst.)$ depends on the Riemannian metric only. 
\end{theorem}
For the present purposes one may assume simply $f(t)=g(t)=t^{\tfrac{n}{2}}$, but in fact $f$ and $g$ can be \emph{much} more general [Op. cit. p.37].

\newpage

\section{Structures of proof for 
\emph{Theorem \ref{thm sols for all time and exponentially decaying}} and \emph{Theorem \ref{Thm MAIN THEOREM}}}
\label{App flowcharts}

\subsection{Proof of Theorem \ref{thm sols for all time and exponentially decaying}}

\tikzstyle{equation} = [circle, draw, text centered, node distance=12em]
\tikzstyle{cited} = [rectangle, draw, text width=11em, text centered, node distance=12em]
\tikzstyle{seta} = [draw, -stealth, thick]
\tikzstyle{linha} = [draw, thick]
\tikzstyle{elli}=[draw, ellipse, minimum height=8mm, text width=10em, text centered, node distance=12em]
\tikzstyle{teo} = [draw, diamond, text width=9em, text centered, minimum height=15mm, node distance=12em]

\begin{tikzpicture}
\begin{scriptsize}

\node [teo] 
(thmsolsforalltime) 
{\emph{Theorem \ref{thm sols for all time and exponentially decaying}}
\\Solutions for all time
\\$H(t) \overset{C^\infty}{\underset{S\rightarrow\infty}{\longrightarrow}}
H_{0},\; \forall t$
};

\node [elli, left of=thmsolsforalltime, xshift=-5em] 
(Propsmoothsolution)
{\emph{Proposition \ref{Prop smooth solution for arbitrary finite time}}
\\ $\exists! H(t),\;\forall t$
\\ $\sup\hat{F}_{H(t)}\leq \sup\hat{F}_{H_0}$
};
\path [seta] (Propsmoothsolution) -- (thmsolsforalltime);

\node [elli, right of=thmsolsforalltime, xshift=11em] 
(CorH->H0exp) 
{\emph{Corollary \ref{Cor H -> H0 exponentially in all derivatives}}
\\ Exponential decay
\\ $H(t) \overset{C^k}{\underset{S\rightarrow\infty}{\longrightarrow}}H_{0},\;\forall k$
};
\path [seta] (CorH->H0exp) -- (thmsolsforalltime);

\node [elli, above of=CorH->H0exp,yshift=-3em,xshift=0em] 
(Propregularityonsmalleropenset) 
{\emph{Proposition \ref{Prop regularity on smaller open set}}
\\ Regularity on $V'\Subset V$
\\ $\left\Vert . \right\Vert _{C^{1}\left( V^{\prime }\right) }
        \leq A\left\Vert . \right\Vert _{C^{0}\left( V\right) }$
};
\path [seta] (Propregularityonsmalleropenset) 
-- coordinate[midway](MID-Propregularityonsmalleropenset) (CorH->H0exp);

\node [equation, above left of=thmsolsforalltime,xshift=3.5em, yshift=-0em] 
(eqboundonLaplacian)
{(\ref{eq bound on Laplacian})};
\path [linha] (eqboundonLaplacian) --  +(4.92,0)  |- (MID-Propregularityonsmalleropenset);

\node [equation, above of=eqboundonLaplacian, yshift=-8em] 
(eqExpDecayC0)
{(\ref{eq exp decay C0})};
\path [linha] (eqExpDecayC0) --  +(4.92,0)  |- (MID-Propregularityonsmalleropenset);

\node [elli, above left of=Propregularityonsmalleropenset, xshift=-2em, yshift=-1.60em]
(LemmaderivativesofFunifbounded)
{\emph{Lemma \ref{Lemma all derivatives of F are unif bounded}}
\\$F_{H_{S}}$ uniformly bounded in $L_{k}^{\infty }\left( W_{S}\right),\;\forall k$
};
\path [seta] (LemmaderivativesofFunifbounded) --  +(0,-1.5)  |- (MID-Propregularityonsmalleropenset);

\node [elli, above of=LemmaderivativesofFunifbounded, yshift=-2em] 
(CorHsTisinLp2)
{\emph{Corollary \ref{Cor Hs(T) is in Lp2}}
\\$H_{S}\left( t\right) $ uniformly 
\\bounded in $L_{2}^{p}\left( W_{S}\right)$
\\ and $H_{S}\left( T\right)$ is $C^{1}$
};
\path [seta] (CorHsTisinLp2) -- (LemmaderivativesofFunifbounded);

\node [cited, right of=LemmaderivativesofFunifbounded, xshift=5em] 
(ThmGaussianbound) 
{\emph{Theorem \ref{Thm Gaussian bound on heat kernel}}
\\Bound on heat kernel
\\$K_{t}\left( x,y\right) \leq \frac{(cst) }{t^{3}}\e{ -\frac{r\left( x,y\right) ^{2}}{Ct}} $
};
\path [seta] (ThmGaussianbound) -- coordinate[midway](MID-ThmGaussianbound) (LemmaderivativesofFunifbounded);

\node [cited, above of=ThmGaussianbound, yshift=-6em] 
(ThmKovalev)
{\emph{Theorem \ref{Thm SU(3) mfds}}
\\$\omega _{\infty } =\kappa _{I}+ds\wedge d\alpha$
};
\path [seta] (ThmKovalev) -- +(-2,0) -| (MID-ThmGaussianbound);

\node [cited, right of=CorHsTisinLp2, xshift=5em] 
(LemmaHiisinLp2)
{\emph{Lemma \ref{Lemma Hi is in Lp2}}};
\path [seta] (LemmaHiisinLp2) -- coordinate[midway](MID-LemmaHiisinLp2) (CorHsTisinLp2);

\node [elli, above of=Propsmoothsolution,yshift=-1.5em] 
(Corlimitmetricissmooth)
{\emph{Corollary \ref{Cor limit metric is smooth}}
\\$H_S(T)$ smooth, \;$\forall S,T>0$
};
\path [seta] (Corlimitmetricissmooth) -- (Propsmoothsolution);

\node [elli, above of=Corlimitmetricissmooth,yshift=-1.5em] 
(Lemmasmoothnessoflimit)
{\emph{Lemma \ref{Lemma smoothness of limit}}
\\$H_{i}\underset{i\rightarrow I}{\overset{C^{\infty}(X)}{\longrightarrow }}H$
};
\path [seta] (Lemmasmoothnessoflimit) -- coordinate[midway](MID-Lemmasmoothnessoflimit) (Corlimitmetricissmooth);
\path [seta] (LemmaderivativesofFunifbounded) -- coordinate[midway](MID-LemmaderivativesofFunifbounded)
 +(-5,0) |- (MID-Lemmasmoothnessoflimit);
 \path [linha] (CorHsTisinLp2) -- (MID-LemmaderivativesofFunifbounded);
\path [seta]  (eqboundonLaplacian) -- (Lemmasmoothnessoflimit);

\node [elli, above of=Lemmasmoothnessoflimit,yshift=7em] 
(CorUniqueness)
{\emph{Corollary \ref{cor uniqueness}}\\
Uniqueness
};
\path [seta]  (CorUniqueness) -- +(-2.5,0)  |-  (Propsmoothsolution);

\node [elli, above of=eqExpDecayC0,yshift=1em] 
(PropeHatDecaysExponentiallyWithS)
{\emph{Proposition \ref{Prop e hat decays exponentially along S}}
\\Decay of $\hat F$\\
$\hat{e}_{t} \leq \left(B\e{t}\right)\epsilon$
};
\path [seta]  (PropeHatDecaysExponentiallyWithS) -- (eqExpDecayC0);

\node [elli, above right of=CorHsTisinLp2, xshift=0.9em, yshift=8em] 
(CorHt->HTinC0)
{\emph{Corollary \ref{Cor H(t) -> H(T) in C0}}
\\ $H_{S}\left( t\right) 
\underset{t\to T}{\overset{C^{0}}{\longrightarrow }}H_{S}\left( T\right)$ and $H_{S}\left( T\right) $ is $C^0$
};
\path [seta] (CorHt->HTinC0) -- coordinate[midway](MID-MorHt->HTinC0) (MID-LemmaHiisinLp2) ;

\node [elli, above of=CorHsTisinLp2,yshift=-5.5em, xshift=0em]
(CorBoundedFhat)
{\emph{Corollary \ref{Cor Bounded F hat}}
\\$\sup_{W_{S}}\vert\hat{F}_{H_{S}(t)}\vert^{2}\leq B$ 
};
\path [linha] (CorBoundedFhat) --   (MID-MorHt->HTinC0);

\node [cited,above of=PropeHatDecaysExponentiallyWithS,yshift=-5.1em,xshift=-0em]
(Lemmamaxprinciplemfdsw)
{\emph{Lemma \ref{Lemma max principle mfds w bdry}
}\\Maximum principle on compact manifolds with boundary
};
\path [seta] (Lemmamaxprinciplemfdsw) -- coordinate[midway](MID-Lemmamaxprinciplemfdsw)
(CorBoundedFhat);
\path [seta]  (Lemmamaxprinciplemfdsw) -- (PropeHatDecaysExponentiallyWithS);
\path [seta]  (Lemmamaxprinciplemfdsw) -- +(0,2.1) |- (CorUniqueness);

\node [cited, above right of=Lemmamaxprinciplemfdsw,xshift=7em,yshift=-4em]
(Weitzenbock)
{Weitzenb\"ock formula};
\path [seta] (Weitzenbock) -- + (-1.8,0) -| (MID-Lemmamaxprinciplemfdsw);

\node [equation, above of=Weitzenbock, yshift=-9em]
(Dirichlet)
{(\ref{Dirichlet})};
\path [linha] (Dirichlet) -- + (-1,0) -| (MID-Lemmamaxprinciplemfdsw);

\node [cited, above of=Lemmamaxprinciplemfdsw, yshift=5em]
(Lemmamaxprinforsigma)
{\emph{Lemma \ref{Lemma max prin for sigma}}
\\$\left( \frac{d}{dt}+\Delta \right) \sigma \leq 0$ 
};
\path [seta] (Lemmamaxprinforsigma) -- +(0,-1)|- coordinate[midway](MID-Lemmamaxprinforsigma)
(CorHt->HTinC0);
\path [linha] (Lemmamaxprinciplemfdsw) -- +(0,1) |- (CorHt->HTinC0);

\node [elli,right of=Lemmamaxprinforsigma, yshift=-1em,xshift=3.5em]
(Prop-shorttimeexistence)
{\emph{Proposition \ref{Prop - short time existence}}
\\$H_S(t)$ exists, for $t\in[0,\varepsilon[$
};
\path [seta] (Prop-shorttimeexistence) -- (CorHt->HTinC0);

\node [elli, above right of=CorHt->HTinC0, yshift=-1.5em,xshift=-3em]
(remarks)
{\emph{\emph{\emph{Remarks \ref{Rem Co norm is complete}} and \emph{\ref{Rem conv sigma => conv lambda}}}}};
\path [seta] (remarks) -- (CorHt->HTinC0);

\end{scriptsize}
\end{tikzpicture}
NB.: Theorem-like results proved in the text are represented by ellipses,  individual equations by circles, cited facts  by rectangles and the main result by a diamond.

\subsection{Proof of \emph{Theorem \ref{Thm MAIN THEOREM}}}

\vspace{-0.7em}
\begin{tikzpicture}
\begin{scriptsize}

\node [teo]
(thmMAINTHEOREM)
{\emph{Theorem \ref{Thm MAIN THEOREM}}\\
$\hat{F}_{H}=0$, \\ $H \tublim H_0$};

\node [teo, right of=thmMAINTHEOREM, xshift=5em] 
(thmsolsforalltime2) 
{\emph{Theorem \ref{thm sols for all time and exponentially decaying}}
\\Solutions for all time
\\$H(t) \overset{C^\infty}{\underset{S\rightarrow\infty}{\longrightarrow}}
H_{0},\; \forall t$
};
\path [seta] (thmsolsforalltime2) -- (thmMAINTHEOREM);

\node [elli, left of=thmMAINTHEOREM, xshift=-13.5em] 
(PropNwIsDecreasing) 
{\emph{Proposition \ref{Prop Nw is defined and bounded above}}
\\$\frac{d}{dt}\mathcal{N}_{W} 
        =-\frac{2}{3}\Vert \hat{F}_{H_{t}}\Vert _
        {L^{2} }^{2}$
};
\path [seta] (PropNwIsDecreasing) -- (thmMAINTHEOREM);

\node [elli, above of=thmMAINTHEOREM, yshift=3em]
(propUniformEnergyBound)
{\emph{Proposition \ref{Prop Uniform upper bound}}\\
$L_t \leq C_\infty ,\;\forall t$
};
\path [seta] (propUniformEnergyBound) -- coordinate[midway](MID-propUniformEnergyBound) (thmMAINTHEOREM);

\node [elli, above right  of=propUniformEnergyBound, xshift=2em,yshift=0em]
(claim)
{\emph{Claim \ref{claim (F hat)_L^2(Dz) > c  over large set}}\\
$\int_{A_t}\Vert \widehat{F_{t\vert z}}\Vert^{2}
\hspace{-0.25em}\geq\hspace{-0.25em} 
\frac{c}{2}\mu_{_{\infty}}\hspace{-0.4em}\left( A_t \right)$\\
$
\frac{\mu_{\infty}(A_t) }{ \sqrt{L_t}}\to c'
$};
\path [seta] (claim) -- (propUniformEnergyBound);

\node [elli, right of=propUniformEnergyBound, xshift=5em,yshift=0em]
(lemmaEnergyoncylindricaldomain)
{\emph{Lemma \ref{lemma Energy on cylindrical domain}}\\
$\int_{\Sigma} \left(\left\vert F_t\right\vert^2
        - \left\vert F_0 \right\vert^2 +R_0\right) 
        $\\$ \geq 
        2\int_{A}\Vert \widehat{F_{t\vert s}}\Vert^{2}
$};
\path [seta] (lemmaEnergyoncylindricaldomain) -- (propUniformEnergyBound);

\node [elli, above left of=propUniformEnergyBound, yshift=0em] 
(eqUnifEnergyBound)
{\emph{Proposition \ref{prop negative energy}}\\
$E(t)\leq 0$};
\path [seta] (eqUnifEnergyBound) -- (propUniformEnergyBound);

\node [equation, left of=propUniformEnergyBound, xshift=-2em, yshift=-8em]
(eqBoundOnSigma)
{(\ref{eq bound on sigma})
};
\path [seta] (eqBoundOnSigma) --  +(1.66,0)  |- (MID-propUniformEnergyBound);

\node [equation, below of=eqBoundOnSigma, yshift=7.2em] 
(eqExpDecayC0)
{(\ref{eq exp decay C0}) 
};
\path [seta] (eqExpDecayC0) --  +(1.66,0)  |- (MID-propUniformEnergyBound);

\node [elli, above of=eqExpDecayC0, yshift=0em, xshift=-3em] 
(CorH->H0exp2) 
{\emph{Corollary \ref{Cor H -> H0 exponentially in all derivatives}}
\\ Exponential decay
\\ $H(t) \overset{C^k}{\underset{S\rightarrow\infty}{\longrightarrow}}H_{0},\;\forall k$
};
\path [seta] (CorH->H0exp2) --  +(2.5,0)  |- (MID-propUniformEnergyBound);
\path [seta] (CorH->H0exp2) -- (eqUnifEnergyBound);

\node [equation, right of=lemmaEnergyoncylindricaldomain,xshift=-2em] 
(EqEnergyBoundOnSurfaceDs)
{ (\ref{Eq energy bound on surface Ds})
};
\path [seta] (EqEnergyBoundOnSurfaceDs) -- (lemmaEnergyoncylindricaldomain);

\node [equation, above of=EqEnergyBoundOnSurfaceDs, yshift=-5em] 
(eqExpDecayC02)
{(\ref{Eq asymptotics of F_H0})
};
\path [seta] (eqExpDecayC02) -- (lemmaEnergyoncylindricaldomain);

\node [equation, below of=EqEnergyBoundOnSurfaceDs,yshift=5em] 
(eqTubularDecompVolumeForms)
{(\ref{eq tubular decomp volume forms})
};
\path [seta] (eqTubularDecompVolumeForms) -- (lemmaEnergyoncylindricaldomain);

\node [equation, right of=eqUnifEnergyBound,yshift=-3em,xshift=-3.5em] 
(eqTubularDecompVolumeForms)
{(\ref{Eq Hodge-Riemann})
};
\path [seta] (eqTubularDecompVolumeForms) -- (eqUnifEnergyBound);

\node [elli, above left of=eqUnifEnergyBound,xshift=-8em,yshift=9em] 
(CortrF2andPhizao)
{\emph{Lemma \ref{lemma: trF^2 and Phizao}}\\
        $
        -\mathbf{i}\bar{\partial}\partial \Phi\left( \ell \right)
        =\frac{d}{d\ell }\func{tr}F_{H_{\ell }}^{2}$
};
\path [seta] (CortrF2andPhizao) -- +(3,0)-|(eqUnifEnergyBound);

\node [elli, above left of=eqUnifEnergyBound, xshift=-3em,yshift=1em] 
(LemmaIntegrationByParts)
{\emph{Lemma \ref{Lemma integration by parts on W}}\\
Integration by parts\\
        $
                \int_{X}
        \Phi \wedge dd^{c}f\wedge \Omega 
        =\int_{X}
        f\left( -\mathbf{i}\bar{\partial}\partial \Phi \right) \wedge \Omega         $\\$+ (\dots)$
};
\path [seta] (LemmaIntegrationByParts) -- (eqUnifEnergyBound);

\node [elli, left of=eqUnifEnergyBound,xshift=-5em, yshift=0em] 
(PropeHatDecaysExponentiallyWithS2)
{\emph{Proposition \ref{Prop e hat decays exponentially along S}}
\\$H_{i}\underset{i\rightarrow I}{\overset{C^{\infty}(X)}{\longrightarrow }}H$
};
\path [seta] (PropeHatDecaysExponentiallyWithS2) -- (eqUnifEnergyBound);
\path [seta] (PropeHatDecaysExponentiallyWithS2) -- (PropNwIsDecreasing);
\node [elli, above right of=CortrF2andPhizao,xshift=0em, yshift=0em] 
(Lemma1stOrderVariationOfCurvature)
{\emph{Lemma \ref{Lemma 1st order var of curvature}}\\
Variation of curvature\\
$F_{H+h}=F_{H}
        +\bar{\partial}\partial _{H}\tau 
        $\\$\qquad
        +O(\tau^{2})$
};
\path [seta] (Lemma1stOrderVariationOfCurvature) -- (CortrF2andPhizao);

\node [elli, above right of=claim, yshift=-2em,xshift=2em] 
(LemmaPnorm1norm)
{\emph{Lemma \ref{lemma p-norm and 1-norm}}\\
$\left\Vert f\right\Vert_p \geq \frac{k_p}{F^x}\left\Vert f^{1+x}         \right\Vert_1$
};
\path [seta] (LemmaPnorm1norm) -- (claim);

\node [elli, above  of=LemmaPnorm1norm,yshift=-6em,xshift=0em] 
(RemDelta)
{\emph{Remark \ref{rem delta+,epsilon,t}}\\
$\frac{\delta^+_{\epsilon,t}}{\sqrt{L_t}}\to c_\epsilon'=\sqrt{\frac{2(1-\epsilon)}{\beta}}
$
};
\path [seta] (RemDelta) -- +(-2.5,0) -| (claim);

\node [elli, above of=claim,yshift=7em,xshift=0em] 
(LemmaMoser)
{\emph{Lemma \ref{lemma Moser}}\\
$\int_{\Sigma_{t}(\epsilon)}\bl_t^{1+x}  
        \geq 
        a_{x,\epsilon} 
        \cdot \delta^+_{\epsilon,t}
        \cdot \left(L_t-b_{x, \epsilon} \right)^{1+x} 
$
};
\path [seta] (LemmaMoser) -- (claim);
\path [seta] (RemDelta) -- (LemmaMoser);

\node [elli, above left of=claim,yshift=-2em,xshift=-1.0em] 
(LemmaNDzBiggerthanXi4/3)
{\emph{Lemma \ref{Lemma NDz > xi_L^4/3}}\\
$\mathcal{N}_{D_{\hspace{-0.1em}z}}
\hspace{-0.3em}\geq \hspace{-0.3em}
c_z (\Vert\xi_t\Vert_{L^\frac{4}{3}\hspace{-0.1em}(\hspace{-0.1em}D_{\hspace{-0.1em}z}\hspace{-0.1em})}
\hspace{-0.8em}-1)
$};
\path [seta] (LemmaNDzBiggerthanXi4/3) -- +(2.4,0) -| (claim);

\node [elli, above of=LemmaNDzBiggerthanXi4/3,yshift=-6em,xshift=0em] 
(LemmaXsiBoundsNDz)
{\emph{Lemma \ref{Lemma Xsi bounds NDz}}\\
$\mathcal{N}_{D_{\hspace{-0.1em}z}} 
\hspace{-0.3em}\leq\hspace{-0.3em}
c_{1} \hspace{-0.2em} L_{t}\hspace{-0.1em}\Vert \widehat{F_{t\vert z}}\Vert_{\hspace{-0.1em}L^2(\hspace{-0.1em}D_{\hspace{-0.1em}z}\hspace{-0.1em})}$
};
\path [seta] (LemmaXsiBoundsNDz) -- +(2.4,0) -| (claim);
\path [seta] (Lemma1stOrderVariationOfCurvature) -- (LemmaXsiBoundsNDz);

\node [elli, above of=RemDelta,yshift=1em,xshift=0em] 
(LemmaParabola)
{\emph{Lemma \ref{lemma parabola}}\\
$\ell_t(s)\geq P_t(s)$
};
\path [seta] (LemmaParabola) -- (LemmaMoser);

\node [elli, above of=LemmaXsiBoundsNDz,yshift=4.3em,xshift=0em] 
(LemmaWeakLap)
{\emph{Lemma \ref{Lemma Delta Lambda weakly bounded}}\\
Weak Laplacian bound\\
$
\Delta\bl\wleq\beta
$
};
\path [seta] (LemmaWeakLap) -- (LemmaMoser);

\node [elli, above of=LemmaParabola,yshift=-4em,xshift=0em] 
(LemmaWeakMaxPrin)
{\emph{Lemma \ref{lemma weak max principle}}\\
\hspace{-0.6em}Weak \\maximum principle\\
$
\left. f \right\vert_{U}\leq \max\limits_{\partial U} f
$
};
\path [seta] (LemmaWeakMaxPrin) -- coordinate[midway](MID-LemmaWeakPrin)
        (LemmaParabola);
\path [seta] (LemmaWeakMaxPrin) -- (LemmaParabola);
\path [seta] (LemmaWeakLap) -- (MID-LemmaWeakPrin);

\end{scriptsize}
\end{tikzpicture}

\newpage




\begin{thebibliography}{99}


\bibitem[B-S]{Bryant-Salamon} Robert L. Bryant $\&$ Simon Salamon, \emph{On
the construction of some complete metrics with exceptional holonomy}, Duke
Math. J. \textbf{58} (1989) 829.

\bibitem[Bae]{Baez} John C. Baez, \emph{The Octonions}, Bull. Amer. Math. Soc. 39 (2002), 145-205; \emph{errata}: B. Am. Math. Soc. 42 (2005), 213.

\bibitem[Bar]{Barth} Wolf P. Barth, \emph{Some properties of stable rank-}$2$%
\emph{\ vector bundles on }$\mathbb{P}_{n}$, Math. Ann. \textbf{226} (1977)
125-150.

\bibitem[BHPV]{BHPV} Wolf P. Barth, Klaus Hulek, Chris A. M. Peters, Antonius Van de Ven, \emph{Compact complex surfaces}, 2nd ed., Springer (2004). 

\bibitem[Brg]{Baraglia} David Baraglia, D.Phil. thesis, University of Oxford (2009).

\bibitem[Ber]{Berger} Marcel Berger, \emph{Sur les groupes d`holonomie des vari%
\'{e}t\'{e}s \`{a} connexion affines et des vari\'{e}t\'{e}s riemanniennes},
B. Soc. Math. Fr. \textbf{83} (1955) 279-330.

\bibitem[Bre]{Brezis} Haim Br\'{e}zis, \emph{Analyse fonctionnelle: th\'{e}%
orie et applications}, Dunod, Paris (1999).

\bibitem[Bry]{Bryant} Robert L. Bryant, \emph{Metrics with holonomy }$G_{2}$%
\emph{\ or }$\func{Spin}\left( 7\right) $, Lect. Notes Math. \textbf{1111%
}, Springer-Verlag, Berlin (1985) 269-277.

\bibitem[But]{Buttler} Michael Buttler, D.Phil. thesis, University of Oxford
(1999).

\bibitem[D-K]{4-manifolds} Simon K. Donaldson \& Peter B. Kronheimer, \emph{%
The geometry of four-manifolds}, Oxford Sci. Pub. (1990).

\bibitem[D-T]{Donaldson-Thomas} Simon K. Donaldson \& Richard Thomas, \emph{%
Gauge theories in higher dimensions}, in ``The Geometric Universe; Science,
Geometry, And The Work Of Roger Penrose", Oxford U. P. (1998).


\bibitem[Don$_{1}$]{ASDYM} Simon K. Donaldson, \emph{Anti self-dual
Yang-Mills connections over complex algebraic surfaces and stable bundles},
P. Lond. Math. Soc. \textbf{50 }(1985) 1-26.

\bibitem[Don$_{2}$]{inf dets stab bdls...} Simon K. Donaldson, \emph{%
Infinite determinants, stable bundles and curvature}, Duke Math. J. \textbf{%
54 }(1987) 231-247.

\bibitem[Don$_{3}$]{Approx instantons} Simon K. Donaldson, \emph{%
Approximation of instantons}, Geom. funct. anal. \textbf{3},
n.2 (1993).

\bibitem[Don$_{4}$]{floer} Simon K. Donaldson, \emph{Floer homology groups
in Yang-Mills theory}, Cambridge U. P. (2002).



\bibitem[E-S]{Eels and Sampson} James Eells Jr. \& J. H. Sampson, \emph{%
Harmonic mappings of Riemannian manifolds}, Am. J. Math. 
\textbf{86}, n.1 (1964) 109-160.

\bibitem[Fe-G]{Fernandez&Gray} Marisa Fern\'andez \& Alfred Gray, \emph{Riemannian manifolds with structure group $G_2$}, Ann. Mat. Pur. Appl. \textbf{32} (1982) 19-45.

\bibitem[Fe-U]{Fernandez&Ugarte} Marisa  Fern\'andez \& Luis Ugarte, \emph{Dolbeault cohomology for $G_2-$manifolds}, Geometriae Dedicata \textbf{70} (1998) 57-86.

\bibitem[Fin]{Joel} Joel Fine, PhD Thesis, University of London, (2004).

\bibitem[Fr-U]{freed-uhl} Daniel S. Freed \& Karen K. Uhlenbeck, \emph{%
Instantons and four-manifolds}, MSRI Publications (1984).

\bibitem[G-H]{Griffiths & Harris} Phillip Griffiths \& Joseph Harris, \emph{%
Principles of algebraic geometry}, Wiley-Interscience (1994).

\bibitem[G-H-J]{Calabi-Yau} Mark Gross, Daniel Huybrechts \& Dominic Joyce, 
\emph{Calabi-Yau manifolds and related geometries}, Springer, Berlin (2003).

\bibitem[G-T]{Gilb-Trud}David Gilbarg \& Neil S. Trudinger, \emph{Elliptic partial differential equations of second order}, Springer (2001).
\bibitem[Gri]{Heat kernel bounds} Alexander Grigor'yan, \emph{Gaussian upper
bounds for the heat kernel on arbitrary manifolds}, J. Differ. Geom. 
\textbf{45 }(1997) 33-52.

\bibitem[Guo]{Guo} Guang-Yuan Guo, \emph{Yang-Mills fields on cylindrical
manifolds and holomorphic bundles I}, Commun. Math. Phys. \textbf{179 }%
(1996) 737-775.

\bibitem[Ham]{Hamilton} Richard Hamilton, \emph{Harmonic Maps of Manifolds
with Boundary}, Lect. Notes  Math. \textbf{471}, Springer-Verlag
(1975).

\bibitem[Har]{Hartshorne} Robin J. Hartshorne, \emph{Algebraic Geometry},
Springer-Verlag (2006).

\bibitem[Huy]{Huybrechts} Daniel Huybrechts, \emph{Complex Geometry},
Springer (2005).

\bibitem[J-S]{FAE} Marcos B. Jardim \& Henrique N. S\'a Earp, \emph{Monad constructions of asymptotically stable bunldes}, in preparation.

\bibitem[Jar]{Jardim Bull Braz} Marcos B. Jardim, \emph{Stable bundles on $3-$fold hypersurfaces}, Bull. Braz. Math. Soc., New series \textbf{38(4)} (2007) 649-659. 


\bibitem[Joy]{Joyce} Dominic D. Joyce, \emph{Compact manifolds with special
holonomy}, Oxford Sci. Pub. (2000).

\bibitem[Kov$_{1}$]{kovalevzinho} Alexei Kovalev, \emph{From Fano threefolds
to compact }$G_{2}-$\emph{manifolds}, in \emph{Strings and Geometry}, Clay
Mathematics Proceedings, vol. 3 (2003).

\bibitem[Kov$_{2}$]{kovalevzao} Alexei Kovalev, \emph{Twisted connected sums
and special Riemannian holonomy}, J. Reine Angew. Math. \textbf{565} (2003)
125-160.

\bibitem[M-S]{Milnor} John W. Milnor \& James D. Stasheff, \emph{%
Characteristic classes}, Princeton U. P. (1974).

\bibitem[Mos]{Moser} J\"{u}rgen Moser, \emph{On Harnack's theorem for
elliptic differential equations}, Communications on Pure and Applied
Mathematics \textbf{14} (1961) 577-591.

\bibitem[Muk$_1$]{MukaiFano3-folds}
S. Mukai, \emph{Fano 3-folds},
In: \emph{Complex Projective Geometry}, LMS Note Series \textbf{179}, Cambridge University Press (1992) 255--263.

\bibitem[Muk$_2$]{MukaiBundlesOfK3}
S. Mukai, \emph{On the moduli space of bundles on K3 surfaces I}, Proceedings of the Bombay Conference 1984, Tata Institute of Fundamental Research Studies \textbf{11}, Oxford University Press (1987) 341-413, 

\bibitem[O-S-S]{Okonek-Spindler-Schneider} Christian Okonek, Michael
Schneider \& Heinz Spindler, \emph{Vector bundles on complex projective
spaces}, Progress in Mathematics \textbf{3 }(1979).

\bibitem[Rey]{Reyez-Carrion} Ram\'on Reyes-Carri\'on, \emph{A generalization of the notion of instanton}, Differ. Geom.  Appl. \textbf{8} (1998) 1-20.

\bibitem[Rud]{Rudin} Walter Rudin, \emph{Principles of mathematical analysis}%
, 3rd Ed. McGraw-Hill (1976).

\bibitem[SaE$_0$]{MyThesis} Henrique N. S\'a Earp, PhD Thesis, Imperial College London (2009).

\bibitem[SaE$_1'$]{G2I} Henrique N. S\'a Earp, \emph{$G_2-$instantons over Kovalev manifolds}, arXiv: 1101.0880.

\bibitem[SaE$_2$]{G2II} Henrique N. S\'a Earp, \emph{$G_2-$instantons over Kovalev manifolds II}, in preparation.

\bibitem[Sal]{Salamon} Simon Salamon, \emph{Riemannian geometry and holonomy
groups}, Pitman Res. Notes Math. \textbf{201}, Longman, Harlow (1989).

\bibitem[Sho]{Shokurov} Vyacheslav V. Shokurov, \emph{Smoothness of a general anticanonical divisor on a Fano variety}, Izv. Akad. Nauk SSSR Ser. Mat. \textbf{43} (1979) 395-405.

\bibitem[Sim]{Simpson} Carlos T. Simpson, \emph{Constructing variations of
Hodge structure using Yang-Mills theory and applications to uniformization},
J. Am. Math. Soc. \textbf{1}, n.4 (1988) 867-918.

\bibitem[Tau]{Taubes} Clifford H. Taubes, \emph{Metrics, connections and gluing theorems}, AMS Reg. Conf. series in Math. \textbf{89} (1996).

\bibitem[Tho]{Thomas} Richard Thomas, D.Phil. Thesis, Univeristy of Oxford
(1997).

\bibitem[Tia]{Tian} Gang Tian, \emph{Gauge theory and calibrated geometry I}, arxiv:math.DG/0010015v1.

\bibitem[Wal]{Walpuski} Thomas Walpuski, \emph{$G_2-$instantons on generalised Kummer constructions. I}, arXiv: 1109.6609v2.

\end{thebibliography}
\end{document}